\documentclass[a4paper,12pt,twoside]{article}
\usepackage[latin1]{inputenc}                 
\usepackage{amsfonts,amssymb,amsmath,exscale} 
\usepackage[dvips]{graphicx,color,psfrag}     
\usepackage{fancyhdr}                         
\usepackage{caption}                          
\usepackage{rotating}
\usepackage{defcml}                           
\usepackage[numbers]{natbib}                  
\Floatboxname{Box}
\usepackage{verbatim}                         

\usepackage{amsmath,amsthm,amsfonts,amssymb,amscd}
\usepackage{subfigure}
\usepackage{graphicx,bm,color}
\usepackage{algorithm}
\usepackage{algpseudocode}
\usepackage{stmaryrd}
\usepackage{booktabs,ctable,multirow,longtable} 
\usepackage[latin1]{inputenc}                 
\usepackage{amsfonts,amssymb,amsmath,exscale} 
\usepackage[dvips]{graphicx,color,psfrag}     
\usepackage{fancyhdr}                         
\usepackage{caption}                          
\usepackage{rotating}
\usepackage{defcml}                            
\usepackage[numbers]{natbib}                  
\Floatboxname{Box}

\renewcommand{\qedsymbol}{\rule{0.7em}{0.7em}}

\definecolor{gray}{gray}{0.6}

\newtheorem{Remark}{Remark}[section]
\newtheorem{theorem}{Proposition}

\newtheorem{form}{Formulation}[section]

\def\R{\mathbb R}

\fancypagestyle{plain}{%
	\fancyhead{}%
	\fancyfoot[c]{\sffamily\thepage}%
}
\makeatletter
\def\cleardoublepage{\clearpage\if@twoside \ifodd\c@page\else
	\hbox{}
	\vspace*{\fill}
	\thispagestyle{empty}
	\newpage
	\if@twocolumn\hbox{}\newpage\fi\fi\fi}

\usepackage{scalerel,stackengine}
\stackMath
\newcommand\reallywidecheck[1]{%
	\savestack{\tmpbox}{\stretchto{%
			\scaleto{%
				\scalerel*[\widthof{\ensuremath{#1}}]{\kern-.6pt\bigwedge\kern-.6pt}%
				{\rule[-\textheight/2]{1ex}{\textheight}}
			}{\textheight}%
		}{0.5ex}}%
	\stackon[1pt]{#1}{\scalebox{-1}{\tmpbox}}%
}

\begin{document}
\Titel{
Adaptive Global-Local Approach for Phase-Field Modeling of Anisotropic Brittle Fracture
      }
\Autor{N. Noii, F. Aldakheel, T. Wick, P. Wriggers}
\Report{02--I--17}
\Journal{

}
%



\thispagestyle{empty}

\ce{\bf\large An Adaptive Global-Local Approach for Phase-Field Modeling}
\vspace*{0.25cm}
\ce{\bf\large of Anisotropic Brittle Fracture}

\vskip .35in

\ce{Nima Noii\(^a\), Fadi Aldakheel\(^{b}\)\footnote{Corresponding author. Phone +49\,511\,762\,4126 \, ; \, Fax +49\,511\,762\,5496\\[1mm]
		E-mail addresses:
		noii@ifam.uni-hannover.de (N. Noii); aldakheel@ikm.uni-hannover.de (F. Aldakheel); thomas.wick@ifam.uni-hannover.de (T. Wick); wriggers@ikm.uni-hannover.de (P. Wriggers).
	}, Thomas Wick\(^{a,c}\), Peter Wriggers\(^{b,c}\)} \vskip .25in
\ce{\(^a\) Institute of Applied Mathematics} \ce{Leibniz Universit\"at Hannover, Welfengarten 1, 30167 Hannover, Germany} \vskip .22in

\ce{\(^b\) Institute of Continuum Mechanics} \ce{Leibniz Universit\"at Hannover, Appelstrasse 11, 30167 Hannover, Germany}\vskip .25in

\ce{\(^c\) Cluster of Excellence PhoenixD (Photonics, Optics, and
	Engineering - Innovation} \ce{Across Disciplines) Leibniz Universit\"at Hannover, Germany}\vskip .25in

\begin{Abstract}
	This work addresses an efficient Global-Local approach supplemented
	with predictor-corrector adaptivity applied to anisotropic phase-field brittle
	fracture. The phase-field formulation is used to resolve the sharp
	crack surface topology on the anisotropic/non-uniform local state in
	the regularized concept. To resolve the crack phase-field by a given
	single preferred direction, second-order structural
	tensors {are} imposed to {both the} bulk and crack surface
	density {functions.} 
 Accordingly, a split in tension and compression modes in anisotropic materials is considered.
	A Global-Local formulation is proposed, in which the full
	displacement/phase-field problem is solved on a {\it lower} (local)
	scale, while dealing with a purely linear elastic problem on an {\it
	upper} (global) scale. Robin-type boundary conditions are
	introduced 
to {relax} the stiff local response at the global scale and
	enhancing 
 {its stabilization.}
{Another important aspect of this contribution is the development of an} adaptive Global-Local approach, where a predictor-corrector scheme is {designed} in which the local domains are {\it dynamically updated} during the computation. To cope with different finite element discretizations {at} the interface between the two nested scales, a {\it non-matching} dual mortar method is formulated. {Hence, more regularity is achieved on the interface.} Several numerical results substantiate our developments.   
	
\end{Abstract} 

\sectpa[Section1]{Introduction}

Heterogeneous materials such as wood, composites and bones are
composed of complicated constituents on different scales. Most of these anisotropic materials, 
even with similar constituents properties at the upper scale, 
can behave differently {on} the lower scale. Such heterogeneous 
responses of solid materials 
{are} related to non-uniform and anisotropic behavior {on} the lower scale. 

The multi-scale family can be classified in two distinct classes denoted as
{\it hierarchical} and {\it concurrent} multi-scale techniques. {These}
are defined by differentiation of the global characteristic length scale
$\frak{L}_{global}$ with its local domain counterpart $\frak{L}_{local}$. In
the hierarchical multi-scale method, the average size of the heterogeneous
local domain is much smaller than its global specimen size,
i.e. $\frak{L}_{local}\ll\frak{L}_{global}$, see \cite{Michel, Fish2014}. This
is often denoted as {scale separation law}, see computational homogenization
approaches based on the Hill-Mandel principal; outlined for
instance in \cite{Hill, Michel} among others. {On the other hand}, the concurrent multi-scale method implies $\frak{L}_{local}\equiv\frak{L}_{global}$, as classified in \cite{Lloberas,Fish2014}. Herein, the local periodicity (which is the underlying assumption of 
classical computational homogenization) is not applicable. Then, the
full resolution of the non-linear response on the local scale must be taken
into account, due to the strain localization effect, 
as outlined in \cite{Fish3}. {These} type of materials require a {different} multi-scale
framework {in which the non-linear response is consistently projected
to the global scale; see for example} \cite{Miehe5,Markovic,Zhang,Hautefeuille}.

In the present contribution, we {develop} a multi-scale
approach \cite{Hughes98,Gendre2,Loehnert07,Allix09,NoiiGL18} when the characteristic length of the local scale is of the same order as its global counterpart. 
This is accomplished by introducing a Global-Local approach based on the idea
of a {\textit{history-dependent algorithm}} at the nodal level,
see \cite{Miehe5} and {references cited} therein. This algorithm refers to the procedure in which the boundary value problem of one scale is solved based on the given information from another scale (as a history variable). Accordingly, the history-dependent algorithm contains both the {\textit{upscaling}} and {\textit{downscaling}} steps. In the upscaling step, the global response is achieved, whereas the lower scale information is retained, representing a {\textit{local-global-transition}} procedure. 
However in the downscaling step, a {\textit{re-localization/re-meshing}} of the coarse domain is performed at the local level, see \cite{allix1, Chevreuil}, thereafter solving a non-linear boundary value problem, based on the information passed from the global scale, representing {\textit{global-local-transition}} procedure.

In this work, the Global-Local approach is employed as a computational
framework 
for solving fracture mechanics problems {as it was first formulated in \cite{NoiiGL18}.} 
{Therein, the following assumptions were made} 
\cite{Fish2014, Lozinski09}: (i) The nonlinear phenomenological constitutive
law 
(e.g. {the} failure mechanism) is embedded {on} the local scale and
linear behavior is assumed {on} the global scale. (ii) The global level is
free from geometrical imperfections and hence heterogeneities exist {\it only}
on the local level. (iii) {On} the local level, we consider a divergence-free
assumption for the stress state, such that it is free from any external
imposed load. Accordingly, an interface energy functional based
on the Localized Lagrange Multiplier method \cite{Park2000,Park2002}
is desired for the coupling {of different domains and scales}. 

{Global-Local approaches easily allow for different spatial discretizations} for the global and local
domains. {This} enables computations {and couplings} with legacy codes for industrial
applications in more efficient settings. 
{In this regard, a flexible choice of {the} discretization scheme can be employed on
each domain independently}; 
e.g. {the} Finite Element Method (FEM) \cite{wriggers06}, Isogeometric
Analysis (IGA) \cite{hughes05} 
and {the} Virtual Element Method (VEM) \cite{Wriggers+aldakheel+blaz19}.
A typical application using a simplified Glocal-Local
model {was} done in \cite{WiSiWhe15}. Therein,
a (phase-field) fracture model (computed with deal.II \cite{BangerthHartmannKanschat2007} in C++) 
was employed as local problem {using finite elements}. {The local
setting} was then coupled to a reservoir simulator (IPARS \cite{IPARS} based
on Fortran) for computing the global
problem. For {this global problem, different discretization schemes, the mainly 
based on finite differences for subsurface fluid flow, were adopted.} 

{In the following, we describe in more detail our main goals.}
{First, we focus on the development} of {a} new Global-Local
formulation 
based on Robin-type boundary conditions 
{\cite{Gendre2, Magoules2006, Magoules2007}.} 
{These conditions} relax the stiff local response {transferred} 
to the global scale and thus enhance the stabilization of the Global-Local
approach. {We briefly recall that Robin-type boundary conditions 
contain both Dirichlet and Neumann conditions.  The formulation
is based on an optimized Schwarz method} in {a} multiplicative manner 
{(see for instance \cite{Magoules2006} or \cite{Magoules2007})}. 

{The second goal of} this contribution is to use the Global-Local scheme 
for the analysis of anisotropic fracture processes. {Specifically}, the 
{\it continuum phase-field approach to brittle fracture} 
is
{employed \cite{FraMar98,BourFraMar00,Bour07,HAKIM2009,amor+marigo+maurini09,miehe+welschinger+hofacker10a,KUHN10}. }
Due to its 
{capability} of capturing complex crack patterns in
various engineering applications, this methodology has attracted a
considerable attention in recent years. 
Using such a variational approach, discontinuities in the displacement field are approximated 
across the lower-dimensional crack surface by an auxiliary phase-field function.
The latter can be viewed as an indicator function, which
introduces a diffusive transition zone between the broken and the unbroken material. 
The essential aspects of a phase-field fracture propagation 
formulation are techniques that must include resolution of the length-scale 
parameter with respect to spatial discretization, efficient and robust
numerical solution procedures,
and the enforcement of the irreversibility of crack growth.
%
%
Recent {studies} on phase-field modeling of isotropic brittle fracture
have been devoted to the multiplicative decomposition of the deformation
gradient into compressive-tensile parts in \cite{hesch+weinberg14}, 
coupled thermo-mechanical and multi-physics problems \cite{miehe+schaenzel+ulmer15},
dynamic cases in \cite{borden+verhoosel+scott+hughes+landis12}, 
a new fast hybrid formulation in \cite{ambati15}, 
different choices of degradation functions in \cite{sargado17}, 
and cohesive fracture in \cite{verhoosel+deborst13}. 
Further applications include hydraulic fracture \cite{MiWheWi14,heider17,Miehe2015186}, 
nonlinear solvers \cite{WickLagrange2014,Wick15Adapt},
linear solvers \cite{FaMau17,HeiWi18_pamm},
crack penetration or deflection at an interface in \cite{paggi17} 
and the virtual element method in \cite{aldakheel+blaz+wriggers18}.

{A} considerable number of materials exhibits 
{\it anisotropic behavior}. There span a wide spectrum of applications such as failure in
rocks \cite{DONATH61,NASSERI08}, 
tearing experiments in thin sheets \cite{Takei13} and
biomechanics \cite{holzapfel+gasser+ogden00,balzani+etal06}. 
{Numerical formulations for anisotropic phase-field modeling of brittle
fracture are investigated, for instance, in
\cite{Fadianiso17,NGUYEN17,GULTEKIN16,LI19,bleyer+alessi18}.}
Anisotropic materials exhibit heterogeneous behavior {on} the local domain
through {a} fiber reinforced structure {allowing for}
a homogeneous resolution on the global {level.} Therefore, 
heterogeneous materials often 
require distinct multi-scale treatments such that the full resolution {on} the
local scale {must} be taken into account. 
{In this paper, we therefore propose a}
phase-field approach to brittle fracture in anisotropic solids based on the
{previously described} Global-Local scheme.

Our third main goal is {the adaptive assignment of the local domain(s)
during 
a computation. This is achieved with adaptivity.}
{The adaptive procedure has two goals: (i) to adjust dynamically
the local domain when fractures are propagating; (ii) to reduce the total 
computational cost because the local domains are tailored to the a priori
unknown fracture path. This procedure is much cheaper than using a large local domain
from the beginning.} 
Our approach is inspired by \cite{Wick15Adapt} in which a dynamic 
update in form of a predictor-corrector scheme of crack-oriented mesh refinement was developed. We now apply this idea to the Global-Local approach.
{In the predictor step, mesh edges are identified below
{a} given threshold value for the
phase-field {variable on} the local level.} 
{On} the global level, neighboring elements are subsequently found, then
re-meshed. {Afterwards,} the old solution is interpolated. 
{In the corrector step, we take the old solution and compute the problem
again, but now on the newly determined local domain. Specifically, the
predictor-corrector
approach is now capable to deal with brutal fracture growth; i.e. when
a complete failure} happens in one load increment.

{The} key requirement for realizing {this} adaptive Global-Local
scheme is a non-matching discretization method on the interface. {To this end,}
a dual mortar method \cite{Wohlmuth,REIS2014168} is implemented, thus 
providing sufficient regularity of the underlying {meshes. 
Consequently, different meshes for the global and local domains can be
employed that allow for a very flexible discretization and mesh generation.}

{In a final step, in addition to} 
the local crack phase-field, we determine the coarse representation of the
crack phase-field at the global level. 
{This is a} \textit{post-processing} step and is computed based on
either (\texttt{a}) solving {the} crack phase-field on the global level,
denoted as \textit{global crack phase-field} solution. {Or,}
(\texttt{b}) by 
means of {a} \textit{homogeneous crack phase-field} solution, 
which is an extension {of the} isotropic formulation 
given in \cite{miehe+schaenzel+ulmer15,miehe+hofacker+schaenzel+aldakheel15} to {our proposed} anisotropic phase-field setting.

In summary, this work contains:
\begin{itemize}
	\item {A modular} framework for a phase-field formulation of fracture in anisotropic solids;
	\item {{A} Global-Local approach in order to capture the full local
	resolution at the global level;}
		
	\item Robin-type boundary conditions {between the local and the global domains}; 
	
	\item A non-matching finite element discretization
        {for achieving} sufficient regularity along {the} coupling interface;
	
	\item A  predictor-corrector adaptive scheme {in which} the local domains
	are dynamically updated during the computation;
	
	\item A coarse representation of the crack phase-field at the global level.
\end{itemize}

The paper is structured as follows: In Section \ref{Section2}, we outline the variational anisotropic phase-field formulation of brittle fracture.
Section \ref{Section3} presents the Global-Local approach to capture the local
heterogeneities and constitutive non-linearities at the global level. This is
augmented by introducing a Robin-type 
boundary {conditions.} Then in Section \ref{Section4}, a robust and efficient
predictor-corrector Global-Local adaptive approach is
developed. Section \ref{Section5} contains numerical results that
demonstrate the modeling capabilities of the proposed approach. Qualitative
and quantitative comparisons {with a single scale phase-field solution
are provided, as well}. Finally, the last section concludes the paper with some remarks.
\sectpa[Section2]{Variational Anisotropic Phase-Field Brittle Fracture}

\begin{figure}[!ht]
	\centering
	{\includegraphics[clip,trim=4cm 16cm 8cm 13cm, width=11cm]{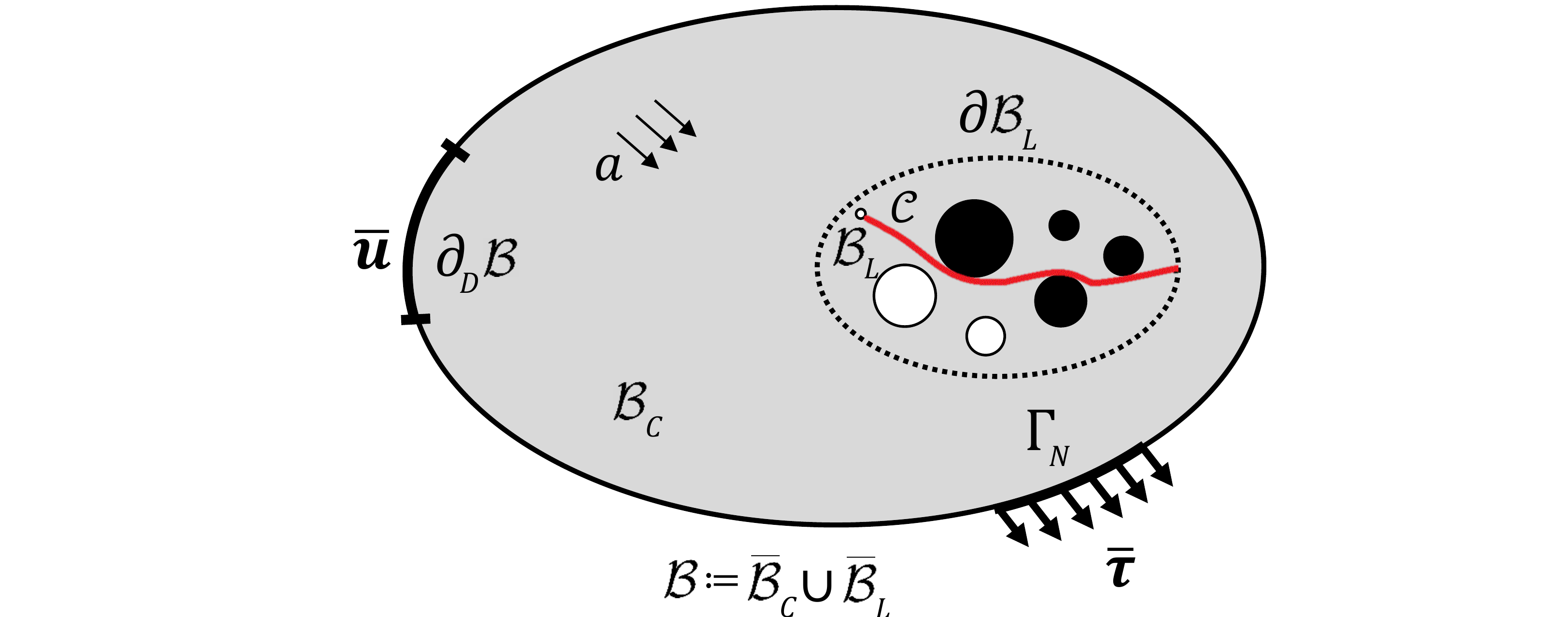}}  
	\caption{ Setup of the notation: the unbroken domain is denoted by
		$\calB_C$ and $\mathcal{C}$ is the crack phase-field. The smeared crack
		phase-field is approximated by the domain $\calB_L$. The whole domain is defined as a close subset as ${\calB:=\bar{\calB}_C\cup\bar{\calB}_L}$. The fracture boundary is $\partial \calB_L$ and the outer boundary of the domain is $\partial
		\calB$.
%
%
	}
	\label{Figure1}
\end{figure}

\sectpb[Section20]{The primary fields of anisotropic brittle solids}
In the following,
let $\calB\subset\R^\delta$, $\delta = 2$ be a smooth open and bounded set with $\partial\calB$
denoted as its boundary. We assume a Dirichlet boundaries conditions
$\partial_D\calB $ and Neumann condition on $\partial_N \calB
:= \Gamma_N \cup \mathcal{C}$, where $\Gamma_N$  denotes the outer domain
boundary and the lower dimensional fracture $\mathcal{C}\in \R^{\delta-1}$ is
the crack boundary, as illustrated in Fig. \ref{Figure1}. Let $I:=(0,t)$
denote the loading/time interval with $t>0$ being the end time value. Using a
phase-field approach, the fracture surface $\mathcal{C}$ is approximated in
$\calB_L\subset\calB \in \mathbb{R}^\delta$ so-called \textit{local
domain}. The intact region with no fracture is denoted
as \textit{complementary domain} 
$\calB_C:=\calB \backslash \calB_L\subset\calB \in \mathbb{R}^\delta$, such
that $\bar{\calB}_C\cup\bar{\calB}_L=:\calB$ 
and $\bar{\calB_C}\cap\bar{\calB_L}=\varnothing$. We note that $\calB_L$,
i.e. the  
domain in which the smeared crack phase-field is approximated, and its
boundary $\partial \calB_L$ 
depend on the choice of the phase-field regularization parameter $l>0$. 
This fracture length scale parameter $l$ is related to the discretization of a domain.  This means in
particular that {$h=o(l)$} (see e.g., \cite{Bou99} for the related
problem of image segmentation) where $h$ denotes the usual spatial
discretization parameter. A simplified numerical analysis on {$h=o(l)$} is 
provided in \cite{MaWi19}.
A detailed computational analysis was performed
in \cite{Wi16_dwr_pff,HeiWi18_pamm}. Moreover, the loading
interval {$\calT := (t_0,T)$} is discretized using the discrete time (loading) points
	\[
	0 = t_0 < t_1 < \ldots < t_n < \ldots < t_N = T,
	\]
{with the end time value $T>0$.} {The parameter $t\in \calT$ denotes for rate-dependent problems the time, for rate-independent problems an incremental loading parameter.}

A phase-field approach to fracture leads to a multi-field problem that depends on the \textit{deformation field} and the \textit{crack phase-field}
\eb
\Bu : 
\left\{
\begin{array}{ll}
	\calB \times \calT \rightarrow \R^\delta \\
	(\Bx, t)  \mapsto \Bu(\Bx,t)
\end{array}
\right.
\AND
d: 
\begin{cases} 
	\begin{array}{l}
		\calB\times \calT \rightarrow [0,1] \\
		(\Bx,t) \mapsto d(\Bx,t),
	\end{array}
\end{cases} 
\label{1_variational_formulation}
\ee
of a material point $\Bx \in \calB$ at time $t\in\calT$. 

Specifically, we deal with a diffusive formulation that interpolates between the intact (unbroken) region with $d=1$ and the fully fractured state of the material with $d=0$ at $\Bx\in\calB$. The Neumann boundary condition $\nabla d.\bn=0$ is imposed on $\partial\calB$ with $\Bn$ being the outward normal to the surface. The strain is assumed to be small, i.e. the norm of the displacement gradient $|| \nabla \Bu || < \epsilon$ is bounded by a small number $\epsilon$. 

%

\sectpb[Section21]{Variational formulation for the multi-field problem}

In this section, we recapitulate a variational approach to brittle fracture in elastic solids at small strains. The energy stored in a bulk strain density for isotropic materials is characterized by the three invariants,
\begin{equation}\label{eq1}
I_1(\bm{\varepsilon})=tr(\bm{\varepsilon}) \quad , \quad I_2(\bm{\varepsilon})=tr(\bm{\varepsilon}^2) \quad , \quad I_3(\bm{\varepsilon})=tr(\bm{\varepsilon}^3).
\end{equation}
Additionally, it is assumed that the solid material is reinforced by only one family of fibers which is denoted as transversely isotropic material. A single preferred
direction at point $\Bx$ is defined by the normal vector $\bm{a}(\bm {x})$ with $\| \Ba \| = 1$ that is called structural director. This type of material has
the highest strength in the direction of the fiber and depicts isotropic response
along its orthogonal direction. Hence, the stress state at {a} material point $\bm {x}$ depends on the deformation and the given single preferred direction. Thus it results to a \textit{deformation-direction-dependent} problem. To do so, a penalty-like parameter $\chi>0$ is defined to restrict a deformation on the normal plane to $\Ba$. The effective bulk free energy now depends on two second-order tensorial quantities, namely the strain $\bm{\varepsilon}$ and structural $\BM$ tensor, defined as
\begin{equation}\label{eq2}
\Bve = \nabla_s \Bu = \sym[ \nabla \Bu ]
\AND
\BM:=\Ba \otimes \Ba,
\end{equation}
they can be represented by additional two deformation-direction-dependent invariants 
\begin{equation}\label{eq3}
I_4(\bm{\varepsilon};\BM)=tr(\bm{\varepsilon.M})=\| \Ba \|_{\bm{\varepsilon}}^{2} = \lambda_{a_0}^2 \quad , \quad I_5(\bm{\varepsilon};\BM)=tr(\bm{\varepsilon}^2.\bm{M}).
\end{equation}
Note, that $I_4$ is nothing else than the quadratic stretch in the direction $\Ba$ of the fiber. Let the effective strain density function, ${\Psi}(\bm{\varepsilon};\BM)$ possesses the property of the transversely isotropic material which has the coordinate-free representation for both matrix and fiber materials. Thus the following holds
\begin{equation}\label{eq4}
{\Psi}(\bm{\varepsilon};\BM)=:{\Psi}(\BQ\bm{\varepsilon}\BQ^T;\BQ\BM\BQ^T) \quad \forall \BQ\in \calG \subset \calO(3),
\end{equation}
that holds for all orthogonal tensor $\BQ$, i.e. $\BQ^T\BQ=\BQ\BQ^T=\BI$, that is a subset of the symmetry group $\calG$ of the anisotropic material. $\BI = \delta_{ij}$ is the second order identity tensor. We denote that ${\Psi}(\bm \varepsilon;\BM)$ is a scalar-valued isotropic tensor function of the symmetric strain tensor $\bm{\varepsilon}$ and the structural tensor $\BM$. Hence, the scalar-valued effective strain density function is an invariant in \textit{space and time} between two pairs of point in the given domain under rotation. Thus ${\Psi}(\bm \varepsilon;\BM)$ can be represented by the principal invariants of $\bm{\varepsilon}$ and $\BM$ as
\begin{equation}\label{eq5}
\begin{aligned}
{\Psi}(\bm\varepsilon;\BM)&={\Psi}\big(I_1(\bm{\varepsilon}),I_2(\bm{\varepsilon}),I_4(\bm{\varepsilon};\BM),I_5(\bm{\varepsilon};\BM)\big)\\
&=\widetilde{\Psi}^{iso}\big(I_1(\bm{\varepsilon}),I_2(\bm{\varepsilon})\big)+\widetilde{\Psi}^{aniso}\big(I_4(\bm{\varepsilon};\BM),I_5(\bm{\varepsilon};\BM)\big).
\end{aligned}
\end{equation}
Herein, the isotropic  free-energy function corresponds to
\begin{equation}\label{eq6}
\widetilde{\Psi}^{iso}\big(I_1(\bm{\varepsilon}),I_2(\bm{\varepsilon})\big):=\frac{\lambda}{2}I^2_1+\mu I_2,
\end{equation}
with $\lambda >0$ and $\mu>0$ being the elastic Lam\'e constants. The anisotropic free-energy function is defined as 
\begin{equation}\label{eq7}
\widetilde{\Psi}^{aniso}\big(I_4(\bm{\varepsilon};\BM),I_5(\bm{\varepsilon};\BM)\big):=\frac{1}{2}\chi I^2_4+2 \,{\Xi}\, I_5,
\end{equation}
with the anisotropic material parameters $\chi$ and $\Xi$. A stress-free condition, i.e. $\bm{\varepsilon}=0$, is required $\widetilde{\Psi}^{iso}\big(I_1(\bm{0}),I_2(\bm{0})\big)=0$. Moreover, $\widetilde{\Psi}^{aniso}\big(I_4(\bm{0};\BM),I_5(\bm{0};\BM)\big)=0$ must hold true.

Using these definitions, to establish variational based anisotropic phase-field approach to brittle fractures,  we define the bulk free energy functional which represents the stored energy in bulk as
\begin{equation}
\begin{aligned}
\mathcal{E}_{bulk}(\Bu;\BM)=\int_{\calB_C}{\Psi}(\bm\varepsilon;\BM)
\mathrm{d}{\bm{x}}-\int_{{\partial_N\calB_C }} {\bm {\bar\tau}} \cdot \bm u\,\mathrm{d}s .
\label{eq8}
\end{aligned}
\end{equation}
Herein, {$\bm{\bar\tau}$} denotes the traction forces on 
the complementary boundaries ${\partial_N\calB_C := \Gamma_N \cup \mathcal{C}}$. 

The total energetic functional is based on both the stored bulk energy as well as the fracture dissipation, defined in the  work of \cite{FraMar98},
\begin{equation}
\mathcal{E}(\Bu,\mathcal{C};\BM)= \mathcal{E}_{bulk}(\Bu;\BM)
+ G_c \mathcal{H}^{\delta-1}  (\mathcal{C}) \, ,
\label{eq9}
\end{equation}
where $G_c$ is the {Griffith's critical elastic energy release rate}
and $\mathcal{H}^{\delta-1}$ is a $\delta-1$ dimensional 
Hausdorff measure. For the numerical treatment we regularize Eq.
\ref{eq9} following \cite{BourFraMar00}. Specifically, 
the crack energy is approximated through a sequence of elliptic 
problems, so-called Ambrosio-Tortorelli functionals, see \cite{AmTo90,AmTo92}.
Therein, $\mathcal{H}^{\delta-1}$ is regularized by the crack phase-field $d$.
Finally, we account for the crack irreversibility constraint meaning the crack can only 
grow:
\begin{equation}\label{eq10}
\dot{d} \leq 0.
\end{equation}
In the incremental version, this condition reads:
\[
d \leq d^{old},
\]
where $d:= d(t_n)$ and $d^{old} := d(t_{n-1})$. For stating the variational formulations, we 
now introduce:
\begin{equation}\label{space1}
\begin{aligned}
\bm{V}&:= \{ {\bf H}^1(\calB)^\delta:\bm u=\bar{\bm u}\; \mathrm{on} \; \partial_D\calB  \}, \quad W:= \text{H}^1(\calB) ,\\
W_{in} &:= \{ d \in \text{H}^1(\calB)^{\delta-1} | \; 0 \leq
d\leq d^{old} \}.
\end{aligned}
\end{equation}
As typical in problems with inequality constraints (see e.g.,
\cite{KiOd88,KiStam00}), $W_{in}$ is a nonempty, closed, convex,
subset of the linear function space $W$. Due to the inequality constraint in Eq. \ref{eq10}, $W_{in}$ is
no longer a linear space.

\sectpb[Section22]{Phase-field approximation of anisotropic crack topologies}

The variational approach of \cite{BourFraMar08} is widely used for fracture failure phenomena in isotropic elastic solids. As a point of departure, in line with \cite{MieWelHof10b,dittmann+fadi+etal18}, let a regularized macro crack topology of a sharp crack be represented by the exponential function $1-\exp^{[ - \vert x \vert /l ]}$ satisfying $d(0)=0$. We define a regularized \textit{isotropic crack surface energy functional} of the solid by, 
\begin{equation}\label{eq11}
G_c \mathcal{H}^{\delta-1}_{iso}  (\mathcal{C}):=G_c\int_\calB\gamma^{iso}_l(d, \nabla d) \,\mathrm{d}{\Bx}
\WITH
\gamma^{iso}_{l}(d,\nabla d):=  \frac{1}{2l} {(1-d)^2} + \frac{l}{2}  \nabla d.\nabla d,
\end{equation}
in terms of the isotropic crack surface density function per unit volume of the solid $\gamma^{iso}_{l}$.


The above representation of a crack surface density function is extended for the class of anisotropic responses; as for instance outlined in \cite{NGUYEN17,Fadianiso17,na+sun18}. Similarly to the deformation field, we define a total crack surface density function. It is additively decomposed into an isotropic and anisotropic crack surface density function, respectively, as follows
\begin{equation}\label{eq13}
\gamma_l(d, \BQ\nabla d) := \gamma^{iso}_l(d, \nabla d)+\gamma^{aniso}_l(\nabla d)
\quad \mbox{for all} \quad
\BQ\in \calG \subset \calO(3),
\end{equation}
where $\calG$ is a given symmetry group of the anisotropic material, i.e. the set of rotation and reflection vectors, and $\gamma^{aniso}_l(d, \nabla d)$ represents an augmented crack surface density for the anisotropic response. Let this function posses  the  property  of  the  transversely  isotropic  material  which  has  the  coordinate-free representation  for  both  matrix and fibers materials. This leads to
\begin{equation}\label{eq14}
\gamma_l(d, \nabla d; \BM) =  \gamma_l(d^2, \nabla d \cdot \nabla d, 
\nabla d \cdot \BM \cdot \nabla d)\;.
\end{equation}
Using this definition, an anisotropic crack surface density response can be defined as
\begin{equation}\label{eq15}
\gamma^{aniso}_{l}(\nabla d;\BM):=  \alpha\frac{l}{2}  \nabla d.\BM.\nabla d.
\end{equation}
This type of fracture function has the highest geometric resistance in
the {fiber direction} and has an {isotropic response along its orthogonal direction}.  Hence, the geometric resistance state at a material point $\bm {x}$ depends on crack phase-field and the given single preferred direction $\bm {a}$. This results in a \textit{crack-direction-dependent} problem.

The anisotropic term $\alpha$ in Eq. \ref{eq15} behaves as a penalty-like parameter and hence for $\alpha \rightarrow \infty$ one obtains $\nabla d \cdot \Ba = 0$, which means that the crack lies parallel to the preferred orientation. 
	For $\alpha\rightarrow0$ the isotropic response will be recovered. 

\begin{form}[Energy functional for the anisotropic crack topology]
	\label{form_1}
	Let $\BM$, $\chi$, $\Xi$ and $\alpha$ be given with the initial conditions $\bm u_0=\bm u(\bm{x},0)$ and  $d_0=d(\bm{x},0)$. For the loading increments 
	$n=1,2,\ldots, N$, find $\bm u:=\bm u^n\in V$ and $d:=d^n\in W_{in}$ such that the functional
	\begin{align*}
	\mathcal{E} ( {\bm u},d;\BM)&:=\mathcal{E}_{bulk}(\bm u,d_+;\BM
	)+\; \mathcal{E}_{frac}(d;\BM)+\;\mathcal{E}_{ext}(\bm u)\\[1mm] 
	&\:=\int_{\calB} w_{bulk}(\bm{\varepsilon},d;\BM) \mathrm{d}{\bm{x}} + \int_\calB w_{frac}(d, \nabla d;\BM) \, \mathrm{d}{\bm{x}}           
	-\int_{\partial_N\calB} {\bm {\bar\tau}} \cdot \bm u\,\mathrm{d}s,
	\end{align*}
	is minimized. The elastic bulk density $w_{bulk}$ along with the fracture contribution $w_{frac}$ both define the so-called total pseudo-energy density function as
	\begin{align*}
&w(\Bve, d, \nabla d;\BM) = w_{bulk} (\Bve, d;\BM) +  w_{frac} (d, \nabla d;\BM), \\
&w_{bulk}(\Bve, d;\BM) = g(d_+){\Psi}(\bm\varepsilon;\BM), \\
&w_{frac}(d, \nabla d;\BM) = G_c \gamma_l(d, \nabla d; \BM).
\end{align*}
\end{form}

\begin{Remark}
	\label{Rem_phi_pos}
	In the case of elastic cracks, it can be shown that the phase field
	satisfies  $0 \leq d \leq 1$. When additional physics are included for
	instance a fluid inside the fracture \cite{MiWheWi18} or non-isothermal effects
	\cite{NoiiWick2019}, the energy functional must be modified to cope with
	negative values of $d$. Hence {in order to allow for future extensions}, we work in the remainder of this paper with $d_+$ rather than $d$.
	A detailed discussion is provided in \cite{MiWheWi18}[Section 3].  
\end{Remark}

\begin{Remark}
	\label{effect_of_phf}
	The comparison of the bulk energy functional in Eq. \ref{eq9}
	and Formulation \ref{form_1} is tow-fold. First, the integration is
	changed from $\calB_C$ to the entire domain $\calB$ due to the presence of
	the phase-field function $d$. 
Second, the presence of $d$ in {the} bulk energy {through the degradation function
	$g(d)$} 
defines the transition state from the unbroken to fracture 
state hence results in the degradation of the solid material as well as the crack propagation.
\end{Remark}

\sectpb[Section23]{Strain-energy decomposition}

Since the fracturing material behaves quite differently in \textit{tension}
and {\textit{compression}}, a consistent split for the strain energy
density function is employed, where we apply the decomposition \textit{only} to the isotropic strain energy function, i.e. $ {\Psi}^{iso}\big(I_1(\bm{\varepsilon}),I_2(\bm{\varepsilon})\big)$. Hence, instead of dealing directly with $\bm\varepsilon(\bm u)$, we perform
additive decomposition of the strain tensor as
\[
\bm\varepsilon(\bm u)=\bm\varepsilon^{+}(\bm u)+\bm\varepsilon^{-}(\bm
u)
\WITH
\bm\varepsilon^{\pm}(\bm
u):=\sum_{i=1}^{\delta} \langle\varepsilon_i\rangle^{\pm} {\textbf{N}_i} \otimes {\textbf{N}_i}
\ ,
\]
with the tension $\bm\varepsilon^{+}$ and compression $\bm\varepsilon^{-}$ strains. Here, $\langle x \rangle_{\pm} := \frac{ x {\pm} |x|}{2}$ is a ramp function of $\R_{\pm}$ expressed by the Macauley bracket. {$\{\varepsilon_i\}$ are the principal strains and $\{\textbf{N}_i\}$ are the principal strain directions.} The tension/compression fourth-order projection tensor is defined as 
\begin{equation}\label{eq16}
\mathbb{P}^\pm_{\bm {\varepsilon}}:=\frac{\partial \bm {\varepsilon}^\pm}{\partial \bm {\varepsilon}}=\frac{\partial \big(\sum_{i=1}^{\delta} \langle\varepsilon_i\rangle^{\pm}  {\textbf{N}_i} \otimes {\textbf{N}_i}\big)}{\partial \bm {\varepsilon}} \ .
\end{equation}
It turns out that, $\mathbb{P}^\pm_{\bm {\varepsilon}}$ projects the total strain into its positive and negative parts accordingly, i.e $\bm {\varepsilon}^{\pm}=\mathbb{P}^\pm_{\bm {\varepsilon}}:\bm {\varepsilon}$. So, a decoupled representation of the strain-energy function into a so-called tension and compression contribution is given as follows,
\begin{equation}\label{eq17}
{\Psi}^{iso}\big(I_1(\bm{\varepsilon}),I_2(\bm{\varepsilon})\big):=\underbrace{\widetilde{\Psi}^{iso,+}\big(I^{+}_1(\bm{\varepsilon}),I^{+}_2(\bm{\varepsilon})\big)}_{\text{tension term}}+\underbrace{\widetilde{\Psi}^{iso,-}\big(I^{-}_1(\bm{\varepsilon}),I^{-}_2(\bm{\varepsilon})\big)}_{\text{compression term}}.
\end{equation}
Herein, the positive and negative principal invariants are
\begin{equation}\label{eq18}
I_1^{\pm}(\bm{\varepsilon}):=\langle{I_1(\bm{\varepsilon})}\rangle_{\pm}, \quad I^{\pm}_2(\bm{\varepsilon}):=I_2(\bm{\varepsilon}^{\pm}).
\end{equation}
\begin{Remark}
	\label{alternative_I1}
	An alternative definition to $I_1^{\pm}(\bm{\varepsilon}):=\langle{I_1(\bm{\varepsilon})}\rangle_{\pm}$, can be defined by using the same description introduced in $I_2^{\pm}(\bm{\varepsilon})$ for the first principal invariant which results in $I^{\pm}_1(\bm{\varepsilon}):=I_1(\bm{\varepsilon}^{\pm})$. This provides a new description for the strain-energy function represented in Eq. \ref{eq17}. However, that is beyond the scope of present paper and will investigated in future work.
\end{Remark}
Physically, it is trivial to assume that the degradation induced by the phase field acts only on the tensile and shear counterpart of the elastic strain density function. Hence, it is that assumed there is no degradation in compression, which also prevents interpenetration of the crack lips during crack closure, see \cite{MieWelHof10b}. It turns out that the bulk work density function for the fracturing material becomes,
\begin{equation}\label{eq19}
w_{bulk}(\Bve, d;\BM):=g(d_+)\Big[\widetilde{\Psi}^{iso,+}(I^{+}_1,I^{+}_2)
+ \widetilde{\Psi}^{aniso} (I_4,I_5) \Big]
+\widetilde{\Psi}^{iso,-}(I^{-}_1,I^{-}_2)
.
\end{equation}
Here a monotonically decreasing quadrature degradation function, i.e. 
\begin{equation}
g(d_+):=(1-\kappa)d_+^2 + \kappa,
\label{eq20}
\end{equation}
describes the degradation of the solid with the evolving crack phase-field
parameter $d$. The small residual stiffness $\kappa$ is introduced to prevent numerical problems. The constitutive stress response corresponding to Eq. \ref{eq19} reads
\begin{equation}
\begin{aligned}
&{\bm \sigma}(\bm{\varepsilon},d ;\BM):=\frac{\partial w_{bulk}(\Bve, d;\BM)}{\partial {\bm \varepsilon}} = \bm {\sigma}_{\bm \varepsilon}^{iso} + \bm {\sigma}_{\bm \varepsilon}^{aniso},
\\[0.5mm]
&\bm {\sigma}_{\bm \varepsilon}^{iso}
=g(d_+)\frac{\partial \widetilde{\Psi}^{iso,+}}{\partial \bm \varepsilon}+\frac{\partial \widetilde{\Psi}^{iso,-}}{\partial \bm \varepsilon}
 =g(d_+){\bm {\widetilde{\sigma}}^{iso,+}_{\bm \varepsilon}}+{\bm {\widetilde{\sigma}}^{iso,-}_{\bm \varepsilon}},
 \\[0.5mm]
&\bm {\sigma}_{\bm \varepsilon}^{aniso}= g(d_+)\frac{\partial \widetilde{\Psi}^{aniso}}{\partial \bm \varepsilon} = g(d_+){\bm {\widetilde{\sigma}}^{aniso}_{\bm \varepsilon}},
\label{eq21}
\end{aligned}
\end{equation}
with,
\begin{equation}
\bm {\widetilde{\sigma}}_{\bm \varepsilon}^{iso,\pm}:=\lambda I_1^{\pm}(\bm{\varepsilon})  {\textbf{I}} 
+ 2\mu \bm\varepsilon_\pm
\AND
{\bm {\widetilde{\sigma}}^{aniso}_{\bm \varepsilon}}=
\chi I_4 \BM + 2\,\Xi\,(\Bve\cdot\BM + \BM\cdot\Bve)
\label{eq24}
\end{equation}
%
%
\begin{form}[Energy functional for the anisotropic crack topology]
	\label{form_2}
	Let $\BM$, $\chi$, $\Xi$ and $\alpha$ be given with the initial conditions $\bm u_0=\bm u(\bm{x},0)$ and  $d_0=d(\bm{x},0)$. For the loading increments 
	$n=1,2,\ldots, N$, find $\bm u:=\bm u^n\in V$ and $d:=d^n\in W_{in}$ such that the functional
	\begin{align*}
	\mathcal{E} ( {\bm u},d;\BM)&=\underbrace{\int_\calB \big[g(d_+) \widetilde{\Psi}^{iso,+} +\widetilde{\Psi}^{iso,-}\big] \mathrm{d}{\textbf{x}}}_{\text{matrix deformation term}}
    +\underbrace{\int_\calB g(d_+) \widetilde{\Psi}^{aniso}   	\mathrm{d}{\textbf{x}}}_{\text{fiber deformtion term}}\\
	&+ G_c\underbrace{\int_\calB\gamma^{iso}_{l}\mathrm{d}{\textbf{x}}}_{\text{matrix fracture term}}
	+ G_c\underbrace{\int_\calB\gamma^{aniso}_{l}\mathrm{d}{\textbf{x}}}_{\text{fiber fracture term}}-\underbrace{\int_{\partial_N\calB } {\bm {\bar\tau}} \cdot \bm u\,\mathrm{d}s}_{\text{external load}},
	\end{align*}\\[0.7mm]
{is minimized.}
\end{form}
The minimization problem for the given energy functional of the anisotropic crack topology in Formulation \ref{form_2} takes the following compact form:
\begin{equation}\label{compat_argmin}
\fterm{ 
	\{ \bm{u}, d \} = 
	\mbox{arg} \{\; 
	\substackrel{\bm{u} \in \bf{V} }{\mbox{min}}
	\substackrel{d \in W_{in} }{\mbox{min}} \,
	[\; \mathcal{E} (\bm u,d;\BM) \; ]
	\; \}.
	}
\end{equation}
The stationary points of the energy functional in Formulation \ref{form_2} are characterized
by the first-order necessary conditions, namely the so-called Euler-Lagrange equations, which are obtained by {differentiation with respect to ${\bm u}$ and $d$.}
\begin{form}[Euler-Lagrange equations]
	\label{form_3}
		Let $\BM$, $\chi$, $\Xi$ and $\alpha$ be given with the initial conditions $\bm u_0=\bm u(\bm{x},0)$ and  $d_0=d(\bm{x},0)$. For the loading increments $n=1,2,\ldots, N$, find $\bm u:=\bm u^n\in V$ and $d:=d^n\in W_{in}:$	
    \begin{equation}
    \begin{aligned}
    &{\mathcal E}_{\bm u}(\bm u, d;\delta\bm u)=
    \int_\calB g(d_+) {\bm {\widetilde{\sigma}}^{iso,+}_{\bm
    	\varepsilon}}(\bm u): {\bm \varepsilon}(\delta\bm u)\mathrm{d}{\textbf{x}}+ 
    	\int_\calB {\bm {\widetilde{\sigma}}^{iso,-}_{\bm \varepsilon}(\bm u)}: {\bm \varepsilon}(\delta\bm u)\mathrm{d}{\textbf{x}}\\
	 &\qquad \qquad \; \; +\int_\calB g(d_+) {\bm {\widetilde{\sigma}}^{aniso}_{\bm
	 		\varepsilon}}(\bm u):{\bm \varepsilon}(\delta\bm u)\mathrm{d}{\textbf{x}}-\int_{\partial_N\calB } {\bm {\bar\tau}} \cdot {\delta\bm u}\,\mathrm{d}s= 0 \quad\forall {\delta\bm u}\in V, \\
	&{\mathcal E}_d(\bm u,d;\delta d-d) =(1-\kappa)\int_\calB d_+ \big[{\bm {\widetilde{\sigma}}^{iso,+}_{\bm \varepsilon}}(\bm u) + {\bm {\widetilde{\sigma}}^{aniso}_{\bm \varepsilon}}(\bm u) 
	\big]
	: {\bm \varepsilon}(\bm u). (\delta d-d)\mathrm{d}{\textbf{x}}\\
	&\qquad \qquad  \qquad  \qquad+ G_c \int_\calB \left( \frac{1}{l}(d-1).(\delta d-d)
	+ l \nabla d.\nabla(\delta d-d)\right)\mathrm{d}{\textbf{x}}\\
	&\qquad \qquad  \qquad  \qquad+ G_c \int_\calB \alpha l \nabla d.\BM.\nabla(\delta d-d) \,\mathrm{d}{\Bx} \geq 0 \quad\forall
	\delta d\in W \cap L^{\infty}.
    \end{aligned}
    \end{equation}
\end{form}
${\mathcal E}_{\bm u}$ and ${\mathcal E}_d$ are the directional derivatives of the energy functional with respect to $\bm u$ and $d$, respectively. Furthermore, $\delta\bm u\in\{ {\bf H}^1(\calB)^\delta: \delta\bm u=\bm 0 \; \mathrm{on} \; \partial_D\calB \}$ is the deformation test function and $\delta d\in H^1(\calB)$ is the phase-field test function.

\sectpb[Section24]{The Euler-Lagrange equations in a strong form}
In order to complete our derivations, the strong form of Formulation
\ref{form_3} will be derived in this section. Using integration by parts, we
obtain a quasi-stationary elliptic system for the displacements and the
phase-field variable, where the latter one is subject to an inequality
constraint in time and therefore needs to be complemented with a complementary
condition:
\begin{form}[Strong form of the Euler-Lagrange equations]
	\label{form_4}
	Let $\BM$, $\chi$, $\Xi$ and $\alpha$ be given with the initial conditions $\bm u_0=\bm u(\bm{x},0)$ and  $d_0=d(\bm{x},0)$. For the loading increments $n=1,2,\ldots, N$, we solve a displacement equation where we seek $\bm u:= \bm u^n: \calB \rightarrow \mathbb{R}^{d}$ such that
	\begin{align*}
	-\div (\Bsigma )= \Bzero \quad &in \; \calB,\\
	{\bm u} =\bar \Bu \quad &on \;  \partial_D\calB,\\
	\Bsigma \cdot \Bn = \bar{\Btau}\quad &on \;  \partial_N\calB,
	\end{align*}
in terms of the stress tensor $\Bsigma$ defined in Eq. \ref{eq21} and the given displacement field $\bar \Bu$. The phase-field system consists of four parts: the PDE, the
	inequality constraint and a compatibility condition (in fracture
	mechanics called Rice condition \cite{ricecond}) along with the Neumann-type boundary conditions.
	Find $d:=d^n : \calB \rightarrow [0,1]$ such that
	\begin{align*}
	-\Big( 2(1-\kappa)d_+ \big[\widetilde{\Psi}^{iso,+}+
	\widetilde{\Psi}^{aniso}
	\big]
	-\frac{G_c}{l}(1-d)
	-G_c l\Delta d
	-G_c \alpha l \div(\nabla d.\BM)\Big)\ge \; 0 \quad &in \; \calB,\\
	\dot{d}\leqslant 0 \quad &in \; \calB,\\
	-\Big( 2(1-\kappa)d_+ \widetilde{\Psi}^{iso,+}(\bm\varepsilon(\bm u))
    -\frac{G_c}{l}(1-d)
	-G_c l\Delta d-G_c \alpha l \div(\nabla d.\BM)\Big) \; \dot{d}= 0\quad &in \; \calB,\\
	\; ({\textbf{I}}+\alpha\BM) {\nabla d \cdot \Bn} = 0 \quad &on \; \partial\calB.
	\end{align*}
\end{form}
The mentioned inequality minimization problem for the phase-field equation can
be resolved through:
(a) fixing the fracture with \textit{Dirichlet conditions} \cite{BourFraMar00},
(b) \textit{Penalty method}, see \cite{MiWheWi15b} (including a mathematical analysis), 
(c) \textit{an Augumented Lagragian penalization} see \cite{WickLagrange2014}, 
(d) \textit{Primal-dual active set method}; see \cite{Wick15Adapt,LeeWheWi16}, 
(e) \textit{Maximum crack driving state function}, see \cite{MieWelHof10b,miehe+schaenzel+ulmer15}. 
In the present work, we consider the maximum crack driving state function to 
prevent the crack healing by having a positive crack dissipation 
known as irreversibility criteria and given in details in next section. 

\sectpb[Section25]{Crack driving force}
In this section, a formulation for the crack phase-field PDE equation in Formulation \ref{form_4} is reformulated based on the crack driving force. A thermodynamical consistency for the preservation {of the energy balance} due to the fracture dissipation results to the Karush-Kuhn-Tucker form, see \cite{miehe+hofacker+schaenzel+aldakheel15}. 
%
As a point of departure, the modular structure of the anisotropic phase-field fracture equation assumes the following form
\begin{equation}\label{eq28}
	\underbrace{\eta \dot{d}}_{{crack \, \,  update}} = \underbrace{-g^{\prime}(d_+)
		\widetilde{D}}_{crack \, \, driving \, \, force} - \underbrace{l \delta_d
		\gamma_l}_{resistance} \leq 0,
\end{equation}
as outlined in the works of {\sc Miehe} and coworkers \cite{miehe+schaenzel+ulmer15,miehe+hofacker+schaenzel+aldakheel15,aldakheel16,aldakheel+wriggers+miehe18}. Here, $\widetilde{D}$ is a {\it{crack driving state function}} which depends on a $state$ array of strain- or stress like quantities. To get rid of the above inequality evolution problem, we maximize the inequality equation given in (\ref{eq28}) for the full process
history $s\in [0,t_n]$,
\begin{equation}\label{eq281}
        -g^{\prime}(d_+)\lmax{s\in [0,t_n]} \widetilde{D}= l \delta_d
		\gamma_l\;.
\end{equation}
We introduce  maximum positive crack driving force $\calH$ in $t\in [0,t_n]$  denoted as,
\begin{equation}
\calH(\Bx,t):= \lmax{{s}\in [0,t_n]}\widetilde{D}\big(state(\Bx,{\bm{\varepsilon}}(s))\big),
\label{7_irreversible_response}
\end{equation}
and hence (\ref{eq28}) can be restated as,
\begin{equation}\label{eq282}
\fterm{
\underbrace{\eta \dot{d}}_{{crack \, \,  update}} = \underbrace{-g^{\prime}(d_+)
	\calH}_{max \, \, crack \, \, driving \, \, force} - \underbrace{l \delta_d.
	\gamma_l}_{resistance}.
}
\end{equation}
Depending on the type of the crack driving state function which can be either without or with threshold, $\widetilde{D}$ can take different description, see \cite{aldakheel+mauthe+miehe14,aldakheel+blaz+wriggers18,aldakheeletal18}.
The crack phase-field evolution in (\ref{eq282}) is defined in the domain $\calB$ that is augmented with an imposed Neumann homogeneous boundary condition as
\begin{equation}\label{eq30}
{({\textbf{I}}+\alpha\BM) {\nabla d \cdot \Bn} = 0} \quad on \; \partial\calB.
\end{equation}
Note, that in Eq. \ref{eq28} the rate-independent case is
recovered for $\eta \rightarrow 0$, where the crack topology is then simply determined by an
equilibrium between the crack driving force and the geometric crack
resistance. These equations are interpreted as generalized Ginzburg-Landau-type
evolution equations for the crack phase-field $d$. Equation \ref{eq282} restated for the rate-independent limit $\eta \rightarrow 0$ to the so-called Karush-Kuhn-Tucker form: 
\begin{equation}\label{eq31}
\fterm{
\dot d \leqslant 0 ; \quad - \delta_{d} w \; \ge \;0 ,\quad \dot d\, [- \delta_{d} w ] = 0.
}
\end{equation}
This condition provides a natural assumption due to the positive fracture dissipation know as crack irreversibility condition. The latter constraint is ensured by a specific constitutive assumption that relates the functional derivative to a positive energetic driving force. The last condition in (\ref{eq31}) is the balance law for the evolution of the crack phase-field which ensures the principal of maximum dissipation during the crack phase-field
evolution (see e.g. \cite{MieWelHof10b}). It is known as compatibility condition.
\begin{Remark}
	\label{KKT}
    Karush-Kuhn-Tucker stated in (\ref{eq31}) {along with (\ref{eq30})} are the modular structure of the Euler-Lagrange equations in a strong form {which are given in Formulation \ref{form_4}}.
\end{Remark}
By defining the maximum positive crack driving force $\calH$ in terms of the the crack driving state function $\widetilde{D}$ at hand, Formulation \ref{form_3} can be stated as an equality minimization. Thus $\calH$ substitutes the corresponding $\big[\widetilde{\Psi}^{iso,+}+
\widetilde{\Psi}^{aniso}
\big]$ term in the original $\mathcal{E}_d$.
%
To derive the crack driving state function we recall the irreversibly inequality condition, i.e. $\dot d \le 0$. It follows that the left hand side of (\ref{eq28}) has to be always positive to avoid the crack healing process
\begin{equation}\label{29_H1}
- \delta_{d} w=(\kappa-1)2d_+ \big[\widetilde{\Psi}^{iso,+}+
\widetilde{\Psi}^{aniso}
\big] - G_c \delta_d\gamma_l(d, \nabla d;\BM)
\leqslant 0.
\end{equation}
Maximization of this inequality in the full process history $s\in [0,t_n]$, yields
\begin{equation}\label{29_H2}
(\kappa-1)2d_+\lmax{s \in [0,t_n]}
\big[\widetilde{\Psi}^{iso,+}+
\widetilde{\Psi}^{aniso}
\big] =G_c \delta_d\gamma_l(d, \nabla d;\BM).
\end{equation} 
To follow the modular structure of the phase-field fracture equation defined in Eq. \ref{eq28}, we multiply  (\ref{29_H2}) by $\frac{l}{G_c}$. With the definition of a positive crack driving force and shown by $\calH$, hence (\ref{29_H2}) is restated as
\begin{equation}\label{eq29_H3}
(\kappa-1)2d_+\calH =l \delta_d\gamma_l \; \;
\mbox{if} \; \calH:= \lmax{s\in [0,t_n]}\widetilde{D} \; \; \mbox{with}\; \; \widetilde{D}:=\frac{l\big[\widetilde{\Psi}^{iso,+}+
	\widetilde{\Psi}^{aniso}
	\big]}{G_c}.
\end{equation}
It is evident that the crack driving state function given by (\ref{eq29_H3}) is directly effected by the regularization parameter $l$. Hence the crack driving state function has the property of length-scale dependency. The functional derivative of $\gamma_l$ with respect to $d$ is obtained as follows,
\begin{equation}\label{eq29_H4}
\int_\calB \delta_{d} \gamma_{l}(d,\nabla d;\BM) \mathrm{d}{\bm{x}}:=\int_\calB \bigg( \frac{\partial \gamma_l}{\partial d} +\frac{\partial \gamma_l}{\partial \nabla d}\bigg) \mathrm{d}{\bm{x}} =
\int_\calB \bigg( \frac{\partial \gamma_l}{\partial d}-\nabla.[\frac{\partial \gamma_l}{\partial \nabla d}]\bigg) \mathrm{d}{\bm{x}}, 
\end{equation}
which leads to,
\begin{equation}\label{eq29_H5}
\int_\calB \delta_{d} \gamma_{l}(d,\nabla d;\BM) \mathrm{d}{\bm{x}}=\int_\calB \frac{1}{l}[(d-1)-l^2 \Delta d-\alpha l^2 \div(\nabla d.\BM)] \mathrm{d}{\bm{x}}.
\end{equation}
Furthermore $\frac{\partial \gamma_{l}(d,\nabla d;\BM)}{\partial 
	{\nabla d}}\cdot\bm n=({\textbf{I}}+\alpha\BM)l \nabla d \cdot\bm n=0$ refers to (\ref{eq30}). 
\\

\begin{form}[Final Euler-Lagrange equations]
	\label{form_3_final}
	Let $\BM$, $\chi$, $\Xi$ and $\alpha$ be given with the initial conditions $\bm u_0=\bm u(\bm{x},0)$ and  $d_0=d(\bm{x},0)$. For the loading increments $n=1,2,\ldots, N$, find $\bm u:=\bm u^n\in V$ and $d:=d^n\in W:$	
	\begin{equation}
	\begin{aligned}
	&{\mathcal E}_{\bm u}(\bm u, d_+;\delta\bm u)=
	\int_\calB g(d_+) {\bm {\widetilde{\sigma}}^{iso,+}_{\bm
			\varepsilon}}(\bm u): {\bm \varepsilon}(\delta\bm u)\mathrm{d}{\textbf{x}}+ 
	\int_\calB {\bm {\widetilde{\sigma}}^{iso,-}_{\bm \varepsilon}(\bm u)}: {\bm \varepsilon}(\delta\bm u)\mathrm{d}{\textbf{x}}\\
	&\qquad \qquad \; \; +\int_\calB g(d_+) {\bm {\widetilde{\sigma}}^{aniso}_{\bm
			\varepsilon}}(\bm u):{\bm \varepsilon}(\delta\bm u)\mathrm{d}{\textbf{x}}-\int_{\partial_N\calB } {\bm {\bar\tau}} \cdot {\delta\bm u}\,\mathrm{d}s= 0 \quad\forall {\delta\bm u}\in V, \\
	&{\mathcal E}_d(\bm u,d;\delta d) =(1-\kappa)\int_\calB 2d_+ \calH \delta d \mathrm{d}{\textbf{x}}\\
	&\qquad \qquad  \qquad + \int_\calB \Big((d-1) \delta d
	+ l^2 \nabla d.\nabla\delta d\Big)\mathrm{d}{\textbf{x}}\\
	&\qquad \qquad  \qquad +  \int_\calB \alpha l^2 \nabla d.\BM.\nabla\delta d  \,\mathrm{d}{\Bx} = 0 \quad\forall
	\delta d\in W.
	\end{aligned}
	\end{equation}
\end{form}
\sectpa[Section3]{Global-Local Formulation Applied to the Anisotropic Phase-Field
	Fracture}

Departing point towards {\em a Global-Local approach} applied to the anisotropic phase-field formulation is the domain decomposition method \cite{DDRey06}. We split the single-scale energy functional $\mathcal{E}$ indicated in Formulation \ref{form_2} to the intact and fractured region, i.e. $\calB_C$ and $\calB_L$, respectively.
%

Accordingly, by introduction of the Fictitious domain $\calB_F$, i.e. a coarse
projection of the local domain into the global domain ({later $\calB_G$
refers to the global} domain in Section \ref{Section32}), we extend the resulting non-overlapping domain decomposition formulation toward a Global-Local formulation applied to the anisotropic phase-field fracture. The Global-Local formulation applied to isotropic phase-field was first proposed by \citet{NoiiGL18}.
The main objective was to introduce an adoption of the phase-field formulation within legacy codes, specifically for industrial applications.

An important definition for the subsequent treatment is the energy
functional. {We recall that the} energy functional for the single-scale problem denoted as $\mathcal{E}$. We further define the energy functional for the domain decomposition by $\widehat{\mathcal{E}}$ and the Global-Local formulations as $\widetilde{\mathcal{E}}$.
\sectpb[Section31]{Non-overlapping domain decomposition formulation}

\begin{figure}[!ht]
	\centering
	{\includegraphics[clip,trim=1cm 14cm 1cm 5cm, width=14cm]{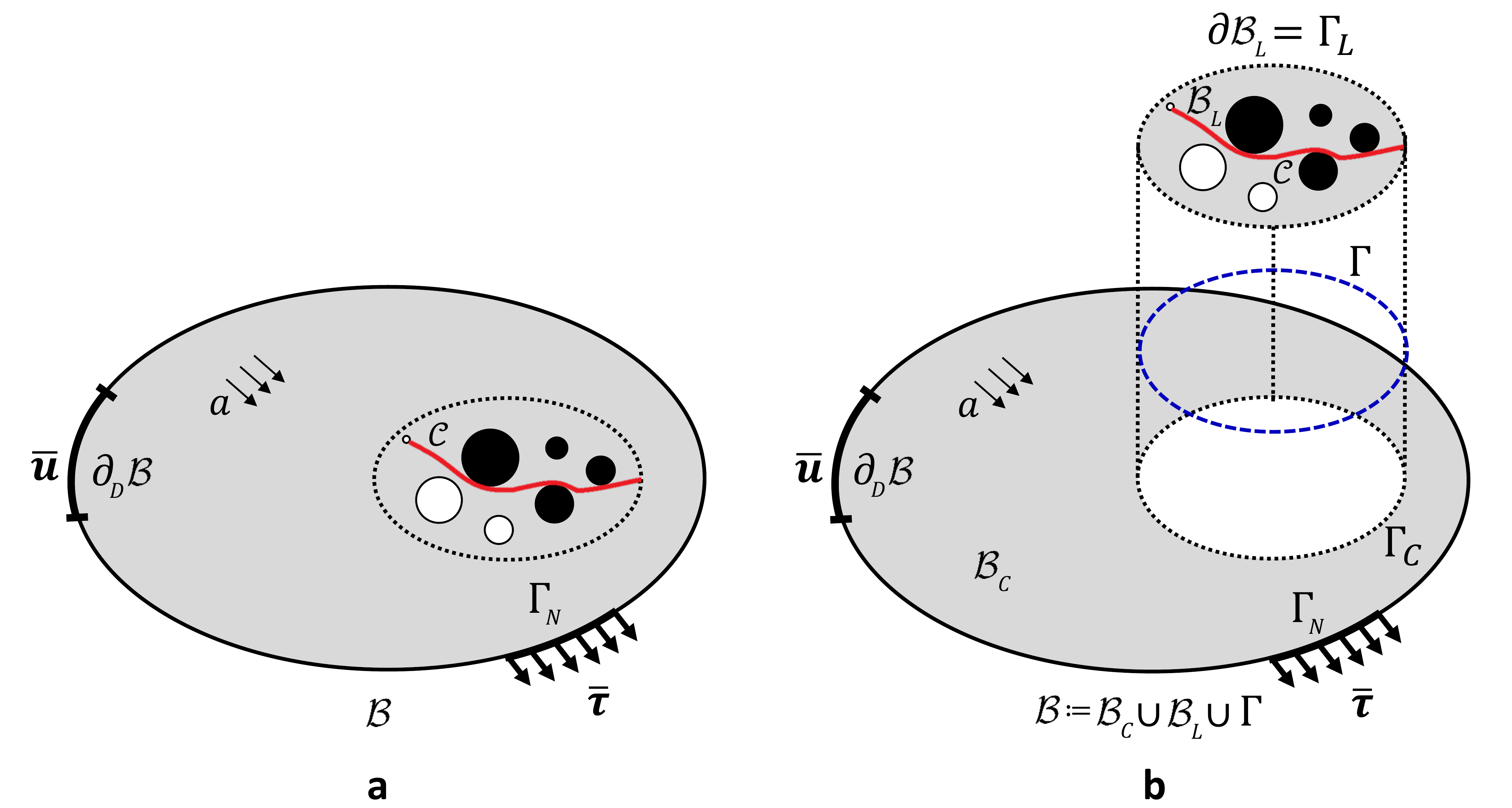}}  
	\caption{Domain decomposition scheme. $(a)$ Geometry and loading setup of the single-scale boundary value problem. $(b)$ Non-overlapping domain decomposition setting whereas $\calB$ is decomposed into the intact and fracture region denoted as complementary and local domains $\calB_L$ and $\calB_C$, respectively.
	}
	\label{Fig2}
\end{figure}

Recall, the complementary domain $\calB_C:=\calB \backslash \calB_L\subset\calB \in \mathbb{R}^d$ corresponds to the intact region and let $\calB_L$ is an open domain, where the fracture surface is approximated in this region, see Fig. \ref{Fig2}(a). It is assumed the fracture surface in $\calB_L$ represents a reasonably small 'fraction' of $\calB$ such that $|\calB_L|\ll|\calB_C|$. We further define an interface between an unfractured domain $\calB_C$ and fractured domain $\calB_L$ by $\Gamma\in \R^{\delta-1}\subset\calB$ in the continuum setting to be the interface between $\calB_L$ and $\calB_C$, such that $\calB\equiv \calB_L\cup\Gamma\cup\calB_C$. We further assume that $\calB_L$ is free from any externally imposed load 
and hence we have prescribed loads {\it only} in $\calB_C$. Such an assumption is standard for the multi-scale setting, see \cite{Fish2014}.

Consider a domain decomposition with geometric sketch in Fig. \ref{Fig2}(b) applied to the single-scale domain plotted in Fig. \ref{Fig2}(a). Two functions on $\calB_L$ and $\calB_C$ are considered, namely $\bm u_L\in{\bf{V}_L}$ and $\bm u_C\in{\bf{V}_C}$, where we introduce additional three sets:
\[
\bm{V}_L := \textbf{H}^1_0 (\calB_L)^\delta, \quad {\bm{V}_C:= {\bf{V}}(\calB_C)}, \quad \text{and} \; \; W_{L}:= W_{in}(\calB_L) \; ,
\]
{referring} to the spaces defined in Eq. \ref{space1}. 

A descriptive motivation of the domain decomposition approach applied to the variational anisotropic phase-field modeling is related to two restriction in the model: (i) the strong coupling scheme that is the strong displacement continuity condition that holds along with (ii) the predefined interface. To this end, one needs to assume that the discrete interfaces for both complementary and local domain do exactly coincide in the strong sense, yielding
\begin{equation}
\bm u_L\overset{!}{=}\bm u_C \quad  \mbox{at} \quad {\bm x}\in\Gamma=\Gamma_C=\Gamma_L .
\label{uCont}
\end{equation}
This displacement continuity is often called {\em primal} approach in the literature, see e.g.\ \cite{Mandel1993}.

Let the {\it single-scale displacement field} $\bm u\in{\bf V}$ be the solution of the multi-field variational problem in (\ref{compat_argmin}). It is decomposed as
\begin{equation}
\bm u(\bm x,t) :=\left\{
\begin{tabular}{ll}
$\bm u_L$ & for ${\Bx}\in\calB_L$, \\[0.1cm]
$\bm u_C$ & for ${\Bx}\in\calB_C$.
\end{tabular}
\right.
\label{uLC}
\end{equation}
Since the fracture surface lives only in $\calB_L$ we introduce scalar-valued function $d_L:\calB_L\rightarrow[0,1]\in W_{L}$. The {\it single-scale phase-field} $d$ is then decomposed in the following form
\begin{equation}
d(\bm x,t)  :=\left\{
\begin{tabular}{ll}
$d_L$ & for ${\Bx}\in\calB_L$, \\[0.1cm]
$1$ & for ${\Bx}\in\calB_C$.
\end{tabular}
\right.
\label{dLC}
\end{equation}
By imposing (\ref{uLC}) and (\ref{dLC}) to the energy functional, indicated in Formulation \ref{form_2}, energy functionals corresponding to $\calB_C$ and $\calB_L$ reads
\begin{equation}
\widehat{\mathcal E}_C(\bm u_C;\BM):
= \int_{\calB_C} w(\Bve_C, 1, 1;\BM) \,\mathrm{d}{\Bx}
- \int_{\partial_N\calB_C} {\bm {\bar\tau}} \cdot \bm u_C \,\mathrm{d}s,
\label{E1}
\end{equation}
and
\begin{equation}
\widehat{\mathcal E}_{L}(\bm u_L,d_L;\BM):
=\int_{\calB_L} w(\Bve_L, d_L, \nabla d_L;\BM) \,\mathrm{d}{\Bx},
\label{E2}
\end{equation}
for the total energy density defined in Formulation \ref{form_1}. With the strong displacement continuity in (\ref{uCont})we obtain
\begin{equation}
{\mathcal E}(\bm u,d;\BM) \equiv \widehat{\mathcal E}(\bm u_C,\bm u_L,d_L;\BM)
:= \widehat{\mathcal E}_C(\bm u_C;\BM)+\widehat{\mathcal E}_{L}(\bm u_L,d_L;\BM),
\label{E1pE2}
\end{equation}
where ${\mathcal E}$ is the original single-scale functional in Formulation \ref{form_2}. As a result, the domain decomposition variational formulation is {\em equivalent} to the single-scale formulation Eq. \ref{compat_argmin}
\begin{equation}\label{DD0}
	\{ \bm u_C, \bm u_L, d_L \} = 
	\mbox{arg} \{\; 
	\substackrel{\bm u_C \in {\bf{V}_C}, \bm u_L \in \bf{V}_L }{\mbox{min}}
	\substackrel{d_L \in W_{L} }{\mbox{min}} \,
	[\; \widehat{\mathcal E} (\bm u_C, \bm u_L, d_L;\BM) \; ]
	\; \}
	\ .
\end{equation}

Note, the major advantage of using this minimization problem instead of the one in (\ref{compat_argmin}) is the reduction of the nonlinearity order of the complementary domain (which is free from the fracture state), and more specifically in small deformation setting that is a linear minimization problem.

\begin{Remark}
	\label{non_match}
	The strong displacement continuity requirement given in Eq. \ref{uCont} is too restrictive from the computational standpoint \cite{Farhat1991}. To resolve the phase field problem, one requires {$h_L \ll h_G$}. However, if we assume $\bm u_L\overset{!}{=}\bm u_C$  on $\Gamma$, this yields $\Gamma_L=\Gamma_C$ in a discretized setting hence $h_L = h_G$ on $\Gamma$ which has the contradiction with {$h_L \ll h_G$}.
\end{Remark}

Following Remark \ref{non_match}, we relax Eq. \ref{uCont} in a weak sense by introducing traction-like terms in the corresponding energy functionals (\ref{E1}) and (\ref{E2}). This results in

\begin{equation}
\widehat{\mathcal E}_C(\bm u_C,\bm\lambda_C;\BM):
= \int_{\calB_C} w(\Bve_C, 1, 1;\BM) \,\mathrm{d}{\Bx}
- \int_{\Gamma_C} \bm\lambda_C\cdot\bm u_C \,\mathrm{d}s
- \int_{\Gamma_{N,C}} {\bm {\bar\tau}} \cdot \bm u_C \,\mathrm{d}s,
\label{E11}
\end{equation}
and
\begin{equation}
\widehat{\mathcal E}_{L}(\bm u_L,d_L,\bm\lambda_L;\BM):
=\int_{\calB_L} w(\Bve_L, d_L, \nabla d_L;\BM) \,\mathrm{d}{\Bx}
- \int_{\Gamma_L} \bm\lambda_L\cdot\bm u_L \,\mathrm{d}s,
\label{E21}
\end{equation}
with $\bm\lambda_C,\bm\lambda_L\in{\bf L}^2(\Gamma)$ being the unknown Lagrange multipliers, which represent traction forces on the interface. The saddle point problem including complementary and local domains assumes the form
\begin{equation*}
\widehat{\mathcal E}(\bm u_C,\bm u_L,d_L,\bm\lambda_L,\bm\lambda_C;\BM)
:= \widehat{\mathcal E}_C(\bm u_C,\bm\lambda_C;\BM)+\widehat{\mathcal E}_L(\bm u_L,d_L,\bm\lambda_L;\BM),
\end{equation*}
which is under-determined, since no relation is yet specified between $\bm u_L$ and $\bm u_C$, nor between $\bm\lambda_L$ and $\bm\lambda_C$. The latter is achieved by introducing the functional
\begin{equation}
\widehat{\mathcal E}_\Gamma(\bm u_\Gamma,\bm\lambda_C,\bm\lambda_L):
= \int_\Gamma \bm u_\Gamma\cdot(\bm\lambda_L+\bm\lambda_C) \,\mathrm{d}s,
\label{E3}
\end{equation}
with $\bm u_\Gamma\in{\bf H}^1(\Gamma)$ representing the (unknown) Lagrange multiplier, which has the dimension of a displacement, called also displacement interface. Summing ${\mathcal E}_C$ and ${\mathcal E}_{L}$ with ${\mathcal E}_\Gamma$, we get
\begin{equation}
\begin{aligned}
\widehat{\mathcal E}(\bm u_C,\bm u_L,&d_L,\bm u_\Gamma,\bm\lambda_C,\bm\lambda_L;\BM):
= \int_{\calB_C} w(\Bve_C, 1, 1;\BM) \,\mathrm{d}{\Bx}
+\int_{\calB_L} w(\Bve_L, d_L, \nabla d_L;\BM) \,\mathrm{d}{\Bx}\\
&+\int_\Gamma \left\{ \bm\lambda_C\cdot(\bm u_\Gamma-\bm u_C) + \bm\lambda_L\cdot(\bm u_\Gamma-\bm u_L) \right\} \mathrm{d}s
- \int_{\Gamma_{N,C}} {\bm {\bar\tau}} \cdot \bm u_C \,\mathrm{d}s.
\label{E11pE21pE3}
\end{aligned}
\end{equation}
Here the introduction of the intermediate displacement $\bm u_\Gamma$ satisfies the weak traction continuity between $\bm\lambda_L$ and $\bm\lambda_C$ along $\Gamma$. This is in addition to the weak displacement continuity between $\bm{u}_L$ and $\bm{u}_C$ across $\Gamma$. Hence, both displacement and traction continuity are imposed implicitly in the weak sense to the energy functional \cite{Park2000}. The coupling interface energy functional used in Eq. \ref{E11pE21pE3} (i.e. third term) is called {\em Localized Lagrange Multipliers}, see e.g. \cite{Park2002,Park2015}. 

The variational formulation of Eq. \ref{E11pE21pE3} is equivalent to the single-scale minimization problem in Eq. \ref{compat_argmin}, such that ${\mathcal E}\approx \widehat{\mathcal E}$, yields

\begin{equation}\label{compact_gl0}
\begin{aligned}
\bm s= 
	\mbox{arg} \{\; 
	\substackrel{\bm u_C  \in {\bm{V}_C}, \bm u_L \in {\bm{V}_L}, \bm u_\Gamma  \in {\bm H}^1(\Gamma), {d_L \in W_{L}} \; \; }{\mbox{min}}
	\substackrel{\bm\lambda_C,\bm\lambda_L \in {\bf L}^2(\Gamma) }{\mbox{max}} \,
	[\; \widehat{\mathcal E} (\bm u_C,\bm u_L,d_L,\bm u_\Gamma,\bm\lambda_C,\bm\lambda_L;\BM) \; ] \},
\end{aligned}
\end{equation}
where $\bm s :=(\bm u_C,\bm u_L,d_L,\bm u_\Gamma,\bm\lambda_C,\bm\lambda_L)$. Accordingly, the displacement field $\bm u$, is decomposed as
\begin{equation}
\bm u =\left\{
\begin{tabular}{ll}
$\bm u_L$ & for ${\Bx}\in\calB_L$, \\[0.1cm]
$\bm u_C$ & for ${\Bx}\in\calB_C$, \\[0.1cm]
$\bm u_\Gamma$ & for ${\Bx}\in\Gamma$,
\end{tabular}
\right.
\label{uLCnew}
\end{equation}
which is based on the solution triple $(\bm u_C,\bm u_L,\bm u_\Gamma)$ as a minimizer of the (\ref{E11pE21pE3}). Note, the representation for $d$ in terms of $d_L$ defined by (\ref{dLC}) remains same.

\sectpb[Section32]{Global-Local formulation}

In this section, the formulation is extended to a Global-Local approach in line with \cite{NoiiGL18}. Specifically in this paper, we extend the Global-Local formulation to the anisotropic crack phase-field which is augmented by Robin-type boundary conditions \cite{Farhat2LM20002, Magoules2006, Magoules2007, Gander2007}. The latter relaxes the stiff local response observed at the global level which is due to the local non-linearity projected to the global level and leads to further reductions of the computational time. Additionally, to have more regularity along the coupling interface, a non-matching finite element discretization is used on the interface. 

\begin{figure}[!ht]
	\centering
	{\includegraphics[clip,trim=1cm 14cm 1cm 5.3cm, width=14cm]{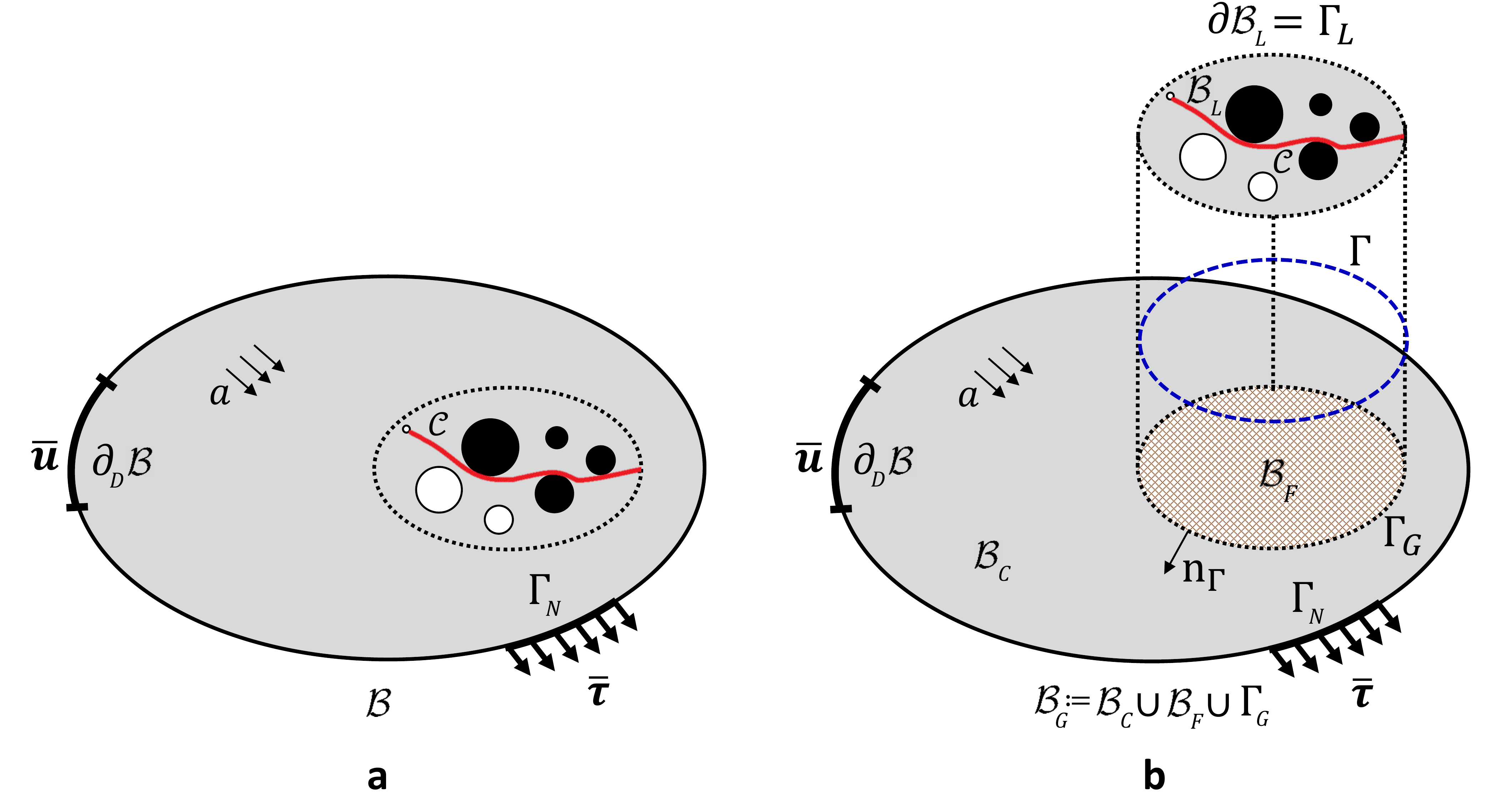}}  
	\caption{Illustration of the Global-Local formulation. $(a)$ Geometry and loading setup of the single-scale boundary value problem. $(b)$ Global-Local setting, by introduction of the fictitious domain $\calB_F$ through prolongation of $\calB_C$ to the entire domain whereas its unification is so-called global domain $\calB_G:=\calB_C\cup\Gamma\cup\calB_F$.
	}
	\label{Fig3}
\end{figure}

Let us define open and bounded \textit{fictitious domain} $\calB_F$ to {\it recover the space} of $\calB$ that is obtained by removing $\calB_L$ from its continuum domain, see Fig. \ref{Fig3}. Indeed, the fictitious domain is prolongation of the $\calB_C$ towards $\calB$. This gives the same constitutive modeling used in $\calB_C$ for $\calB_F$. Thus, the energy functional of the complementary and fictitious domain is the same. We also use the identical discretization space for both $\calB_F$ and $\calB_C$, which results $h_F:=h_C$. We further define, an open and bounded \textit{global domain} $\calB_G$ such that $\calB_G = \calB_F\cup\Gamma\cup\calB_C$. It yields the same energy functional for $\calB_C$, $\calB_F$ and $\calB_G$. Hence, the material parameters are identical for $\calB_C$, $\calB_F$ and $\calB_G$. Additionally, this unification yields on identical discretization space for the global domains $\calB_F$ and $\calB_C$, and results in $h_G\approx h_F \approx h_C$ referring to the element size. 

Note that, the fictitious domain $\calB_F$ is assumed to be free from geometrical {\it imperfections} which may be present in $\calB_L$, see Fig. \ref{Fig3}(b). Thus, the global domain is assumed to be free from any given imperfection. Let us also define, global and local interfaces denoted as $\Gamma_G\subset\calB_G$ and $\Gamma_L\subset\calB_L$, such that in the continuum setting we have $\Gamma=\Gamma_G=\Gamma_L$. However in a discrete setting we might have $\Gamma\neq\Gamma_G\neq\Gamma_L$ due to the presence of different meshing schemes (i.e. different element size/type used in $\calB_G$ and $\calB_L$ such that $h \neq h_L \neq h_G$ on $\Gamma$).

It is assumed that there exists a continuous prolongation of $\bm u_C$ into $\calB_F$. Hence, we introduce a function $\bm u_G\in{\bf V}(\calB_G)$ such that $\bm u_G|_{\calB_C}\equiv\bm u_C$ and $\bm u_G=\bm u_C$ on $\Gamma$ in the sense of a trace. Thus, the boundary conditions for $\calB_G$ is same as the $\calB_C$, therefore it holds $\bm u_G=\bar{\bm u}$ on $\partial_{D} \calB$ and $\bm t=\bar{\bm t}$ on $\Gamma_{N,G}$.
By means of the fictitious domain, the first term in Eq. \ref{E11pE21pE3} is recast as follows
\begin{equation}
\begin{aligned}
\int_{\calB_C} w(\bm\varepsilon(\bm u_C), 1, 1;\BM) \,\mathrm{d}{\Bx}&=
\int_{\calB_C} w(\bm\varepsilon(\bm u_G), 1, 1;\BM) \,\mathrm{d}{\Bx}\\
&=\int_{\calB_G} w(\bm\varepsilon(\bm u_G), 1, 1;\BM) \,\mathrm{d}{\Bx}-\int_{\calB_F} w(\bm\varepsilon(\bm u_G), 1, 1;\BM) \,\mathrm{d}{\Bx}.
\end{aligned}
\label{ugfc}
\end{equation}
Note, we substitute $\bm u_G$ for $\bm u_C$ in the second and fourth integrals in Eq. \ref{ugfc}. That is trivial by means of the prolongation concept such that $\bm u_G|_{\calB_F}\equiv\bm u_F$ and $\bm u_G=\bm u_F$ on $\Gamma$.

This provides the {\em Global-Local} approximation of the single-scale energy functional $\mathcal{E}$ indicated in Formulation  \ref{form_2} by
\begin{equation}
\begin{aligned}
\widetilde{\mathcal E}(\bm u_G,\bm u_L,d_L,\bm u_\Gamma,\bm\lambda_C,\bm\lambda_L;\BM):
&= \underbrace{\int_{\calB_G} w(\Bve_G, 1, 1;\BM) \,\mathrm{d}{\Bx}
- \int_{\calB_F} w(\Bve_F, 1, 1;\BM) \,\mathrm{d}{\Bx}
- \int_{\Gamma_{N,G}} {\bm {\bar\tau}} \cdot \bm u_G\,\mathrm{d}s}_{\text{global terms}} \,\\
&+\underbrace{\int_{\calB_L} w(\Bve_L, d_L, \nabla d_L;\BM) \,\mathrm{d}{\Bx}}_{\text{local term}}\\
&+\underbrace{\int_\Gamma \left\{ \bm\lambda_C\cdot(\bm u_\Gamma-\bm u_G) + \bm\lambda_L\cdot(\bm u_\Gamma-\bm u_L) \right\} \mathrm{d}s}_{\text{coupling terms}}.
\label{EGL}
\end{aligned}
\end{equation}
where the approximation ${\mathcal E} \equiv \widetilde{\mathcal E}$ holds.

\begin{form}[Global-Local energy functional applied to the anisotropic crack topology]
	\label{form_5}
	Let $\BM$, $\chi$, $\Xi$ and $\alpha$ be given with initial conditions $\bm u_0=\bm u(\bm{x},0)$ and  $d_0=d(\bm{x},0)$. For the loading increments 
	$n=1,2,\ldots, N$, find $\bm u_G:=\bm u_G^n\in {\bf V_G}$, $\bm u_L:=\bm u_L^n\in {\bf V_L}$, $\bm u_\Gamma:=\bm u_\Gamma^n\in {\bf H}^1(\Gamma) $, $d_L:=d_L^n\in W_{L}$, $\bm \lambda_C:=\bm \lambda_C^n\in {\bf L}^2(\Gamma) $ and $\bm \lambda_L:=\bm \lambda_L^n\in {\bf L}^2(\Gamma) $, such that the functional
	\begin{align*}
	&\widetilde{\mathcal{E}} ( \bm u_G,\bm u_L,d_L,\bm u_\Gamma,\bm\lambda_C,\bm\lambda_L;\BM)=\\
	&\underbrace{\int_{\calB_G}({\frac{\lambda}{2}I^{2}_1+\mu I_2})\mathrm{d}{\textbf{x}}}_{\text{global matrix deformation term}}+\underbrace{\int_{\calB_G}({\frac{1}{2}\chi I^2_4+2\Xi I_5})\mathrm{d}{\textbf{x}}}_{\text{global fiber deformtion term}}-\underbrace{\int_{\Gamma_{N,G} } {\bm {\bar\tau}} \cdot \bm u_G\,\mathrm{d}s}_{\text{global external load}}\\
	&-\underbrace{\int_{\calB_F}({\frac{\lambda}{2}I^{2}_1+\mu I_2})\mathrm{d}{\textbf{x}}}_{\text{fictitious matrix deformation term}}-\underbrace{\int_{\calB_F}({\frac{1}{2}\chi I^2_4+2\Xi I_5})\mathrm{d}{\textbf{x}}}_{\text{fictitious fiber deformtion term}}\\
	&+\underbrace{\int_{\calB_L}g(d_{L+})\big({\frac{\lambda}{2}I^{+2}_1+\mu I^+_2}\big)\mathrm{d}{\textbf{x}}+\int_{\calB_L}({\frac{\lambda}{2}I^{-2}_1+\mu I^-_2})\mathrm{d}{\textbf{x}}}_{\text{local matrix deformation term}}
	+\underbrace{\int_{\calB_L}g(d_{L+})\big({\frac{1}{2}\chi I^2_4+2\Xi I_5}\big)\mathrm{d}{\textbf{x}}}_{\text{local fiber deformtion term}}\\
	&+G_c\underbrace{\int_{\calB_L}\{\frac{1}{2l} (1-d_L)^2 + \frac{l}{2}  \nabla d_L.\nabla d_L\}\mathrm{d}{\textbf{x}}}_{\text{local matrix fracture term}}
	+G_c\underbrace{\int_{\calB_L}(\alpha\frac{l}{2}  \nabla d_L.\BM.\nabla d_L)\mathrm{d}{\textbf{x}}}_{\text{local fiber fracture term}}\\
	&+\underbrace{\int_\Gamma \left\{ \bm\lambda_C\cdot(\bm u_\Gamma-\bm u_G) + \bm\lambda_L\cdot(\bm u_\Gamma-\bm u_L) \right\} \mathrm{d}s}_{\text{interface coupling term}},
	\end{align*}
	is minimized.
\end{form}
Note, we are not any more using $\partial_N\calB$ for the applied surface load and hence $\Gamma_{N,G}$ is considered. This is because the global domain is free from any fracture state.
The minimization problem for the Global-Local energy functional given in Formulation \ref{form_5} that is applied to the anisotropic crack topology takes the following compact form,
\begin{equation}\label{GL}
\bm s= 
\begin{aligned}
\mbox{arg} \{ 
\substackrel{\bm u_G \in {\bm{V}_G}, \bm u_L \in {\bm{V}_L}, \bm u_\Gamma  \in {\textbf{ H}}^1(\Gamma), {d_L \in W_{L}} \; \; }{\mbox{min}}
\substackrel{\bm\lambda_C,\bm\lambda_L \in {\bf L}^2(\Gamma) }{\mbox{max}} \,
[\; \widetilde{\mathcal{E}}(\bm u_G,\bm u_L,d_L,\bm u_\Gamma,\bm\lambda_C,\bm\lambda_L;\BM) ] \}
\end{aligned}
\end{equation}
where $\bm s :=(\bm u_G,\bm u_L,d_L,\bm u_\Gamma,\bm\lambda_C,\bm\lambda_L)$.

The relation between the solution $\bm u$ of the minimization problem in Eq. \ref{compat_argmin} and the solution triple $(\bm u_G,\bm u_L,\bm u_\Gamma)$ of Eq. \ref{GL} reads
\begin{equation*}
\bm u =\left\{
\begin{tabular}{ll}
$\bm u_L,$ & for ${\Bx}\in\calB_L$, \\[0.1cm]
$\bm u_G,$ & for ${\Bx}\in\calB_G$, \\[0.1cm]
$\bm u_\Gamma,$ & for ${\Bx}\in\Gamma$.
\end{tabular}
\right.
\end{equation*}

\begin{Remark}
	\label{diff_gl_ss}
	When using standard single scale phase-field modeling, we are most of the time not dealing with a	uniform mesh and hence the domain is divided
	into coarser and finer mesh elements. To
	resolve the crack phase-field, we need to have $l \ge h$ {must hold} at every point of the domain such that $l\ge h_c\ge h_f$ ($c$ and $f$ refers to the coarse and fine region in domain, respectively) satisfied. This typically leads to a finer mesh even for the area which is sufficiently far from the fracture zone, and therefore increases the computational time considerably. However, this is not the case for the Global-Local approach where the phase-field formulation is only embedded within the local domain and not the entire domain. Hence the computational time is reduced drastically.    
\end{Remark}

\sectpb[Section33]{Variational formulation for the Global-Local coupling system}

Now we consider the weak formulation of Eq. \ref{GL}. The directional derivatives of the functional $\widetilde{\mathcal E}$ yield for the global weak form 
\begin{align*}
\widetilde{\mathcal E}_{\bm u_G}(\bm s;\delta {\bm u}_G)&:=
\int_{\calB_G}\bm\sigma(\bm u_G;\BM):\bm\varepsilon(\delta {\bm u}_G)\,\mathrm{d}{\Bx}
-\int_{\calB_F}\bm\sigma(\bm u_G;\BM):\bm\varepsilon(\delta {\bm u}_G)\,\mathrm{d}{\Bx}\\
&-\int_{\Gamma_G} \bm\lambda_C \cdot \delta {\bm u}_G\,\mathrm{d}s
-\int_{\Gamma_{N,G} } {\bm {\bar\tau}} \cdot \delta{\bm u}_G\,\mathrm{d}s=0,
\label{Global}
\tag{G}
\end{align*}
where $\bm\sigma(\bm u_G):=\partial_{\bm\varepsilon}w(\bm\varepsilon(\bm u_G),1,1;\BM)$ and $\delta {\bm u}_G\in\{ {\bf H}^1(\calB_G): \delta {\bm u}_G=\bm 0 \; \mathrm{on} \; \partial_D\calB \}$ is the test function. The local weak formulations assumes the form
\begin{equation*}\label{alLocL}
\left\{
\begin{tabular}{l}
$\widetilde{\mathcal E}_{\bm u_L}(\bm s;\delta {\bm u}_L):=\displaystyle \int_{\calB_L}\bm\sigma(\bm u_L,d_L;\BM):\bm\varepsilon(\delta {\bm u}_L)\,\mathrm{d}{\Bx}
-\int_{\Gamma_L} \bm\lambda_L \cdot \delta {\bm u}_L\,\mathrm{d}s=0$, \\[0.1cm]
$\displaystyle \widetilde{\mathcal E}_{d_L}(\bm s;\delta d_L):=(1-\kappa)\int_{\calB_L} 2d_{L+} \calH(\bm\varepsilon(\bm u_L);\BM). \delta d_L\mathrm{d}{\Bx}+  \int_{\calB_L}  (d_L-1).\delta d_L\mathrm{d}{\Bx}$\\[0.1cm]
\qquad \qquad \; \: $ \displaystyle + \int_{\calB_L} l^2 \nabla d_L.\nabla(\delta d_L)\,\mathrm{d}{\Bx}+  \int_{\calB_L} \alpha l^2 \nabla d_L.\BM.\nabla(\delta d_L)\,\mathrm{d}{\Bx}= 0$,\\[0.1cm]
\end{tabular}
\right.
\tag{L}
\label{alLoc}
\end{equation*}

where $\bm\sigma(\bm u_L,d_L;\BM)=\partial_{\bm\varepsilon}w(\Bve_L, d_L, \nabla d_L;\BM)=
{{\bm \sigma}}_{\bm \varepsilon}^{iso}({\bm \varepsilon}_L,d_L)+
{{\bm \sigma}}_{\bm \varepsilon}^{aniso}({\bm \varepsilon}_L,d_L,\BM)$ is defined in Eq. \ref{eq21}, $\delta {\bm u}_L\in {\bf H}^1(\calB_L)$ is the local test function and $\delta d_L\in \text{H}^1(\calB_L)$ is the local phase-field test function.

The variational derivatives of $\widetilde{\mathcal E}$ with respect to $(\bm u_\Gamma,\bm\lambda_C,\bm\lambda_L)$ provide kinematic equations due to weak coupling between global and local form
\begin{equation*}
\widetilde{\mathcal E}_{\bm u_\Gamma}(\bm s;\delta {\bm u}_\Gamma):=\int_\Gamma (\bm\lambda_C+\bm\lambda_L) \cdot \delta {\bm u}_\Gamma \,\mathrm{d}s=0,
\label{Coupl1}
\tag{C$_1$}
\end{equation*}
\begin{equation*}
\widetilde{\mathcal E}_{\bm\lambda_C}(\bm s;\delta {\bm \lambda}_C):=\int_\Gamma (\bm u_\Gamma-\bm u_G) \cdot \delta {\bm \lambda}_C \,\mathrm{d}s=0,
\label{Coupl2}
\tag{C$_2$}
\end{equation*}
\begin{equation*}
\widetilde{\mathcal E}_{\bm\lambda_L}(\bm s;\delta {\bm \lambda}_L):=\int_\Gamma (\bm u_\Gamma-\bm u_L) \cdot \delta {\bm \lambda}_L \,\mathrm{d}s=0.
\label{Coupl3}
\tag{C$_3$}
\end{equation*}
Herein $\delta {\bm u}_\Gamma\in{\bf H}^1(\Gamma)$ and $\delta {\bm \lambda}_C,\delta {\bm \lambda}_L\in{\bf L}^2(\Gamma)$ are the corresponding test functions.
%

Let us now focus on the global variational in (\ref{Global}). The presence of the two domain integrals over $\calB_G$ and $\calB_F$ would imply in this case the need to simultaneously access the corresponding stiffness matrices. Avoiding this can be done as follows: We focus on the domain integral over $\calB_F$ in (\ref{Global}). The idea is to transform the domain integral in $\calB_F$ to the global interface $\Gamma_G$. The divergence theorem leads to
\begin{equation}
\int_{\calB_F}\bm\sigma(\bm u_G):\bm\varepsilon(\delta {\bm u}_G)\,\mathrm{d}{\Bx}
=-\int_{\calB_F}\mathrm{div}(\bm\sigma(\bm u_G))\cdot\delta {\bm u}_G\,\mathrm{d}{\Bx}
+\int_{\partial\calB_F} \bm\sigma(\bm u_G)\cdot\bm n_{\partial\calB_F} \cdot \delta {\bm u}_G\,\mathrm{d}s,
\label{OmegaF}
\end{equation}
where $\bm n_{\partial\calB_F}$ is the unit outward normal vector to $\partial\calB_F$. 

The first term in the right-hand side of in Eq. \ref{OmegaF} can be canceled by using the divergence-free assumption for the stress (no body forces in $\calB_F$). Following a detailed argument in \citet{NoiiGL18}, the second term can be further simplified 
\begin{equation*}
\int_{\partial\calB_F} \bm\sigma(\bm u_G)\cdot\bm n_{\partial\calB_F} \cdot \delta {\bm u}_G\,\mathrm{d}s
=\int_{\Gamma_G} \bm\sigma(\bm u_G)\cdot\bm n_\Gamma \cdot \delta {\bm u}_G\,\mathrm{d}s
+\int_{\partial\calB_F\cap\Gamma_{N,G}} {\bm {\bar\tau}} \cdot \delta {\bm u}_G\,\mathrm{d}s.
\end{equation*}
{Here, $\bm n_\Gamma:=\bm n_{\partial\calB_F}$ denotes the normal vector} on $\Gamma_G$, outward of $\calB_F$, as illustrated in Fig. \ref{Fig3}. Furthermore, it is possible to choose $\calB_L$ and its coarse representation into the global level as $\calB_F$ such that $\partial\calB_F\cap\Gamma_{N,G}=\emptyset$. This is in line with the assumption introduced in section \ref{Section31} that the local domain $\calB_L$ and additionally $\calB_F$ is free from any applied external load. Thus, the last surface integral cancels and (\ref{OmegaF}) can be restated as,
\begin{equation}
\int_{\calB_F}\bm\sigma(\bm u_G):\bm\varepsilon(\delta {\bm u}_G)\,\mathrm{d}{\Bx}
=\int_{\Gamma_G} \bm\sigma(\bm u_G)\cdot\bm n_\Gamma \cdot \delta {\bm u}_G\,\mathrm{d}s,
\label{OmegaF1}
\end{equation}
such that there exists a fictitious Lagrange multiplier $\bm\lambda_F\in{\bf L}^2(\Gamma)$ with
\begin{equation}
\int_{\Gamma_G} \bm\sigma(\bm u_G)\cdot\bm n_\Gamma \cdot \delta {\bm u}_G\,\mathrm{d}s
\label{lambdaF} =:
\int_{\Gamma_G} \bm\lambda_F \cdot \delta {\bm u}_G\,\mathrm{d}s .
\end{equation}
Here,  ${\bm\lambda_F}$ is a traction-like quantity on ${\Gamma_G}$. Due to (\ref{OmegaF1})$-$(\ref{lambdaF}), the {\em partitioned} representation of equation (\ref{Global}) takes the following form
\begin{equation*}
\int_{\calB_G}\bm\sigma(\bm u_G):\bm\varepsilon(\delta {\bm u}_G)\,\mathrm{d}{\Bx}
-\int_{\Gamma_G} \bm\lambda_F \cdot \delta {\bm u}_G\,\mathrm{d}s
-\int_{\Gamma_G} \bm\lambda_C \cdot \delta {\bm u}_G\,\mathrm{d}s
-\int_{\Gamma_{N,G} } {\bm {\bar\tau}} \cdot \delta {\bm u}_G\,\mathrm{d}s=0,
\label{Global1}
\tag{G$_1$}
\end{equation*}
with $\bm\lambda_F$ satisfying
\begin{equation*}
\int_{\Gamma_G} \bm\lambda_F \cdot \delta {\bm u}_G\,\mathrm{d}s
=\int_{\calB_F}\bm\sigma(\bm u_G):\bm\varepsilon(\delta {\bm u}_G)\,\mathrm{d}{\Bx}.
\label{lambdaF1}
\tag{G$_2$}
\end{equation*}
Equations (\ref{Global1}), (\ref{lambdaF1}) refer to the global system of equations. 
The system of equations (\ref{alLocL}) is called a local variational equation and additionally  (\ref{Coupl1}), (\ref{Coupl2}), (\ref{Coupl3}) refer to the coupling terms. The entire system is the basis for the {\em Global-Local approach}. 

\sectpb[Section34]{Dirichlet-Neumann type boundary conditions}
To accommodate a Global-Local computational scheme, instead of finding the stationary solution of the (\ref{Global1}), (\ref{lambdaF1}), (\ref{alLocL}) along with (\ref{Coupl1}), (\ref{Coupl2}), (\ref{Coupl3}) in the monolithic sense, an alternate minimization is used. {This} is in line with \cite{NoiiGL18}, which leads to the Global-Local formulation through the concept of non-intrusiveness. Here the global and local level are solved in a multiplicative manner according to the idea of Schwarz' alternating method \cite{alternatMota17}.

Let $k\geq0$ be the Global-Local iteration index at a fixed loading step $n$. The iterative solution procedure for Global-Local computational scheme is as {follows:}
\begin{itemize}
	\item Dirichlet local problem: solution of local problem (\ref{alLocL}) coupled with (\ref{Coupl3}),
	\item Pre-processing global level: recovery phase using (\ref{Coupl1}) and (\ref{lambdaF1}),
	\item Neumann global problem: solution of global problem (\ref{Global1}),
	\item Post-processing global level: recovery phase using (\ref{Coupl2}).
\end{itemize}
The detailed scheme for applying the Dirichlet-Neumann type boundary conditions to the isotropic phase-field fracture modeling is described in \cite{NoiiGL18}.

\begin{figure}[!ht]
	\centering
	{\includegraphics[clip,trim=1cm 14cm 0cm 11cm, width=16cm]{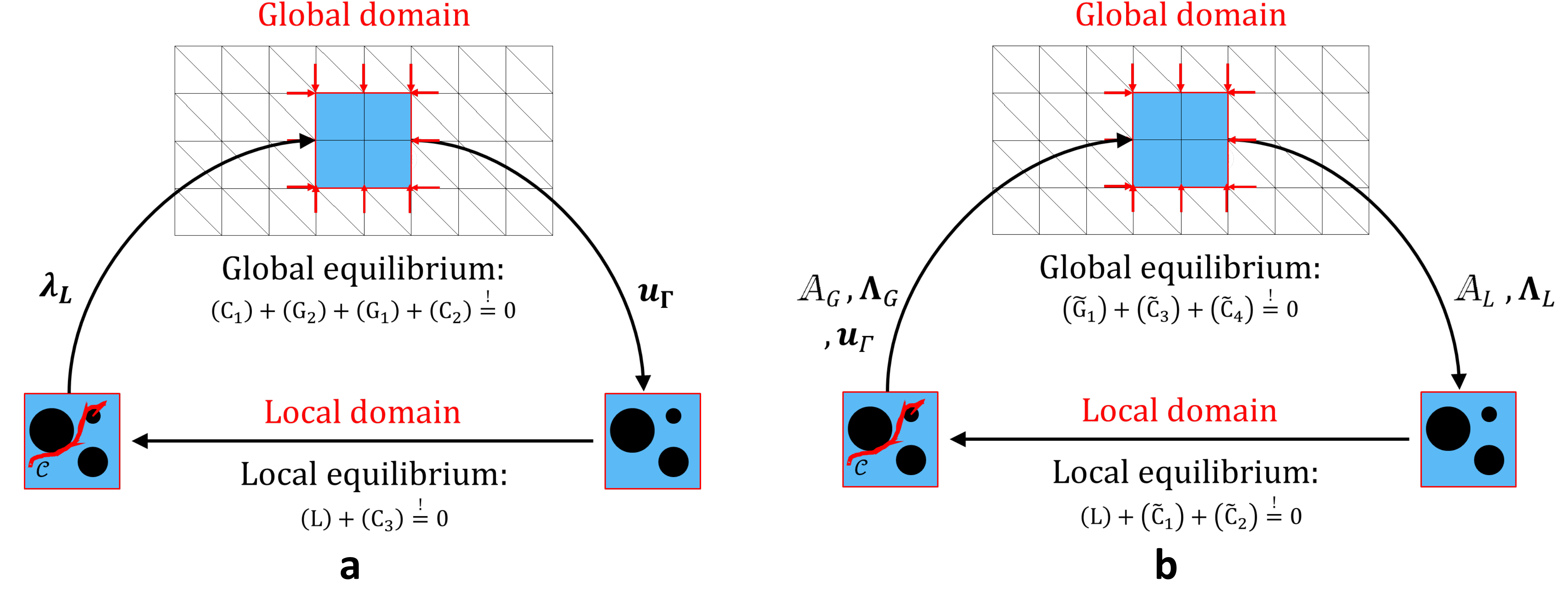}}  
	\caption{Global-Local iterative scheme with $(a)$ Dirichlet-Neumann type boundary conditions; $(b)$ Robin-type boundary conditions.
	}
	\label{Fig4}
\end{figure}

Despite of its strong non-intrusiveness implementation point of view \cite{Allix09}, there are two shortcomings embedded in the system which have to be resolved. ({a}) Due to the extreme difference in stiffness between the local domain and its projection to the global level, i.e. fictitious domain, the relaxation/acceleration techniques has to be used, see \cite{NoiiGL18}. ({b}) Additionally, it turns out that if the solution vector  $(\bm u_G^k,\bm u_L^k,d_L^k,\bm u_\Gamma^k,\bm\lambda_C^k,\bm\lambda_L^k)$ is plugged into equations (\ref{Global1}), (\ref{lambdaF1}), (\ref{alLocL}), (\ref{Coupl1}), (\ref{Coupl2}), (\ref{Coupl3}), the imbalanced quantities follow
\begin{equation}
\int_\Gamma (\bm u_\Gamma^k-\bm u_L^k) \cdot \delta {\bm \lambda}_L  \,\mathrm{d}{s} \neq 0 \quad \texttt{and} \quad \int_\Gamma \bm\lambda_F^k \cdot \delta {\bm u}_G \,\mathrm{d}{s} \neq
\int_{\calB_F}\bm\sigma(\bm u_G^k):\bm\varepsilon(\delta {\bm u}_G)\,\mathrm{d}{\Bx} ,
\label{t2x}
\end{equation}
{resulting in the} {\it iterative} Global-Local computation scheme. Figure \ref{Fig4}a depicts one iteration of the Global-Local approach by means of the Dirichlet-Neumann type boundary conditions. The aforementioned difficulties motivate us to provide an alternative coupling conditions that overcome these challenges, which are explained in the following section.

\sectpb[Section35]{Robin-type boundary conditions}

In this section, {the} Global-Local formulation is enhanced using Robin-type boundary conditions to relax the stiff local response that is observed at the global level (due to the local non-linearity). Furthermore the computational time is reduced. This improves the resolution of the imbalanced quantities in (\ref{t2x}) and it accelerates the Global-Local computational iterations.

Recall, the coupling equations denoted in (\ref{Coupl1}), (\ref{Coupl2}) and (\ref{Coupl3}) arise from the stationary of the interface energy functional. That provides the boundary conditions which have to be imposed on the global and local levels. At that level the Robin-type boundary conditions are formulated.
\begin{itemize}
	\item {\em \textbf{Robin-type boundary conditions at the local level}}		
\end{itemize}
At the local level the new coupling term is introduced as a combination of (\ref{Coupl1}) and (\ref{Coupl2})
\begin{equation}
\widetilde{\mathcal E}_{\bm u_\Gamma}(\bm s;\delta {\bm u}_\Gamma)+\KIA_L\widetilde{\mathcal E}_{\bm\lambda_C}(\bm s;\delta {\bm \lambda}_C)=\int_\Gamma (\bm\lambda_C+\bm\lambda_L) \cdot \delta {\bm u}_\Gamma \,\mathrm{d}s+\KIA_L\int_\Gamma (\bm u_\Gamma-\bm u_G) \cdot \delta {\bm \lambda}_C \,\mathrm{d}s=0 .
\end{equation}
This leads for iteration $k$ to
\begin{equation}
\int_\Gamma (\bm\lambda^{k-1}_C+\bm\lambda^{k}_L) \cdot \delta {\bm u}_\Gamma \,\mathrm{d}s+\KIA_L\int_\Gamma (\bm u_\Gamma^{k,\frac{1}{2}}-\bm u^{k-1}_G) \cdot \delta {\bm \lambda}_C \,\mathrm{d}s=0.
\label{Coup4L}
\end{equation}
Herein, $\KIA_L$ is a local augmented stiffness matrix applied on the interface which serves as regularization of the local Jacobian matrix. By means of  (\ref{Coup4L}) at iteration $k$, the local system of equations results in the following boundary conditions
\begin{equation*}
\int_\Gamma \bm\lambda^{k}_L\cdot \delta {\bm u}_\Gamma \,\mathrm{d}s+\KIA_L\int_\Gamma \bm u_\Gamma^{k,\frac{1}{2}} \cdot \delta {\bm \lambda}_C \,\mathrm{d}s={\bf \Lambda}^{k-1}_L,
\label{Coupl7nL}
\tag{$\widetilde {\text{C}}_1$}
\end{equation*}
\begin{equation*}
\int_\Gamma (\bm u_\Gamma^{k,\frac{1}{2}}-\bm u^{k}_L) \cdot \delta {\bm \lambda}_L \,\mathrm{d}s=0,
\label{Coupl6nL}
\tag{$\widetilde {\text{C}}_2$}
\end{equation*}  
with
\begin{equation}
{\bf \Lambda}^{k-1}_L:=\Lambda_L(\bm\lambda^{k-1}_C,\bm u_G^{k-1};\KIA_L)=
\KIA_L\int_\Gamma \bm u_G^{k-1} \cdot \delta {\bm \lambda}_C \,\mathrm{d}s-\int_\Gamma \bm\lambda^{k-1}_C\cdot \delta {\bm u}_\Gamma \,\mathrm{d}s.
\label{rhs_robin_G}
\end{equation}
Along with (\ref{alLocL}), the local system of equations has to be solved for $(\bm u^{k}_L,\bm\lambda^{k}_L,\bm u_\Gamma^{k,\frac{1}{2}})$ for given local Robin-type parameters $({\bf \Lambda}^{k-1}_L, \KIA_L)$.

\begin{itemize}
\item {\em \textbf{Robin-type boundary conditions at the global level}}		
\end{itemize}
Accordingly, at the global level, the new coupling term is stated as a combination of (\ref{Coupl1}) and (\ref{Coupl3})
\begin{equation}
\widetilde{\mathcal E}_{\bm u_\Gamma}(\bm s;\delta {\bm u}_\Gamma)+\KIA_G\widetilde{\mathcal E}_{\bm\lambda_L}(\bm s;\delta {\bm \lambda}_L)=\int_\Gamma (\bm\lambda_C+\bm\lambda_L) \cdot \delta {\bm u}_\Gamma \,\mathrm{d}s+\KIA_G\int_\Gamma (\bm u_\Gamma-\bm u_L) \cdot \delta {\bm \lambda}_L \,\mathrm{d}s=0 .
\label{Coupl4G}
\end{equation}
This leads for iteration $k$ to
\begin{equation*}
\int_\Gamma (\bm\lambda^{k}_C+\bm\lambda^{k}_L) \cdot \delta {\bm u}_\Gamma \,\mathrm{d}s+\KIA_G\int_\Gamma (\bm u_\Gamma^{k}-\bm u^{k}_L) \cdot \delta {\bm \lambda}_L \,\mathrm{d}s=0,
\label{Coupl5nG}
\end{equation*}
where, $\KIA_G$ is a global augmented stiffness matrix applied on the interface.

Through (\ref{Coupl4G}) at the iteration $k$, the Robin-type boundary condition at the global level follows

\begin{equation*}
\int_\Gamma \bm\lambda^{k}_C\cdot \delta {\bm u}_\Gamma \,\mathrm{d}s+\KIA_G\int_\Gamma \bm u_\Gamma^{k} \cdot \delta {\bm \lambda}_L \,\mathrm{d}s={\bf \Lambda}^{k}_G,
\label{Coupl7nG}
\tag{$\widetilde {\text{C}}_3$}
\end{equation*}
\begin{equation*}
\int_\Gamma (\bm u_\Gamma^{k,\frac{1}{2}}-\bm u^{k}_G) \cdot \delta {\bm \lambda}_C \,\mathrm{d}s=0,
\label{Coupl6nG}
\tag{$\widetilde {\text{C}}_4$}
\end{equation*}
with
\begin{equation}
{\bf \Lambda}^{k}_G:=\Lambda_L(\bm\lambda^{k}_G,\bm u_L^{k};\KIA_G)=
\KIA_G\int_\Gamma \bm u_L^{k} \cdot \delta {\bm \lambda}_L \,\mathrm{d}s-\int_\Gamma \bm\lambda^{k}_L\cdot \delta {\bm u}_\Gamma \,\mathrm{d}s.
\label{rhs_robin_L}
\end{equation}
Together with (\ref{Global1}) and (\ref{lambdaF1}), the global system of equations has to be solved for $(\bm u^{k}_G,\bm\lambda^{k}_C,\bm u_\Gamma^{k})$ for a given $({\bf \Lambda}^{k}_G, \KIA_G,\bm u_\Gamma^{k,\frac{1}{2}})$. Here, $\KIA_G$ and ${\bf \Lambda}^{k}_G$ stand for global Robin-type parameters. 

Based on the new boundary conditions provided in ($\widetilde {\text{C}}_1$), ($\widetilde {\text{C}}_2$), ($\widetilde {\text{C}}_3$) and ($\widetilde {\text{C}}_4$) the imbalanced quantities in the Global-Local iterations read
\begin{equation}
\int_\Gamma (\bm u_\Gamma^k-\bm u_{\Gamma}^{k,\frac{1}{2}}) \cdot \delta {\bm \lambda}_L \,\mathrm{d}{s} \neq 0 \quad \texttt{and} \quad \int_\Gamma \bm\lambda_F^k \cdot \delta {\bm u}_G \,\mathrm{d}{s} \neq
\int_{\calB_F}\bm\sigma(\bm u_G^k):\bm\varepsilon(\delta {\bm u}_G) \,\mathrm{d}{\Bx},
\label{t2}
\end{equation}
For the specific Robin-type boundary conditions, we can resolve Eq. \ref{t2}$_1$ such that this term does not produce any error in the iterative procedure. To do so, following Appendix B, the global and local augmented stiffness matrices within the Robin-type boundary conditions are given by
\begin{equation}\label{eq_aug_stiff}
\fterm{
	\KIA_G={\bm L^T_L}{\bm T^{-T}_L}{\bm {\mathcal{S}}}_L \quad \mbox{and} \quad \KIA_L:={\bm {\mathcal{S}}}_C.
}
\end{equation}
$\KIA_G$ and $\KIA_L$ can be seen as augmented stiffness matrices regularize the Jacobian stiffness matrix at the global and local levels, respectively. 

\begin{Remark}
	\label{diff_scheme_from_robin}
	In the Robin-type boundary condition given in  (\ref{Coupl7nL}) and (\ref{Coupl7nG}), we can extract different criteria, e.g.
	\begin{itemize}
		\item $\KIA_L\rightarrow \infty$: Dirichlet boundary conditions and $\KIA_G\rightarrow 0$: Neumann boundary conditions;
		\item $\KIA_L\rightarrow 0$: Neumann boundary conditions and $\KIA_G\rightarrow \infty$: Dirichlet boundary conditions;
		\item $\KIA_L={\bm {\mathcal{S}}}_C$: Robin-type boundary conditions and $\KIA_G\rightarrow \infty$: Dirichlet boundary conditions;		
	\end{itemize}
	{Hence, depending on the Robin-type parameters, a family of boundary conditions can be formulated.}
\end{Remark}

Additionally to achieve a balance state of Eq. \ref{t2}$_2$, the following partitioned representation of equation (\ref{Global})
\begin{equation*}
\int_{\calB_G}\bm\sigma(\bm u_G):\bm\varepsilon(\delta {\bm u}_G)\,\mathrm{d}{\Bx}
-\int_\Gamma \bm\lambda_F \cdot \delta {\bm u}_G\,\mathrm{d}s
-\int_\Gamma \bm\lambda_C \cdot \delta {\bm u}_G\,\mathrm{d}s
-\int_{\Gamma_{N,G} } {\bm {\bar\tau}} \cdot \delta {\bm u}_G\,\mathrm{d}s=0,
\label{GlobalN1}
\tag{$\widetilde {\text{G}}_1$}
\end{equation*}
is equipped with a linearized $\Delta\bm\lambda_F$ satisfying
\begin{equation*}
\int_\Gamma \Delta\bm\lambda_F \cdot \delta {\bm u}_G\,\mathrm{d}s
=\int_{\Gamma} {\bm{\mathcal{S}}_F}{\Delta{{\bm u}}}_{G}\cdot{\delta{{\bm u}}}_{G} \,\mathrm{d}{\Bx},
\label{lambdaFN1}
\tag{$\widetilde {\text{G}}_2$}
\end{equation*}
where (\ref{eqdefA51}) in Appendix B is used with ${{{\bm u}}_F}={{{\bm u}}_G}$. Note that now the second imbalance quantity shown in (\ref{t2}) does not anymore produce an error. We are not solving for $\bm\lambda_F$ and in the linearized equation of (\ref{GlobalN1}) this term is replaced by (\ref{lambdaFN1}). The  linearized equation of (\ref{GlobalN1}) is solved within a single iteration, because we are dealing with a linear elastic constitutive equations.

The detailed Global-Local formulation using Robin-type boundary 
conditions is depicted in {Algorithm} \ref{TableGL}. Accordingly, Fig. \ref{Fig4}b depicts one iteration of the Global-Local coupling scheme by means of the Robin-type boundary conditions. The Global-Local setting provides a generic two-scale finite element algorithms that enables capturing local non-linearities.

\sectpb[Section36]{Spatial discretization}
The computational domain is subdivided into bilinear quadrilateral elements denoted as ${Q}_1$. Both subproblems are discretized with a Galerkin finite element method 
using $H^1$-conforming bilinear (2D) elements, i.e., the
ansatz and test space uses $Q_1^c$--finite elements, e.g., for details, we refer readers to the \cite{Cia87}. Consequently, the discrete spaces have the property $V_h \subset V$ and $W_h \subset W$. Here, $h$ refers to the finite element size. Accordingly, a finite element discretization is illustrated in Appendix A for the primal fields refers to the  $(\bm u_G,\bm u_L,d_L,\bm u_\Gamma,\bm\lambda_C,\bm\lambda_L)$ and its constitutive state variables represented by $(\Bve_{G}, \Bve_{L}, d_L, \nabla d_L)$.

\newpage
\begin{algorithm}[H]\small
	\caption{\em Global-Local iterative scheme combined with Robin-type boundary conditions.}
	\label{TableGL}
	{
		\begin{tabular}{l}
			{\bf Input:} loading data $(\bar{\bm u}_{n},\bar{\bm \tau}_{n})$ on $\partial_D \calB$ and $\Gamma_{N}$, respectively; \\
			\hspace{1.37cm}solution $(\bm u_{G,n-1},\bm u_{L,n-1},d_{L,n-1},\bm u_{\Gamma,n-1},\bm\lambda_{C,n-1},\bm\lambda_{L,n-1})$ and $\calH_{L,n-1}$ from step $n-1$. \\ [0.5cm]
			
			{{\bf Global-Local iteration} $k\geq1$}:\\[0.2cm]
			
			\quad\quad \underline{Local boundary value problem}:\\
			
			\quad\quad\quad\;\; \textbullet\; given $\KIA_L,{\bf \Lambda}^{k-1}_L, \calH_{L,n-1}$; solve \\
			\qquad\qquad\qquad\;\;${\text{phase-field part:}}$\\
			\qquad\qquad\qquad\;\;\;\;$\displaystyle (1-\kappa)\int_{\calB_L} d_{L+} \calH(\bm\varepsilon(\bm u_L)). \delta d_L\mathrm{d}{\bm{x}}+ G_c \int_{\calB_L}  \frac{1}{l}(d_L-1).\delta d_L \mathrm{d}{\bm{x}}$\\[0.1cm]
			\qquad \qquad \qquad \qquad $ \displaystyle + G_c \int_{\calB_L} l \nabla d_L.\nabla(\delta d_L) \mathrm{d}{\bm{x}} + G_c \int_{\calB_L} \alpha l \nabla d_L.\BM.\nabla(\delta d_L) \mathrm{d}{\bm{x}}= 0$,\\
			
			\qquad\qquad\qquad\;\;${\text{mechanical part:}}$\\
			\hspace{2.4cm}\quad
			$\begin{cases} 
			\displaystyle \int_{\Omega_L}\bm\sigma(\bm u_L,d_L):\bm\varepsilon(\delta {\bm u}_L)\,\mathrm{d}{\Bx}-\int_\Gamma \bm\lambda_L \cdot \delta {\bm u}_L \,\mathrm{d}{s}=0, \\
			\displaystyle \int_\Gamma \bm\lambda_L\cdot \delta {\bm u}_\Gamma \,\mathrm{d}s+\KIA_L\int_\Gamma \bm u_\Gamma \cdot \delta {\bm \lambda}_C \,\mathrm{d}s={\bf \Lambda}^{k-1}_L , \\
			\displaystyle \int_\Gamma (\bm u_\Gamma-\bm u_L) \cdot \delta {\bm \lambda}_L \,\mathrm{d}s=0,		
			\end{cases}$ \\
			\quad\quad\quad \;\;{\color{white}\textbullet}\; \;set $(\bm u_L,d_L,\bm u_\Gamma,\bm\lambda_L)=:(\bm u_L^k,d_L^k,\bm u_\Gamma^{k,\frac{1}{2}},\bm\lambda_L^k)$, \\
			
			\quad\quad\qquad \textbullet\; given $({\bm u_L^{k}},{\bm\lambda^{k}_L};\KIA_G)$, set\\
			\hspace{2.7cm} $\displaystyle {\bf \Lambda}^{k}_G
			=\KIA_G\int_\Gamma \bm u_L^{k} \cdot \delta {\bm \lambda}_C \,\mathrm{d}s-\int_\Gamma \bm\lambda^{k}_L\cdot \delta {\bm u}_\Gamma \,\mathrm{d}s$. \\[0.5cm]
			
			\quad\quad \underline{Global boundary value problem}:\\
			
			\quad\quad\quad\;\; \textbullet\; given $\KIA_G,{\bf \Lambda}^{k}_G,\bm u_\Gamma^{k,\frac{1}{2}}$, solve \\	
			\hspace{2.4cm}\quad
			$\begin{cases} 
			\displaystyle \int_{\Omega_G}\bm\sigma(\bm u_G):\bm\varepsilon(\delta {\bm u}_G) \,\mathrm{d}{\Bx}
			-\int_\Gamma \bm\lambda_F \cdot \delta {\bm u}_G \,\mathrm{d}{s}
			-\int_\Gamma \bm\lambda_C \cdot \delta {\bm u}_G \,\mathrm{d}{s}
			-\int_{\Gamma_{N} } \bar{\bm \tau} \cdot \delta {\bm u}_G \,\mathrm{d}{s}=0,  \\
			\displaystyle \int_\Gamma \bm\lambda_C\cdot \delta {\bm u}_\Gamma \,\mathrm{d}s+\KIA_G\int_\Gamma \bm u_\Gamma \cdot \delta {\bm \lambda}_C \,\mathrm{d}s={\bf \Lambda}^{k}_G , \\
			\displaystyle \int_\Gamma (\bm u_\Gamma^{k,\frac{1}{2}}-\bm u_G) \cdot \delta {\bm \lambda}_C \,\mathrm{d}s=0,	
			\end{cases}$ \\
			
			\quad\quad\quad\;\; {\color{white}\textbullet}\; set $(\bm u_G,\bm u_\Gamma,\bm\lambda_C)=:(\bm u_G^k,\bm u_\Gamma^{k},\bm\lambda_C^k)$, \\

			\quad\quad\qquad \textbullet\; given $({\bm u_G^{k}},{\bm\lambda^{k}_C};\KIA_L)$, set \\
			\hspace{2.7cm} $\displaystyle {\bf \Lambda}^{k}_L
			=\KIA_L\int_\Gamma \bm u_G^{k} \cdot \delta {\bm \lambda}_C \,\mathrm{d}s-\int_\Gamma \bm\lambda^{k}_C\cdot \delta {\bm u}_\Gamma \,\mathrm{d}s$. \\

			\quad\quad\quad\;\; \textbullet\; if fulfilled, set $(\bm u_G^k,\bm u_L^k,d_L^k,\bm u_\Gamma^k,\bm\lambda_C^k,\bm\lambda_L^k)=:(\bm u_{G,n},\bm u_{L,n},d_{L,n},\bm u_{\Gamma,n},\bm\lambda_{C,n},\bm\lambda_{L,n})$ and stop; \\
			\quad\quad\quad\;\; {\color{white}\textbullet}\; else $k+1\rightarrow k$. \\[0.5cm]
			
			{\bf Output:} solution $(\bm u_{G,n},\bm u_{L,n},d_{L,n},\bm u_{\Gamma,n},\bm\lambda_{C,n},\bm\lambda_{L,n})$ and $\calH_{L,n}$.
%
		\end{tabular}
	}
\end{algorithm}

\sectpa[Section4]{Predictor-Corrector Adaptivity Applied to the Global-Local
	Formulation}

We assume the Global-Local formulation is at the converged state, which is denoted as $(\bm u_{G,n},\bm u_{L,n},d_{L,n}$
$,\bm u_{\Gamma,n},\bm\lambda_{C,n},\bm\lambda_{L,n})$. The
Global-Local approach is augmented by a {\it dynamic allocation} of a local state using an adaptive scheme which has to be performed at
time step $t_{n}$. By the adaptivity procedure, we  mean: ({a}) to
determine which global elements need to be refined and identified by
$\calB^{\;adapt}_G\subset\calB_G$; ({b}) to create the new fictitious
domain $\calB^{\;new}_F:=\calB^{\;old}_F\cup\calB^{\;adapt}_F$ with $\calB^{\;adapt}_F:= \calB^{\;adapt}_G$ and as a result
a new local domain is defined as $\calB^{\;new}_L:=\calB^{\;old}_L\cup\calB^{\;adapt}_L$, see Fig. \ref{adapt}; ({c}) to determine a new local interface denoted as $\Gamma_L$;
({d}) to interpolate the old global solution in $\calB^{\;adapt}_L$. {All these steps refer to predictor steps. 
The corrector step is explained in Section \ref{Section4_a2}.}
{We briefly notice that the principal idea of this adaptive scheme 
is inspired from \cite{Wick15Adapt} in which a predictor-corrector scheme for 
mesh refinement in the crack zone was proposed.}
 
\sectpb[Section4_a1]{Predictor step}
{In this section, we start explaining the predictor step.}
\begin{itemize}
 	\item {\em \textbf{Determining global elements which have to be refined}}		
\end{itemize}
Recall that the {interfaces} at the global, fictitious and local domains are denoted by $\Gamma_G \subset \calB_G$, $\Gamma_F \subset \calB_F$ and $\Gamma_L \subset \calB_L$. We denote $e_{G}\in\calB_G$, $e_{F}\in\calB_F$ and $e_{L}\in\calB_L$ as the elements in the global, fictitious and local domain. Let $E^1_{G,i}, E^2_{G,i}, E^3_{G,i}$ and $E^4_{G,i}$ refer to the left, top, right and bottom global edges for the $i^{th}$ element $e_{G,i}$, respectively (because, it is quadrilateral hence it has four edges). Accordingly, $E^q_{F,i}\in e_{F,i}$ and $E^q_{L,i}\in e_{L,i}$ with $\;q=(1,2,3,4)$ refer to the fictitious and local edges, see Fig. \ref{adapt}.  


We now develop a procedure, to determine the global elements $e_{G,i}$
which have at least one edge $E^{q}_{G,i}$ such that their fine resolution at the
local level, i.e. $E^{q}_{L,i}$, reaches to the crack phase-field threshold
value. Thus $e_{G,i}$ has to be refined. The predictor step for
the adaptive scheme of the Global-Local formulation {is explained} in Algorithm \ref{alg_1}.

\begin{algorithm}[H]\small
	\caption{\em Predictor step for the adaptive procedure.}	
	\label{alg_1}
	Let $0<\texttt{TOL}_{d}<1$ be given. For the $d_L:=d_{L,n}\;\in\calB_L$, 
find corresponding $e_{G}\in\calB_G$ which {must be} refined using the following steps:
	\begin{enumerate}
		\item Find $\bm x_L\in\Gamma_L$ such that $d_L({\bm
	x}_L)<\texttt{TOL}_{d}$ on $\Gamma_L$:\\[2mm] 
{\it Checking criterion:} If "\texttt{Yes}" proceed {to step 2.}
If "\texttt{No}" stop,
		\item find $E^{q}_{L,i}\in\Gamma_L$ such that ${\bm x_L}\in E^{q}_{L,i}$,
		\item find $E^{q}_{G,i}=\calP^{-1}E^{q}_{L,i}$ (corresponding edge {in} $\calB_G$),
		\item find $e_{G}\in\calB_G$ and $e_{G}\not\in\calB_F$ such that $E^{q}_{G,i}\in e_{G}$.
    \end{enumerate}
\end{algorithm}

Here, $\calP:\Gamma_G\rightarrow\Gamma_L$ is denoted as a projection/geometrical operator which maps the global to the local interface by $E^{q}_{L,i}:=\calP E^{q}_{G,i}$. 

\begin{itemize}
	\item {\em \textbf{Creating new fictitious and local domains: $(\calB^{\;new}_F,\calB^{\;new}_L)$}}		
\end{itemize}
We are now able to determine a new fictitious and local domains. Knowing $e_{G}$ from the previous step, the new fictitious domain is $\calB^{\;new}_F:=\calB^{\;old}_F\cup\calB^{\;adapt}_F$ such that $\calB^{\;adapt}_F:=e_{G}$. As a result, a new local domain is $\calB^{\;new}_L:=\calB^{\;old}_L\cup\calB^{\;adapt}_L$ such that $\calB^{\;adapt}_L$ is a fine discretization (including heterogeneity as well) of $\calB^{\;adapt}_F$, see Fig. \ref{adapt}.
\begin{itemize}
	\item {\em \textbf{Determining the coupling interface: $(\Gamma_G,\;\Gamma_F,\;\Gamma_L)$}}		
\end{itemize}
{So far, we have identified {\it new} fictitious and local domains. Next, we determine $\Gamma_F\subset\calB_F$, due to its coarse discretization. Afterwards, we find the local interface $\Gamma_L \subset \calB_L$ by projecting $\Gamma_F$ to $\calB_L$. Finally, it is trivial that $\Gamma_G:=\Gamma_F$, because $\calB_F$ and $\calB_G$ are in the same discretization space. 

The edge $E^q_{F,i}\in\Gamma_F$ is on the interface if 
\begin{equation}\label{inter1}
E^{q}_{F,i}\in\Gamma_F \;\; \mbox{if} \;\;\nexists \;\;i,j\;:E^{q}_{F,i}=E^{r}_{F,j} ,\quad \mbox{for}\;q,r=(1,2,3,4),
\end{equation}
{which} means if an edge $E^q_{F,i}$ is shared between two elements in $\calB_F$, it is not on the interface (inner edge) and if it belongs only to one element, 
then {it must also} belong to interface (outer edge). As a result, we define the fictitious interface as $\Gamma_F = \substackrel{i,q}{{\huge \texttt{A}}} \, E^{q}_{F,i}$, hence $\Gamma_L = \substackrel{i,k}{{\huge \texttt{A}}} \, E^{q}_{L,i}$ with $ E^{q}_{L,i}:=\calP E^{q}_{F,i}$. 
\begin{itemize}
	\item {\em \textbf{Interpolating the old solution at $t^n$ {from
	the global to the local mesh}}}		
\end{itemize}
Given a continuous function $\bm u^{adapt}_G$ that is $\bm u_G$ in $\calB^{\;adapt}_G$, we define the linear interpolation operator  $\mathcal{\pi}:\calB^{\;adapt}_G\rightarrow\calB^{\;adapt}_L$ to $\bm u^{adapt}_G$ by
\begin{equation}\label{inter3}
{\bm u^{adapt}_{L,n}}(\bm{x}):=\mathcal{\pi}{\bm u^{adapt}_{G,n}}=\bm N_u^G(\bm{x}){\hat{\bm u}^{adapt}_{G,n}}\quad\mbox{for}\;\; {\bm{x}}\in\calB^{\;adapt}_L. 
\end{equation}
where $\bm N_u^G(x)$ is defined in Appendix A. Hence, ${\bm u^{new}_{L,n}}:={\bm u_{L,n}}\cup{\bm u^{adapt}_{L,n}}$.

\sectpb[Section4_a2]{{Corrector step}}


We introduce a corrector step {in which the computation is rerun on the newly determined local mesh.} To this end, we compute Global-Local
solutions until {the} checking criterion in Algorithm \ref{alg_1} is not
satisfied (\texttt{No}). 
That means we could not find additional local edges on the interface such that $d_L({\bm x}_L)<\texttt{TOL}_{d}$ on $\Gamma_L$ holds.

{Let us write Algorithm \ref{TableGL} in} the following abstract form
\begin{equation}\label{inter4}
{\bm s_n}=\texttt{GL}({\bm s_{n-1}}),
\end{equation}
with ${\bm s}=(\bm u_{G},\bm u_{L},d_{L},\bm u_{\Gamma},\bm\lambda_{C},\bm\lambda_{L})$. We define an intermediate solution ${\bm s^{j}_{n-1}}$ at fixed $t^{n}$ such that the corrector step for adaptive scheme reads,
\begin{equation}\label{inter5}
{\bm s^{j}_{n-1}}=\texttt{GL}({\bm s_{n-1}}). 
\end{equation}
Perform Algorithm \ref{alg_1}, if the \textit{checking criterion} in step 1 is
satisfied. But, if this is not the case, then the corrector step is fulfilled, thus set ${\bm
s^{j}_{n-1}}=:{\bm s_{n}}$ and stop; else $j+1\rightarrow j$.

\sectpb[Section4_a3]{{The final predictor-corrector scheme}}

The aforementioned predictor-corrector adaptivity procedure is summarized in Algorithm \ref{alg_2} as follows,
\begin{algorithm}[H]\small
	\caption{\em Predictor-corrector steps for the adaptive procedure.}
	\label{alg_2}
	Let $0<\texttt{TOL}_{d}<1$ be given. For the $d_L:=d_{L,n}\;\in\calB_L$, find corresponding $e_{G}\in\calB_G$ which have to be refined using the following steps:
	\begin{enumerate}
		\item Compute the intermediate solution by ${\bm s^{j}_{n-1}}=\texttt{GL}({\bm s_{n-1}})$,
		\item perform Algorithm \ref{alg_1} if \textit{Checking criterion} is satisfied (\texttt{Yes}),
		\item if \textit{checking criterion} in Algorithm \ref{alg_1} is not satisfied, then the corrector step is fulfilled, thus set ${\bm s^{j}_{n-1}}=:{\bm s_{n}}$ and stop; else $j+1\rightarrow j$.
	\end{enumerate}
\end{algorithm}	
Fig. \ref{adapt} depicts one iteration of the predictor-corrector steps for the adaptive procedure which is illustrated in Algorithm \ref{alg_2}.
\begin{figure}[!ht]
	\centering
	{\includegraphics[clip,trim=0cm 3cm 0cm 4.5cm, width=14cm]{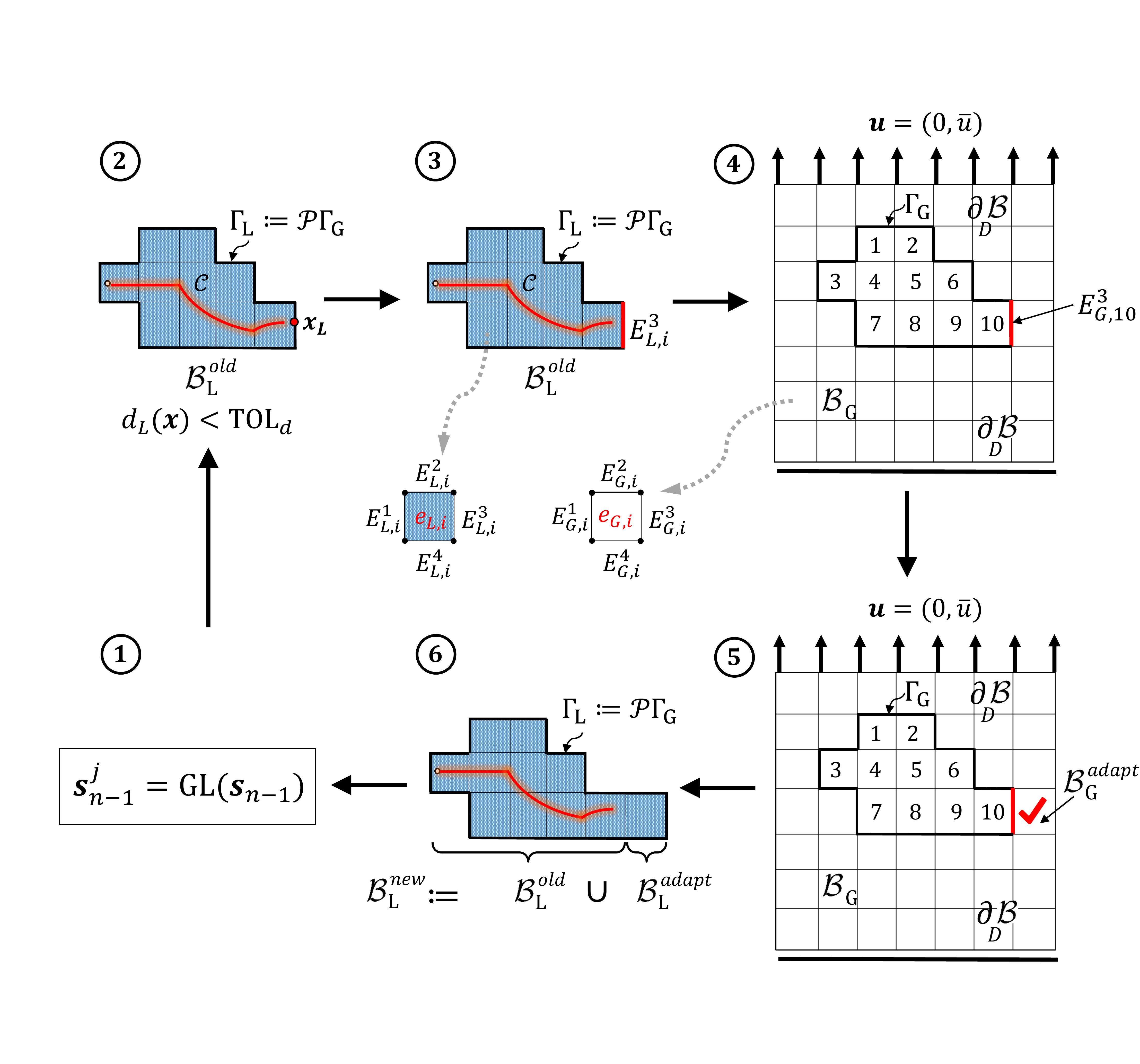}}  
	\caption{{Explanation of the predictor-corrector adaptive scheme.}
	}
	\label{adapt}
\end{figure}

We notice that for the brutal fracture behavior, where a complete
failure {happens} in one load increment, 
Algorithm \ref{alg_2} has an important effect. This is mainly because the corrector step is performed until there is no nodal point on the interface such that $d_L({\bm x}_L)<\texttt{TOL}_{d}$ holds. If this is fulfilled, the adaptivity procedure stops and goes to the next load increment. Thus brutal fracture can be observed in the next time step. This is illustrated in the numerical example of the Section \ref{Section522}. 

The performance of Algorithm \ref{alg_1} and Algorithm \ref{alg_2} is depicted in Fig. \ref{A1} and Fig. \ref{A2}. This example refers to the isotropic single-edge notch under tension. The numerical setup is given in \cite{MieWelHof10b}. By applying predictor-corrector steps, we will have a more regularized fracture surface. This is observed for instance in Fig. \ref{A2} step 58, step 60 and step 62 (right figures). But that is not the case, if we only apply predictor step. For instance in Fig. \ref{A1}, at step 58, after the predictor step (in the absence of the corrector stage). Here we do not have a regularized fracture surface. It will be regularized in the next load increment which is shown in step 59 . That is also observed in steps 60 and 63 in Fig. \ref{A1}, as well. To compare these effects on the global level, we refer to the load-displacement curve in Fig. \ref{B4}a. It was observed that the corrector procedure applied to the predictor step improved the Global-Local results.

\begin{figure}[!ht]
	\centering
	{\includegraphics[clip,trim=0cm 9cm 0cm 10cm, width=16cm]{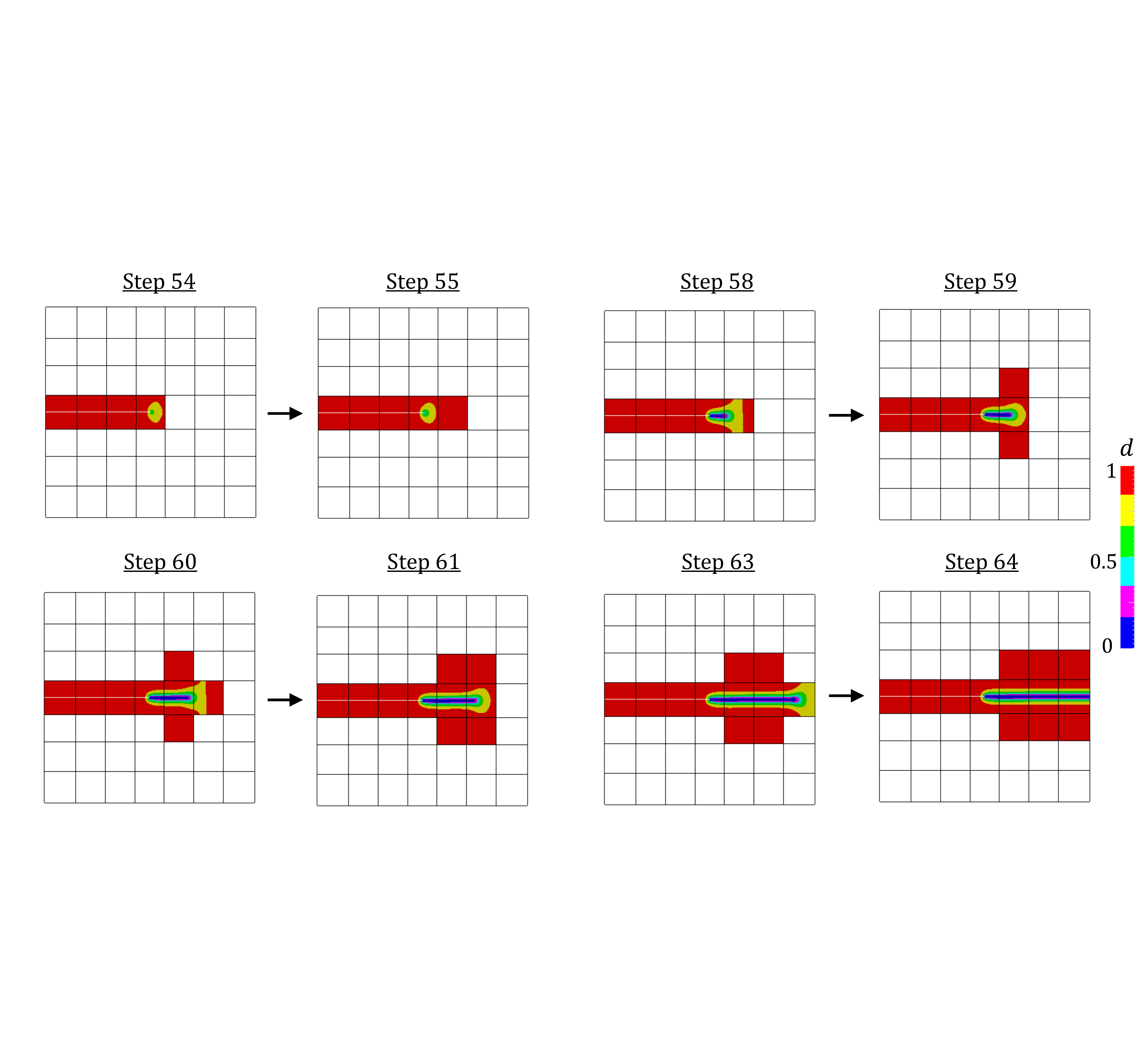}}  
	\caption{Global-Local approach augmented with the predictor adaptive scheme given in Algorithm \ref{alg_1}; Fracture of isotropic single-edge-notched plate under tension per loading steps.
	}
	\label{A1}
\end{figure}

\begin{figure}[!ht]
	\centering
	{\includegraphics[clip,trim=0cm 9cm 0cm 8.8cm, width=16cm]{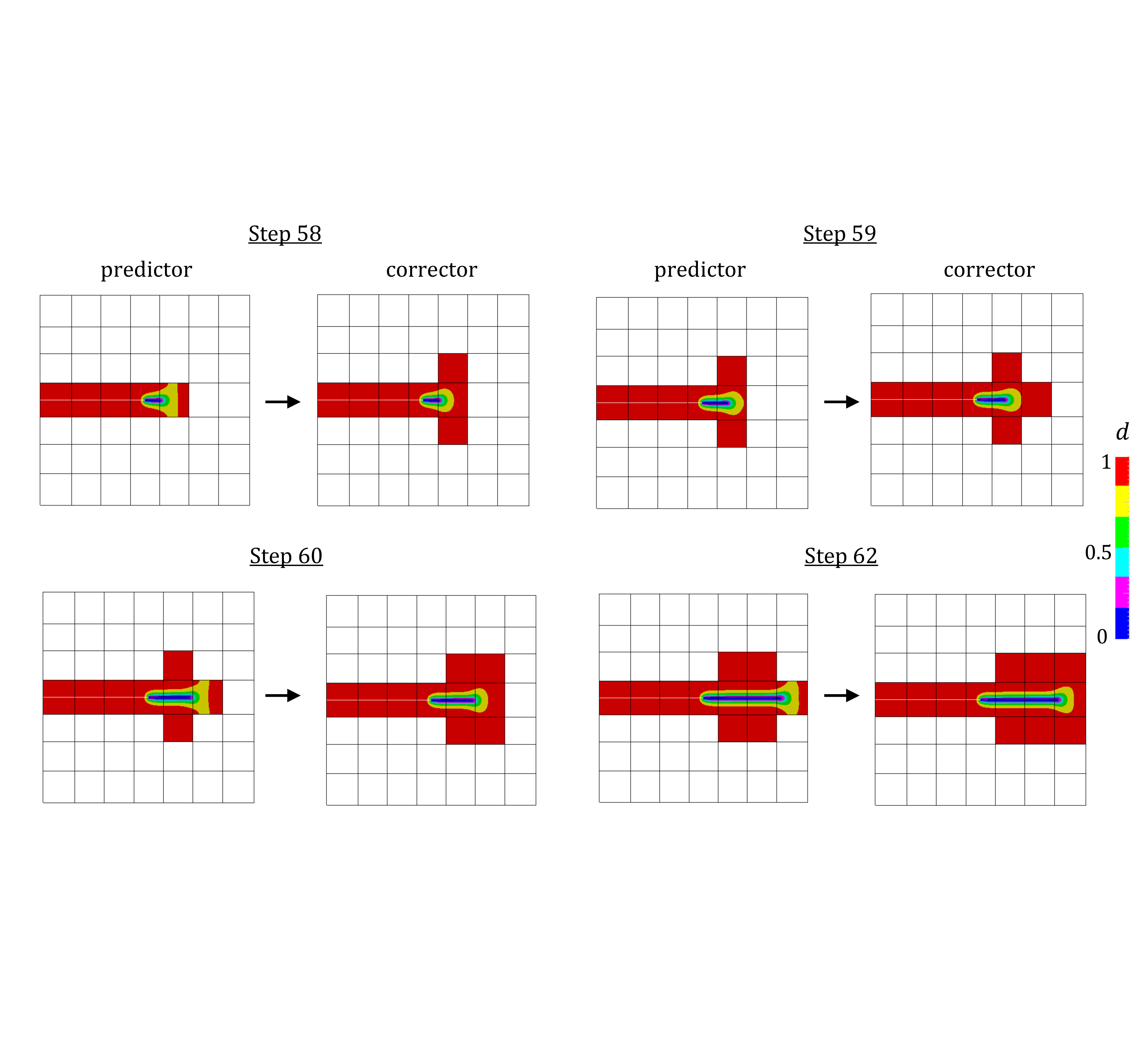}}  
	\caption{Global-Local approach augmented with the predictor-corrector adaptive scheme indicated in Algorithm \ref{alg_2}; Fracture of isotropic single-edge-notched plate under tension per loading steps.
	}
	\label{A2}
\end{figure}

\sectpb[Section4_homoPhF]{Homogenized {phase-field solution on the global level}}
We determine the coarse representation of the crack phase-field. Here we denote $(\bm u_{G,n},\bm u_{L,n},d_{L,n}$
$,\bm u_{\Gamma,n},\bm\lambda_{C,n},\bm\lambda_{L,n})$ to be the converged solution of the Global-Local approach. We emphasis that the computation of the global phase-field $d_G$ is a post-processing step. The \textit{homogenized global crack phase-field} solution can be determine based on the following ways.

(\textit{a}) \textbf{Global crack phase-field solution.} We solve the crack phase-field given in (\ref{eq282}) after having obtained the converged global solution by
\begin{equation}\label{global_crack1}
	d_G =\substackrel{d_G \in W }{\mbox{argmin}} \,
	[-g^{\prime}({d_{G+})\calH(\bm u_{G,n};\BM_G) - l_G \delta_d.\gamma_l}],
\end{equation}
with,
\begin{equation}\label{global_crack2}
l_G := l_L \frac{h_G}{h_L}\quad \mbox{and}\quad\BM_G := \BM_L. 
\end{equation}
where $\calH$ is described by (\ref{eq29_H3}) and formulated on the global level. The last condition in (\ref{global_crack2})$_2$ holds because we assumed the structural tensor at the global level inclined with identical angle as the local level. This holds in the case of the transverse isotropic setting.

(\textit{b}) \textbf{Homogeneous crack phase-field solution.} Assume that at the global level the transition zone of the crack phase-field vanishes. That results in a free isotropic and anisotropic Laplacian operator $\Delta (\bullet)=0$ in (\ref{eq282}). Hence, the second and third terms in (\ref{eq29_H5}) become zero. Following that Eq. \ref{eq282} at the quasi-static stationery state is restated by  
\begin{equation}\label{homo_phf_sol}
d_G=\frac{1}{1+2(1-\kappa)\calH(\bm u_{G,n};\BM_G)} \;\in[0,1],
\end{equation}
{such} that if $\calH\rightarrow\infty \;\mbox{then}\;d_G=0$ and for $\calH\rightarrow 0  \;\mbox{then}\;d_G=1$ holds. The homogeneous crack phase-field solution $d_G$ in (\ref{homo_phf_sol}) is independent of the $l_G$ but the global crack phase-field solution in (\ref{global_crack1}) depends on $l_G$.
At the global level we do not have any given imperfection (e.g. notch shaped). That is located \textit{only} at the local level. However, (\ref{global_crack1}) or (\ref{homo_phf_sol}) still provides the desired crack direction because of $\calH(\bm u_{G,n};\BM_G)$. This is due to the fact that $\bm u_{G,n}$ is determined based on the given $\bm u_\Gamma^{k,\frac{1}{2}}$ which is upscaled from local level, see (\ref{Coupl6nG}). As a result, the crack driving force at the global level is the \textit{true projection} of the constitutive non-linearity at the local level. That is varitionally consistent and resulting from the upscaling procedure (i.e. information which is transformed from the local level to the global scale).

%
%
%
%
%
\sectpa[Section5]{Numerical Examples}

The section presents the performance of the proposed adaptive
Global-Local approach applied to the phase-field modeling of anisotropic
brittle fracture. We consider {four} numerical model problems. The first example
deals with an isotropic single-edge-notched shear test in which we set {the} directional
tensor to be zero. The next {three} examples deal with transverse isotropic
setting with different directional {tensors}.

\sectpb[Section5_a]{{Goals of the computations}}

For comparison purposes, we compute quantitative and qualitative single scale and Global-Local resolutions. In detail, we investigate: 

\begin{itemize}
	\item Crack patterns on the local scale at the complete fracture state in order to evaluate the down-scaling procedure (i.e. transition of external loading increments from the global scale to the local level);
	\item Load-displacement curves to evaluate the up-scaling procedure during the Global-Local coupling approach (i.e. transition of local non-linearity and heterogeneity
	responses to the global level);
	\item Investigations of the thermodynamically consistency between the single scale strain-energy and its Global-Local energy functional;
	\item Efficiency of the overall response resulting from {the} predictor-corrector adaptive scheme;
	\item Effect of the given threshold phase-field value in the adaptive process for the derivation of the fracture zone;
	\item Evaluating the homogenized phase-field solution at the global level, when the Global-Local scheme is in the converged state.
\end{itemize}

The outlined constitutive formulation is considered to be a canonically
consistent and robust scheme for capturing the non-linearities on the
lower level and its projection on to the global level.

\sectpb[Section5_b]{{Geometry, data and solution procedures}}

As a setup for the numerical investigations, we use:

\begin{itemize}
	\item \textbf{\underline{Geometries and parameters:}} In the first two examples, a boundary value problem applied to the square plate is shown in Fig. \ref{B1}. We set $A=0.5\;mm$ hence $\calB=(0,1)^2$ $mm^2$ that includes a predefined single notch from  the left edge to the body center, as depicted in Fig. \ref{B1}. The predefined crack is in the $y=A$ plane and is restricted in $0\le|\mathcal{C}|\le l_0$ and we set $l_0=A=0.5$. As a loading setup, we set the initial values for displacement and phase-field as $\bm u_0:=0 \in\calB$ and $d_0:=1 \in \calB_C$ and $\calB_L$. The finite element discretization is explained in Section \ref{Section36}. {Details of the last two examples are accordingly given.}
	
	\item \textbf{\underline{Material parameters:}} In the first two examples, the constitutive parameters for the isotropic and transversely isotropic material are the same as in \cite{miehe+welschinger+hofacker10a} and given as $\lambda = 121.15$ kN/mm$^2$, $\mu = 80.77$ kN/mm$^2$. Griffith's critical elastic energy release rate is set as $G_c = 2.7 \times 10^{-3}$ kN/mm. In the first example, the preferred fiber direction is set to zero ${\bm{a}}=0$ which represents the standard isotropic setting. Whereas, the three other examples represent transversely isotropic behavior that is characterized by the symmetric second-order structural tensor $\BM={\bm{a}}\otimes{\bm{a}}$. For a given angle $\phi$ of the preferred fiber direction, the normal vector is defined as ${\bm{a}}:=~\big[cos(\phi) \quad sin(\phi) \big]^T$. For the second example, the preferred fiber direction is given by the structural director $\Ba$ which is inclined by $\theta = +30^\circ$ and $\theta = -30^\circ$ with respect to the $x$-axis of a fixed Cartesian coordinate system in the Sections \ref{Section521} and \ref{Section522}, respectively. Additionally, the anisotropy penalty-like parameters for the deformation part and fracture contribution are set to $\alpha= \chi = 50$ and $\Xi =0$. All material properties are fixed for the following numerical examples, unless indicated otherwise.
 
    \item \textbf{\underline{Model parameters:}} The phase-field parameters
    are chosen as $\kappa = 10^{-10}$ and $l = 2h$. The threshold value for
    the Global-Local predictor-corrector mesh refinement scheme is
    $\texttt{TOL}_{d} = 0.85$. This threshold value $\texttt{TOL}_{d}$ is a
    fixed value except for {the} compression {cases} in which we {use a} different $\texttt{TOL}_{d}$. 
    
    \item \textbf{\underline{Solution of the nonlinear problems:}}
    
    {An} alternate minimzation scheme is used for solving the local boundary value problem indicated in Algorithm \ref{TableGL}. Thus, we alternately {solve} for $d_L$ by fixing ${\bm u}_L$ and then solving for $({\bm u}_L,{\bm u}_{\Gamma},{\bm \lambda_L})$ by fixing $d_L$ until convergence is reached.  
    An iterative Newton solver is used in which the linear equation systems are solved
    with  a  generalized  minimal
    residual {method}. The {stopping criterion of the single scale
    and local Newton methods is $\texttt{Tol}_\texttt{N-R}=10^{-10}$. Specifically,
    the relative residual norm is given by $\texttt{Residual}: \| \bm F(\bm x_{k+1}) \|
\leq \texttt{Tol}_\texttt{N-R} \| \bm F(\bm x_{k}) \|$.}
Here, $\bm F$ refers to the residual of the equilibrium equation of {the
nonlinear} single scale and local boundary value problems. 

    \item \textbf{\underline{Software:}}    
    
    The implementation is based on \textsc{MATLAB R2018}b \cite{MATLAB18b}
    and \textsc{Fortran 90} \cite{Chapman03}. The user elements including the constitutive
    modeling at each Gaussian quadrature points are written in \textsc{Fortran
    90}. The general framework for the Global-Local approach {is} implemented in \textsc{MATLAB} as a parent/main program such that all subprograms in \textsc{Fortran 90} {are} called as a \textsc{Mex-file}. 
    
\end{itemize}

\begin{figure}[!ht]
	\centering
	{\includegraphics[clip,trim=0cm 13cm 2cm 11cm, width=16cm]{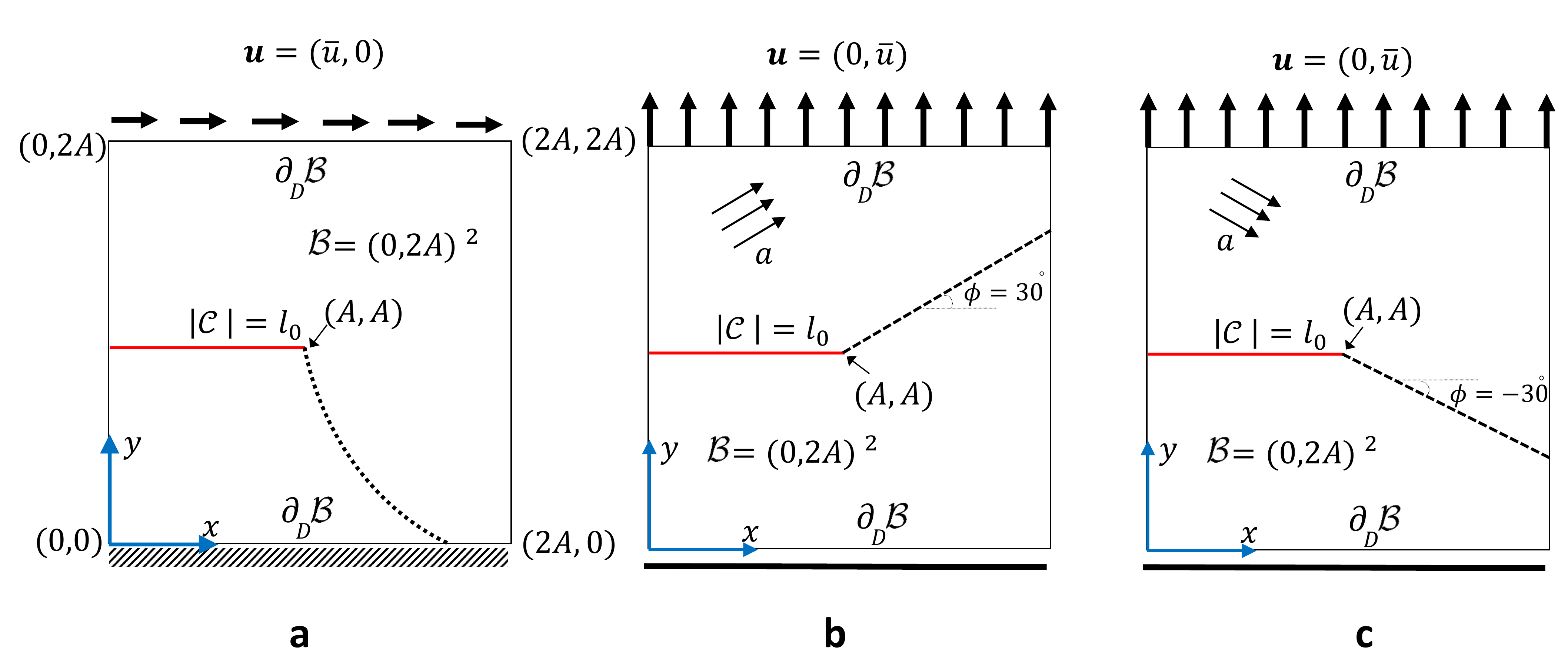}}  
	\caption{Geometry and loading setup for the single-edge-notched shear test in (a) and tensile test with the structural director $\bm{a}$ inclined under an angle $\phi=+30$ in (b) and $\phi=-30$ in (c) with respect to the $x$-axis.
	}
	\label{B1}
\end{figure}

\sectpb[Section51]{{Example 1:} Isotropic single-edge-notched shear test}

In this example, attention is restricted to pure isotropic crack-propagation by letting ${\bm{a}}=0$. In this setting, we consider a shear test such that fracture response exhibits a curved surface. The numerical computation is performed by applying a monotonic displacement $\bar{\Bu}=5\times10^{-5}$ in  horizontal direction at the top boundary of the specimen.

An important aspect that has to be verified at the local level
is the fracture state. Thus, we look at {the} crack phase field pattern
at the local scale to investigate the transition of external loading
increments (i.e. {the} down-scaling procedure) from the global domain to the local level. The Global-Local adaptive scheme to capture the curved surface is evaluated for different  $\texttt{TOL}_{d}$. Hereby, $\texttt{TOL}_{d}$ leads to different fracture zones and hence different local domains. Four different values of $\texttt{TOL}_{d}:=(0.9,0.85,0.80,0.70)$ are considered. Fig. \ref{B2} shows the evolution of the local domains for different $\texttt{TOL}_{d}$. The global mesh is only used to show a clear representation for the evolution of the local domain. Since global and local domains are performed independently, and therefore we deal with a two-scale finite element algorithm.

\begin{figure}[!ht]
	\centering
	{\includegraphics[clip,trim=1cm 0cm 3cm 0cm, width=16cm]{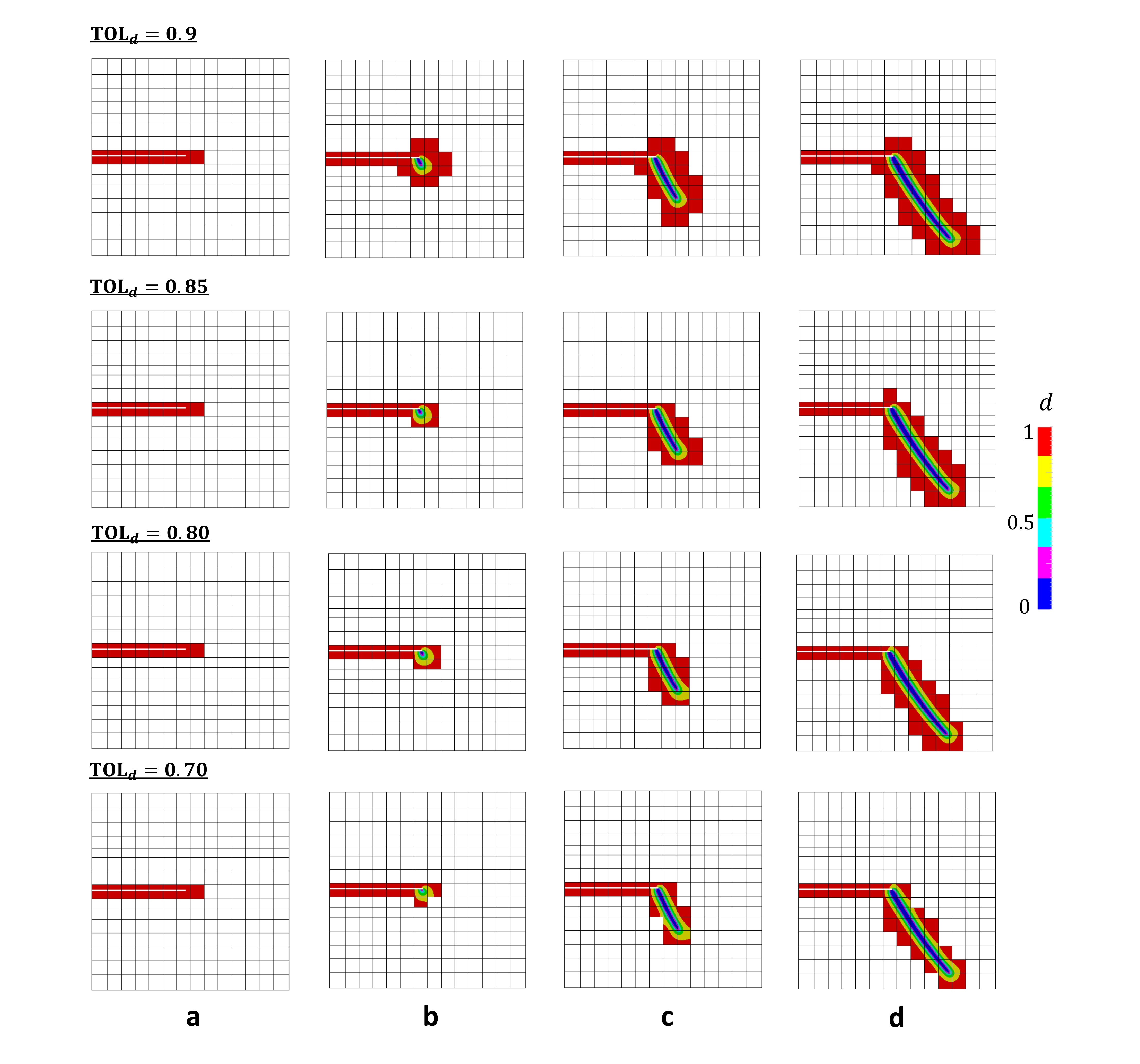}}  
	\caption{Example 1. Evolution of the local domain for different threshold values of $\texttt{TOL}_{d}$ and four deformation states up to final failure, (a) $\bar{\bm u}= 0.0058\;mm$ (b) $\bar{\bm u}= 0.0103\;mm$ (c) $\bar{\bm u}= 0.0117\;mm$ (d) $\bar{\bm u}= 0.0180\;mm$.
	}
	\label{B2}
\end{figure}

By comparing e.g. the first and fourth row in Fig. \ref{B2}, it is trivial that a smaller value of the $\texttt{TOL}_{d}$ yields a narrow fracture zone. Hence, if $\texttt{TOL}_{d_1}<\texttt{TOL}_{d_2}$ then $|\calB_{L,1}|<|\calB_{L,2}|$. The resultant narrow local domain due to the adaptivity approach in Fig. \ref{B2} demonstrates the great efficiency of the proposed method. 

The above observation plays an important role by constraining the diffusivity zone of the crack phase-field in a narrow fracture region. Whereas in the standard single scale phase field modeling, the fracture zone is spread over more areas and hence a wider diffusive zone. Thus, by the Global-Local approach we limit the effect of diffusivity on the local level and not in the entire domain. 

The influence of local effects (non-linear constitutive responses) on the global scale are described based on the load-displacement response, depicted in Fig. \ref{B4}b.
\begin{figure}[!ht]
	\centering
	{\includegraphics[clip,trim=1cm 12cm 0cm 13cm, width=16cm]{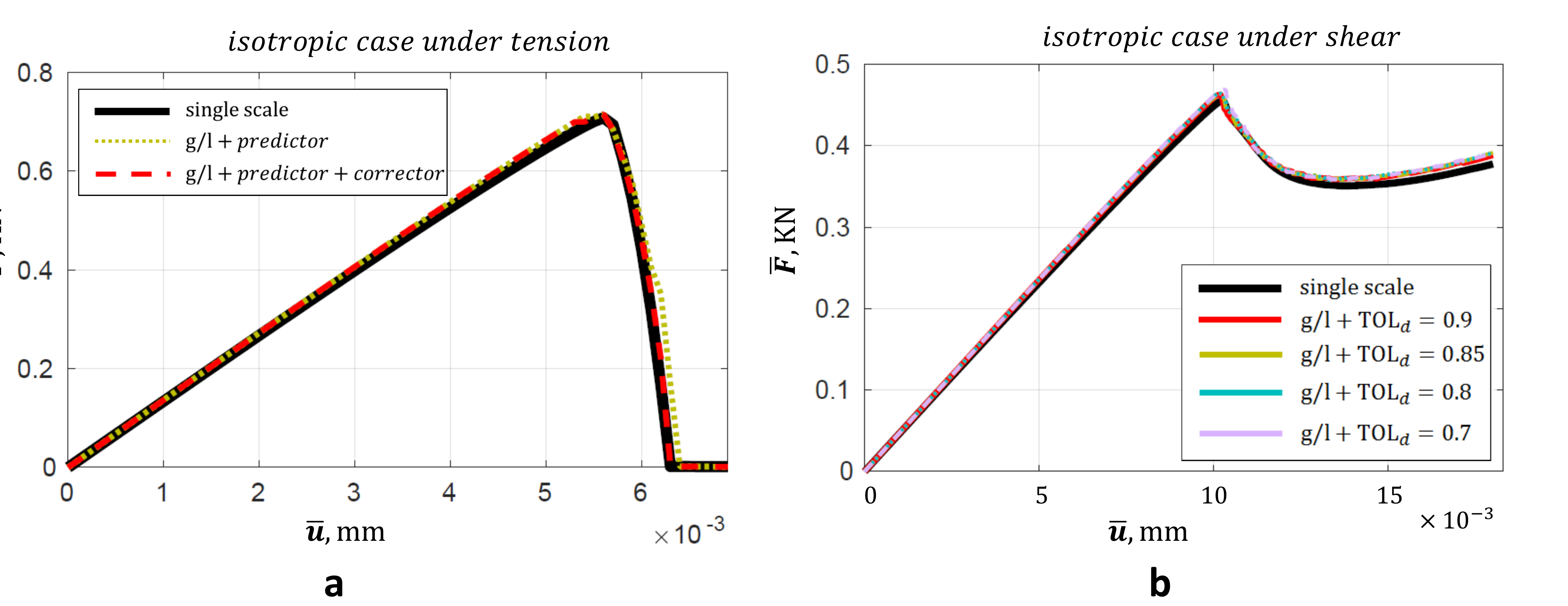}}  
	\caption{Load-displacement curves for the isotropic single-edge-notched test, (a) Example in  Section \ref{Section4}, specimen under tension test and (b) Example 1, shear test based on Global-Local approach with different $\texttt{TOL}_{d}$ verses single scale.}
	\label{B4}
\end{figure}
These curves are in a good agreement with the single scale solution. As it is excepted the Global-Local approach with a higher value of the $\texttt{TOL}_{d}=0.90$ is in very good agreement with the single scale solution. That is mainly because within single scale phase-field modeling, we deal with a wider diffusivity zone and hence more elements with fracture state are involved. 

The use of adaptivity leads to a narrow fracture zone and hence to a reduction of
degrees of freedom. That is shown in Fig. \ref{B5}a
for different values of $\texttt{TOL}_{d}$. It turns out that a smaller values of $\texttt{TOL}_{d}$ lead to a reduction of the active degrees of freedom and the computational time. At every jump which appears in
Fig. \ref{B5}a the predictor-corrector adaptive scheme is applied to the Global-Local scheme hence the number of degrees of freedom is increased.

More specifically, we show that the adaptive scheme applied to the Global-Local approach {\it considerably} reduces the computational cost in comparison with the single scale solution, as indicated in Fig. \ref{B5}b. Note that, at load step $\bar{\bm u}=10\times 10^{-3}$ (the loading point where the fracture initiates) a higher computational time is observed, see Fig. \ref{B5}. That is due to the alternate minimization approach used for solving the variational phase-field formulation which needs more iteration at cracking to reach the equilibrium state.

\begin{figure}[!ht]
	\centering
	{\includegraphics[clip,trim=1cm 12cm 0cm 14cm, width=16cm]{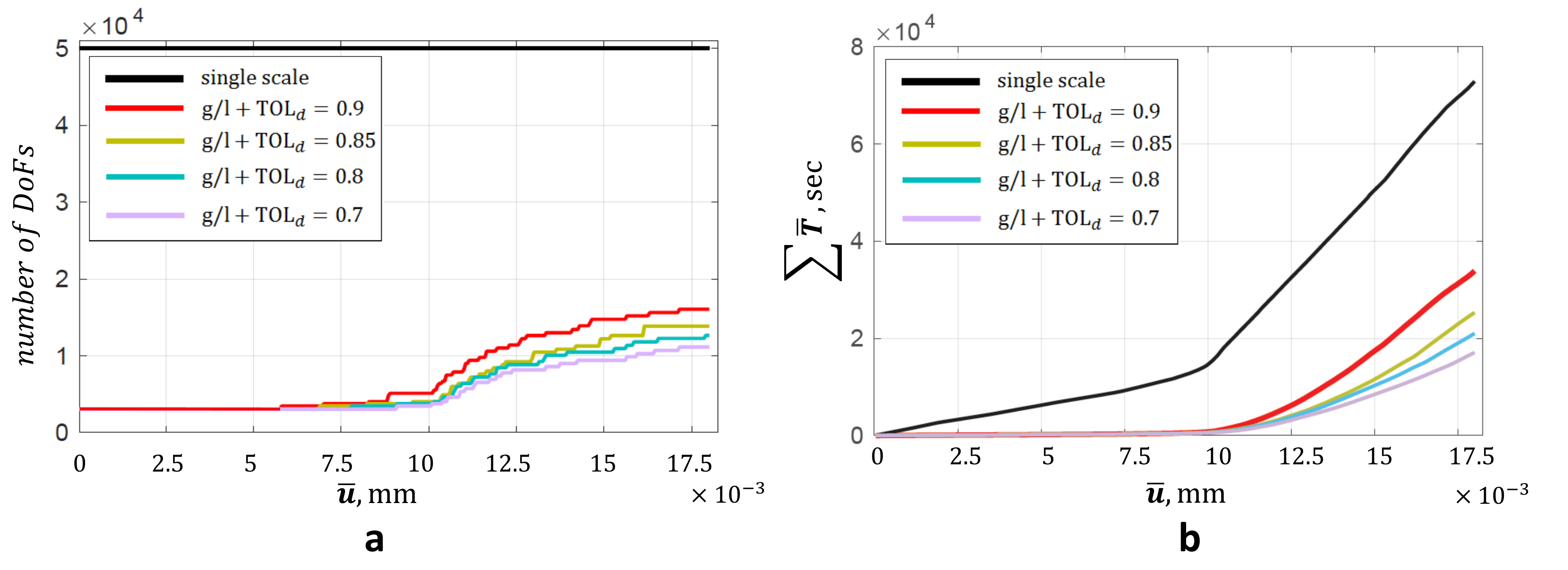}}  
	\caption{Example 1. Isotropic single-edge-notched shear test. (a) Number of degrees of freedom which have to solved for and (b) Time-displacement curves in terms of the accumulated time.
	}
	\label{B5}
\end{figure}

We now aim to investigate the energy response when solving a problem as single scale problem and as Global-Local. Recall, the consistency of the energy functional (departing point of the Global-Local approximation)
\begin{equation*}
	{\mathcal E}(\bm u,d;\BM) \equiv \widetilde{\mathcal E}(\bm u_G,\bm u_L,d_L,\bm u_\Gamma,\bm\lambda_C,\bm\lambda_L;\BM),
\end{equation*}
between the single scale and the Global-Local functional indicated in Formulation \ref{form_2} and (\ref{EGL}), respectively. We investigate this approximation by means of the evolution of the total stored elastic strain energy plotted in Fig. \ref{B6}a and the dissipated fracture energy shown in Fig. \ref{B6}b during load increments. These Global-Local simulation results show very good agreement with the single scale scheme yet with its efficiency in time shown in Fig. \ref{B5}b. 
\begin{figure}[!ht]
	\centering
	{\includegraphics[clip,trim=1cm 14cm 0cm 12cm, width=16cm]{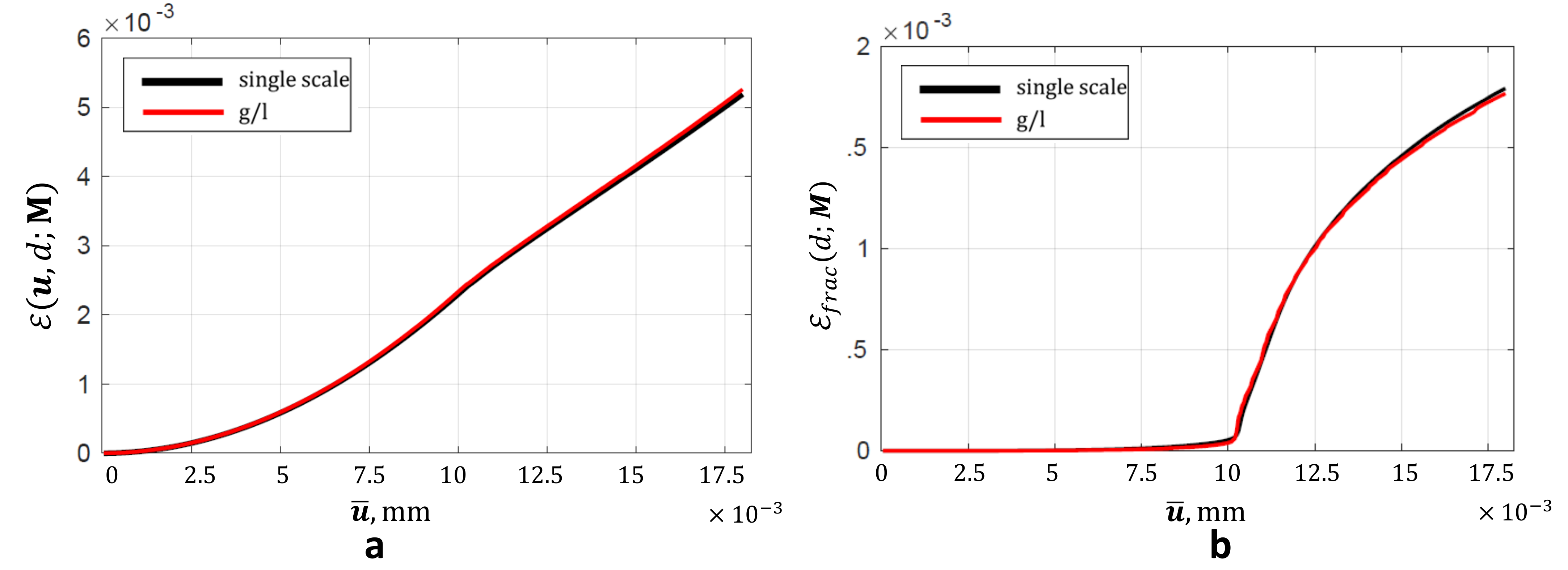}}  
	\caption{Example 1. Comparison of the energy response between the single scale domain and the Global-Local scheme for (a) the total free-energy functional and (b) the dissipated fracture energy.
	}
	\label{B6}
\end{figure}

At the converged Global-Local state, we obtain the following updated fields:\\ $(\bm u_{G,n}, \bm u_{L,n}$ $,d_{L,n}, \bm u_{\Gamma,n}, \bm\lambda_{C,n}, \bm\lambda_{L,n})$. Based on that, the homogenized global phase-field solution at the complete failure state is illustrated in Fig. \ref{B3} for different values of $\texttt{TOL}_{d}$. We emphasis that in the global level there is no pre-defined crack (i.e. globally there is no notch). However, it is interesting to note that the homogenized global phase-field solution is able to capture the crack direction which is indeed a consistent projection of the local response. Figure \ref{B3}a provides the global phase-field solution by means of (\ref{global_crack1}) at the converged Global-Local state. Accordingly, Fig. \ref{B3}b provides the homogeneous phase-field solution based on (\ref{homo_phf_sol}). Note that, the homogenized global phase-field is slightly affected by $\texttt{TOL}_{d}$.
\begin{figure}[!ht]
	\centering
	{\includegraphics[clip,trim=1cm 23cm 0cm 1cm, width=16cm]{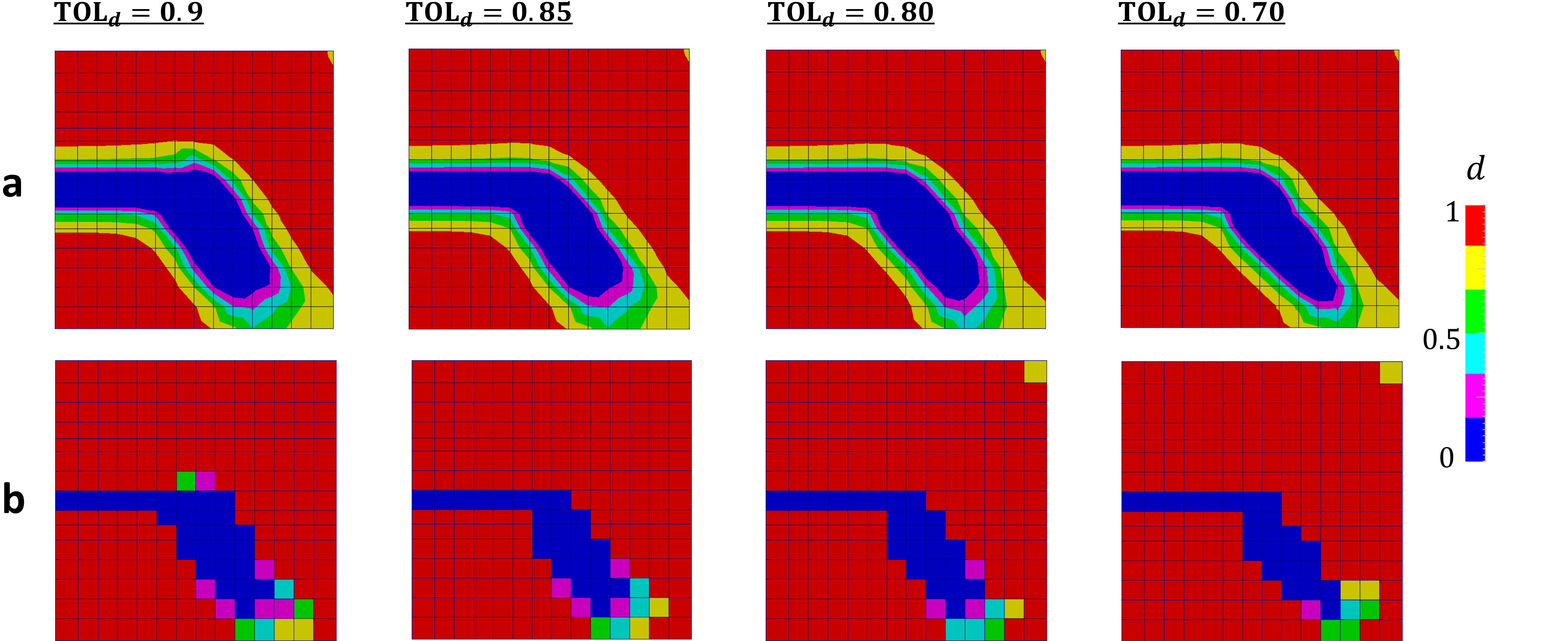}}  
	\caption{Example 1. Homogenized global crack phase-field solution for the isotropic single-edge-notched shear test (a) Global crack phase-field (b) Homogeneous crack phase-field solution at $\bar{\bm u}= 0.0180\;mm$ for different $\texttt{TOL}_{d}$.  
	}
	\label{B3}
\end{figure}

\sectpb[Section52]{{Example 2:} Analysis of transversely isotropic single-edge-notched tension test}

The second example deals with transversely isotropic material responses under
tension. It is based on different fiber directions given by the structural director
$\bm{a}$ which is inclined under an angle $\phi$ with respect to the $x$-axis
of a fixed Cartesian coordinate system. The numerical simulation is performed by applying a monotonic displacement $\bar{\Bu}=5\times10^{-5}$ in vertical direction at the top of the specimen with a linearly increasing
displacement. This loading setting is applied to the rest of numerical examples.

\sectpc[Section521]{Fiber direction of $\phi=+30^\circ$}
Here we investigate the transversely isotropic single-edge-notched tension test based on the fiber direction angle $\phi=+30^\circ$. We apply the Global-Local approach as follows:

\begin{figure}[!ht]
	\centering
	{\includegraphics[clip,trim=1cm 9cm 0.8cm 8.5cm, width=14.5cm]{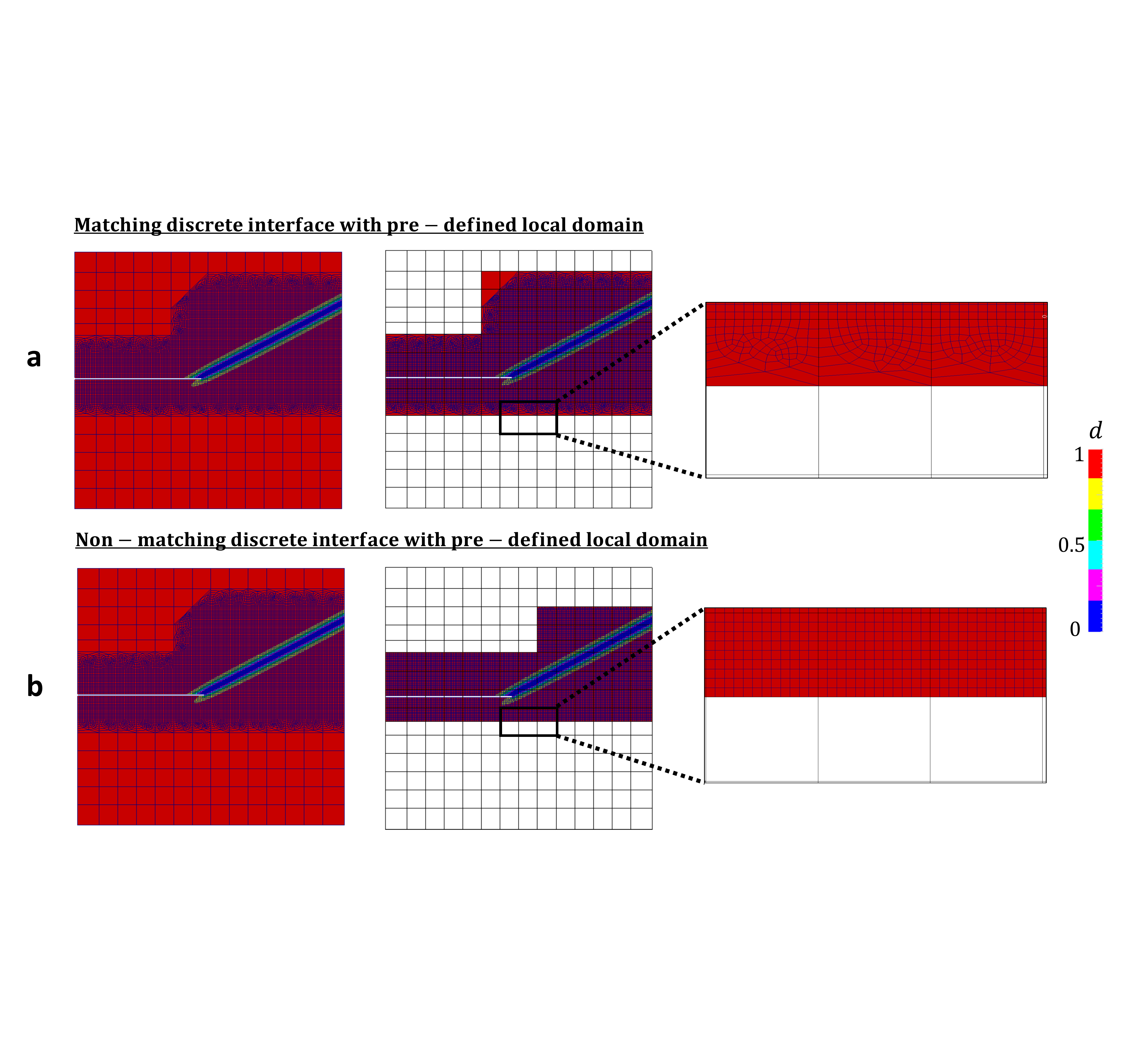}}  
	\caption{Example 2 (\ref{Section521}). Transversely isotropic single-edge-notched tension with $\phi=+30^\circ$. Complete fracture state at $\bar{\bm u}= 0.010\;mm$ for (a) matching discrete interface (\textbf{Case a}) $\Gamma_G=\Gamma_L$ and (b) non-matching discrete interface $\Gamma_G\neq\Gamma_L$ with predefined local domain (\textbf{Case b}).
	}
	\label{C1}
\end{figure}

\begin{itemize}
	\item \textbf{Case a.} {\textit{Without}} non-matching discrete interface and \textit{without} adaptive scheme (a pre-defined local domain),	
\end{itemize}

In this case, see Fig. \ref{C1}a, we aim to evaluate the
Proposition \ref{theorem1} such that
$\calB=\calB_C\cup\Gamma\cup\calB_L$. Here, the local domain is
predefined {and} no adaptivity is applied. Also, the discrete interface
between the global and local domains are one to one such that $\Gamma_G=\Gamma_L$
(see Fig. \ref{C1}a last column). It is clear that the two finite 
element meshes used for Global-Local approach represent "precisely"
the same as a single scale domain. Hence, we expect an identical
Global-Local response compared with the single scale solutions (see Remark 1 in
Appendix B). The complete fracture state is shown in Fig. \ref{C1}a. 
Accordingly, {a} comparison of the load-displacement curves of
the proposed formulation {is} demonstrated in Fig. \ref{C2}a {and}
shows a very good agreement compared with the single scale problem.

\begin{Remark}
Note that this case is similar to the work of \citet{NoiiGL18} if Robin-type boundary conditions are not taken into account (i.e. Global-Local approach based on Dirichlet-Neumann type boundary conditions), as sketched in Fig. \ref{Fig4}a. Therein, the corresponding cumulative computational time is higher compared with the reference single scale solution (see Figure 10 in \cite{NoiiGL18}), due to the slow convergence of the Global-Local procedure. That motivated the introduction of the Robin-type boundary conditions, resulting in a reduction of the computational cost, see Fig. \ref{B5}b.
\end{Remark}

\begin{itemize}	
	\item \textbf{Case b.} {\textit{With}} non-matching discrete interface and \textit{without} adaptive scheme (a pre-defined local domain).	
\end{itemize}

In the second case, we assume $\calB\neq\calB_C\cup\Gamma\cup\calB_L$  and that the local domain is pre-defined hence no adaptivity is applied. Furthermore, the discrete interface between global and local domains are non-matching such that $\Gamma_G\neq\Gamma_L$ (see Fig. \ref{C1}b last column). This removes one restriction applied in \textbf{Case a}, that is the matching discrete interface criteria.

The complete fracture state is shown in Fig. \ref{C1}b. Compared with the first case, by the non-matching discrete interface, we are able to have an arbitrary mesh at the local domain (including interface) without any given interface conditions (and to avoid having distorted mesh between fine and coarse discretizations). The interface conditions refer to the identical discretization space for $\Gamma_L$ and $\Gamma_G$, see Remark \ref{non_match}.  The importance can be observed when the fracture reaches the interface, see e.g. Fig. \ref{B2}.

The resulting load-displacement curve in Fig. \ref{C2}b has a very good
agreement when compared with {the} single scale {problem.}

\begin{figure}[!ht]
	\centering
	{\includegraphics[clip,trim=1cm 11.5cm 0cm 14.5cm, width=16cm]{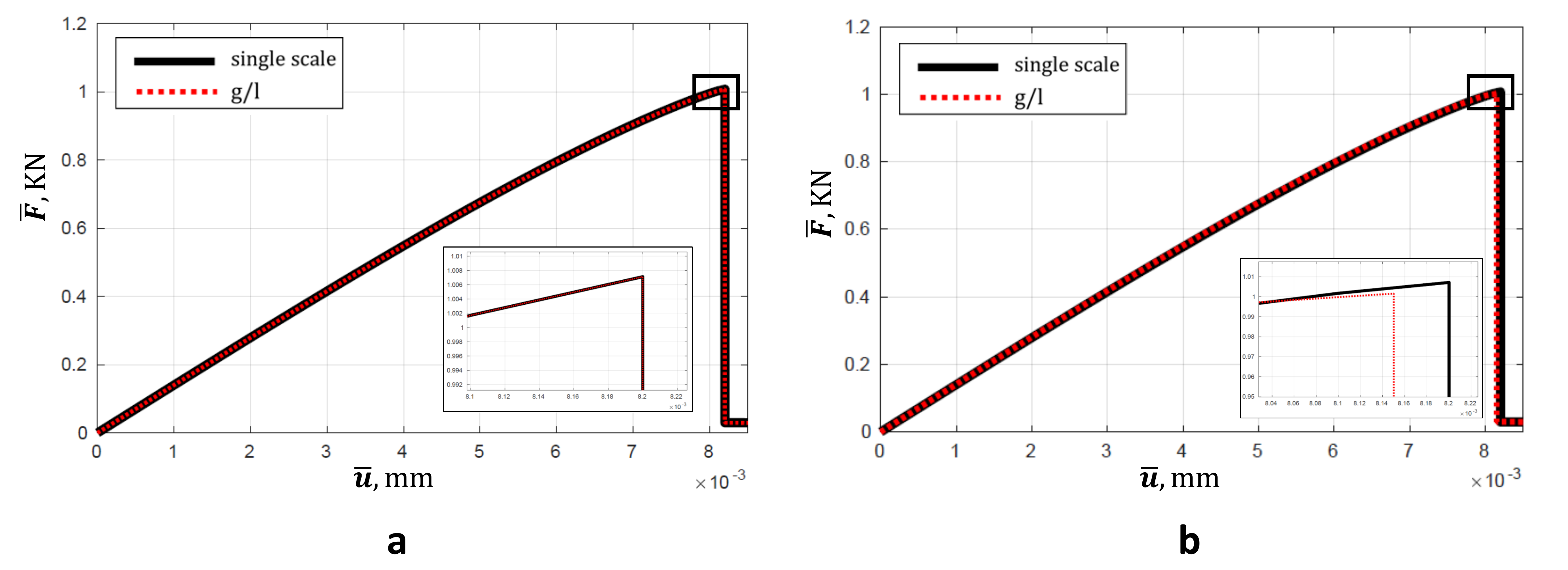}}  
	\caption{Example 2 (\ref{Section521}). Transversely isotropic single-edge-notched tension with $\phi=+30^\circ$. Comparison of the load-displacement curve between single scale problem and Global-Local formulation. (a) Matching discrete interface with $\calB=\calB_C\cup\Gamma\cup\calB_L$ (\textbf{Case a}) and (b) non-matching discrete interface with $\calB\neq\calB_C\cup\Gamma\cup\calB_L$ (\textbf{Case b}).
	}
	\label{C2}
\end{figure}

\begin{itemize}	
	\item \textbf{Case c.} {\textit{With}} non-matching discrete interface and \textit{with} adaptive scheme.
\end{itemize}
In the third case, we consider a non-matching discrete interface along with an adaptive scheme. This case removes all restrictions applied in \textbf{Case a} (matching interface and predefined local domain). 

Fig. \ref{D1} illustrates the evolution of the crack phase-field along with the local domain and the corresponding Global-Local interface. The local domain and its coupling interface must be computed at each stage. The second row of Fig. \ref{D1} presents the local mesh evolution such that the non-matching discrete interface between global and local mesh is examined.
\begin{figure}[!ht]
	\centering
	{\includegraphics[clip,trim=1cm 17cm 2cm 8cm, width=16cm]{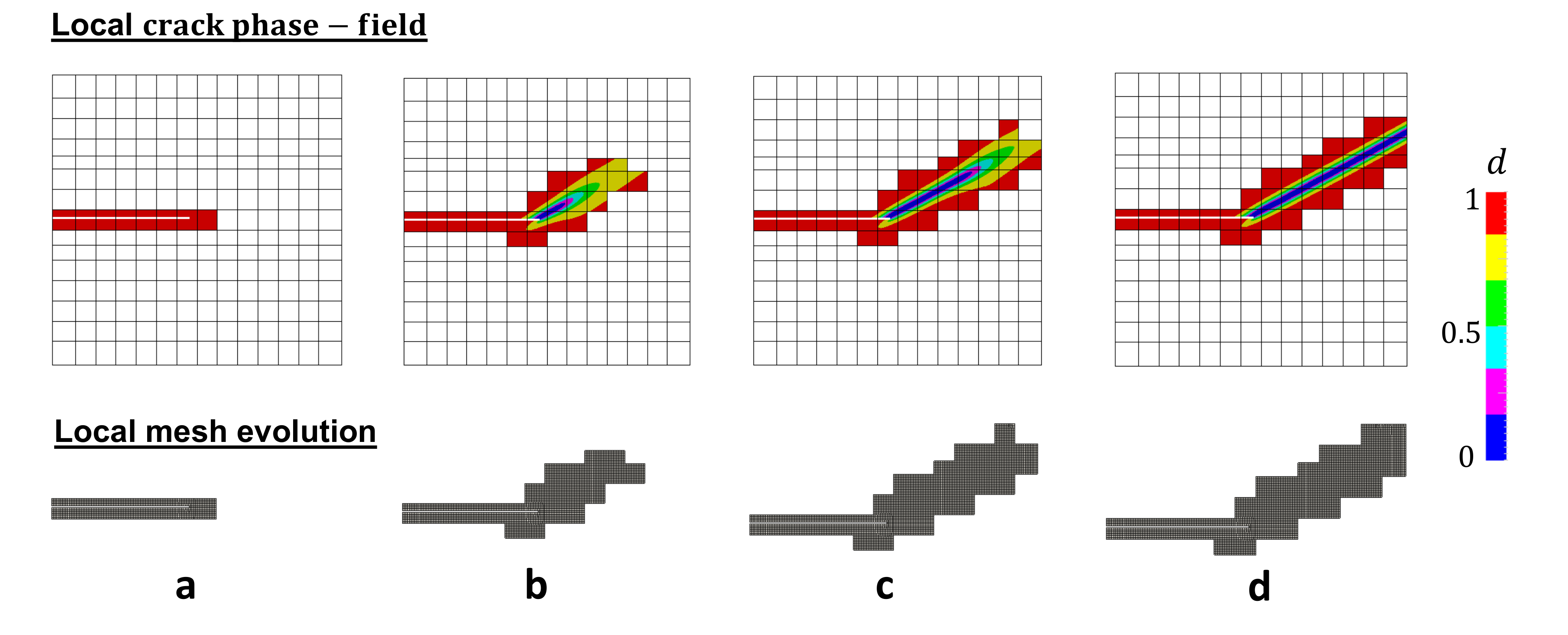}}  
	\caption{Example 2 (\ref{Section521} -- \textbf{Case c}). Transversely isotropic single-edge-notched tension with $\phi=+30^\circ$. First row indicates the local crack phase-field resolution and second row represents the evolution of the local domain per time for different deformation states as follows: (a) $\bar{\bm u}= 0.0030\;mm$ (b) $\bar{\bm u}= 0.0088\;mm$ (c) $\bar{\bm u}= 0.0092\;mm$ and (d) $\bar{\bm u}= 0.010\;mm$.}
	\label{D1}
\end{figure}
To evaluate the solution related to the local to global transition, the load-displacement curve is shown in Fig. \ref{D3}a. 
Since the single scale problem produces a very diffusive transition zone for the phase-field, more elements are involved (this is not the case in the sharp crack limit). This results in a small difference in the load-displacement curves between the  Global-Local formulations and the reference single scale. Fig. \ref{D3}b illustrates a reduction of the number of degrees of freedom.

\begin{figure}[!ht]
	\centering
	{\includegraphics[clip,trim=1cm 11cm 0cm 16cm, width=16cm]{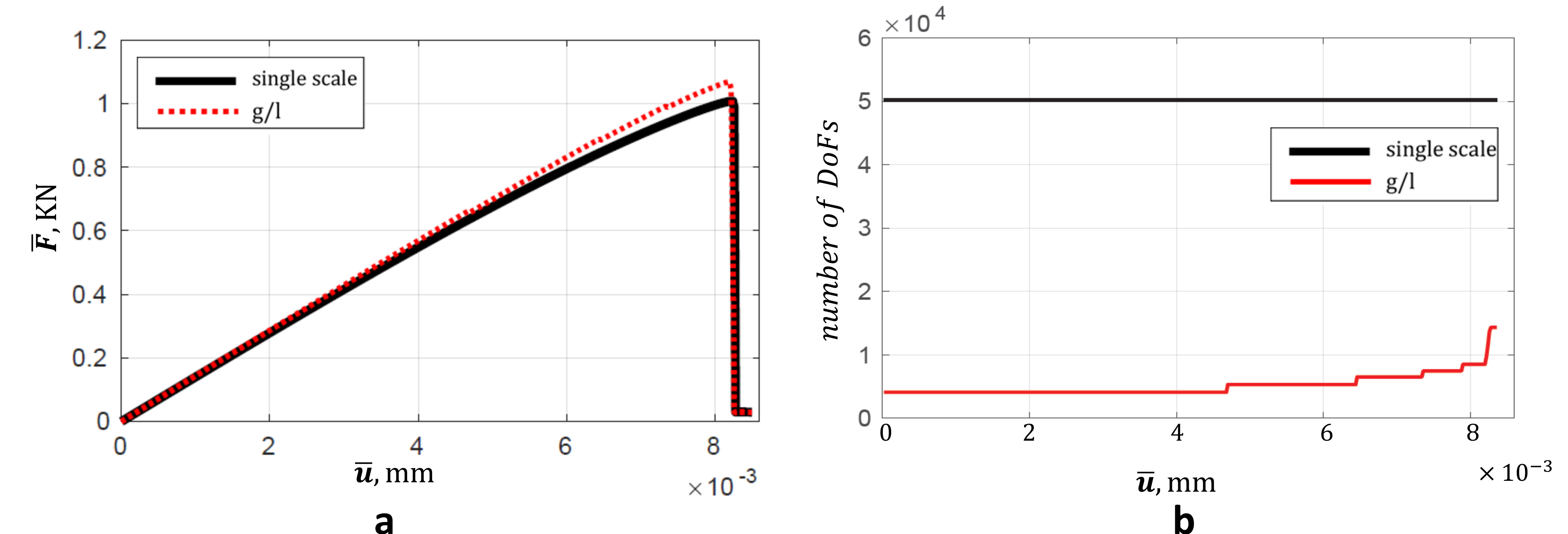}}  
	\caption{Example 2 (\ref{Section521} -- \textbf{Case c}). Transversely isotropic single-edge-notched tension with $\phi=+30^\circ$. (a) Comparison of the load-displacement curve; (b) Number of degrees of freedom.
	}
	\label{D3}
\end{figure}

The Global-Local approach, besides its feasibility for having two ad-hoc finite
element models for the global and local domain, enables computations with legacy codes. Additionally, the reduction of unknowns leads to a reduction of
the computational time. To illustrate the time efficiency, the simulation time
ratio between single scale and the Global-Local approach are shown in
Fig. \ref{D4}a. 
It {can be observed that in average, the} Global-Local formulations 
{perform} $12$ times faster. Furthermore, Fig. \ref{D4}b demonstrates the corresponding accumulative computational time, which underlines the efficiency of the predictor-corrector adaptive scheme.

\begin{figure}[!ht]
	\centering
	{\includegraphics[clip,trim=0cm 13.5cm 0cm 12cm, width=16cm]{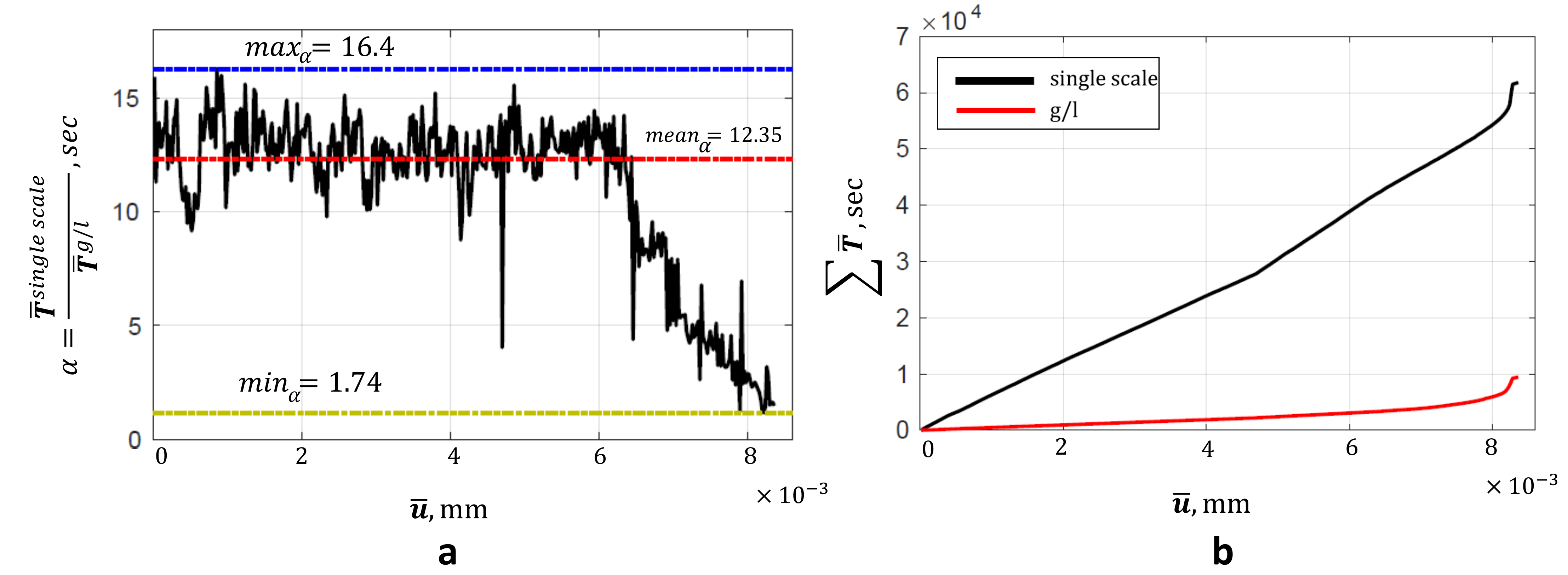}}  
	\caption{Example 2 (\ref{Section521} -- \textbf{Case c}). Transversely isotropic single-edge-notched tension with $\phi=+30^\circ$. (a) Time ratio between the computed single scale and the Global-Local time per loading steps and (b) accumulated time-displacement curves.
	}
	\label{D4}
\end{figure}

Fig. \ref{D5}a presents the total elastic strain energy per load increments. The resulting Global-Local curve is in a very good agreement with the single scale approach.

\begin{figure}[!ht]
	\centering
	{\includegraphics[clip,trim=1cm 11cm 0cm 14.5cm, width=16cm]{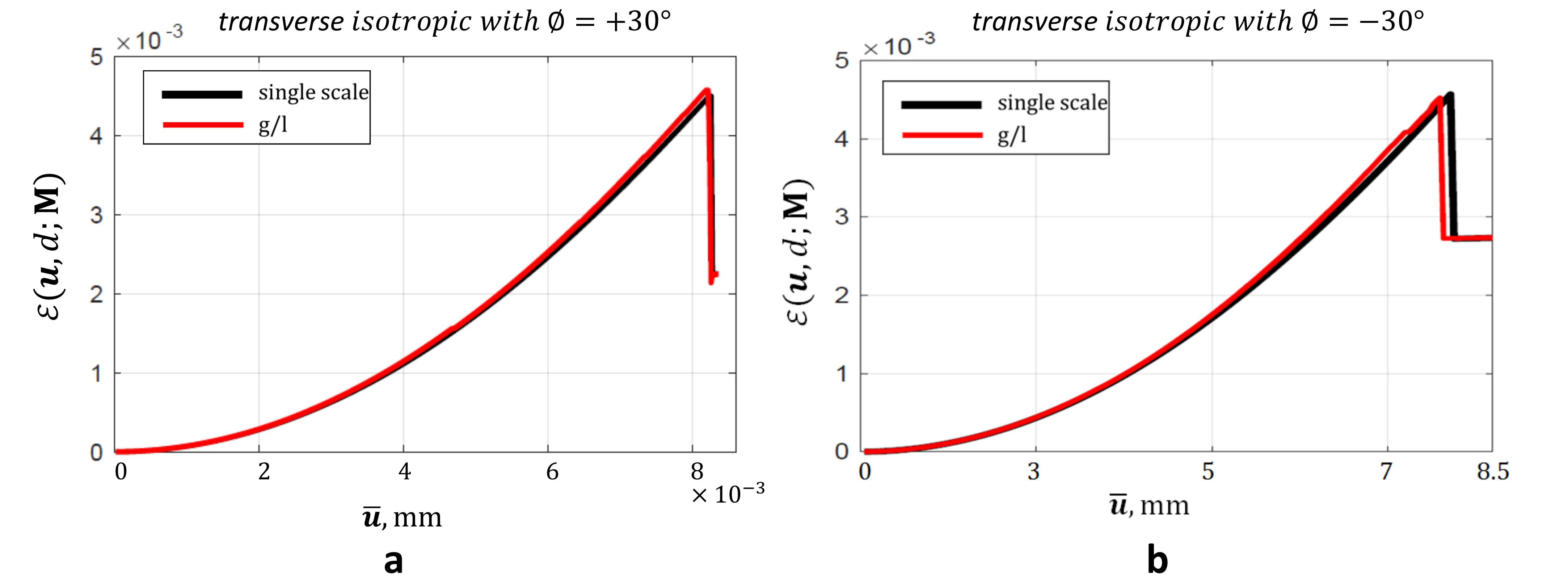}}  
	\caption{Example 2 (\ref{Section521} -- \textbf{Case c}). Transversely isotropic single-edge-notched tension test. Comparison of the total energy functional between single scale domain and Global-Local scheme for fiber direction with $\phi=+30^\circ$ in (a) and fiber direction with $\phi=-30^\circ$ in (b).}
	\label{D5}
\end{figure}
Accordingly, the homogenized global phase-field solutions for different fracture states are depicted in Fig. \ref{D2}. First and second row of Fig. \ref{D2} are based on approach (\textit{a}) and (\textit{b}) outlined in section \ref{Section4_homoPhF}. For comparison purposes, the single scale resolution is also plotted in the third row of Fig. \ref{D2}. It is observed that, the homogenized phase-field solution in the case of the anisotropic setting, is able to capture the crack direction. The global phase-field solution is affected by the global element size and also $\texttt{TOL}_{d}$ (local scale).
\begin{figure}[!ht]
	\centering
	{\includegraphics[clip,trim=1cm 7cm 7cm 3cm, width=14cm]{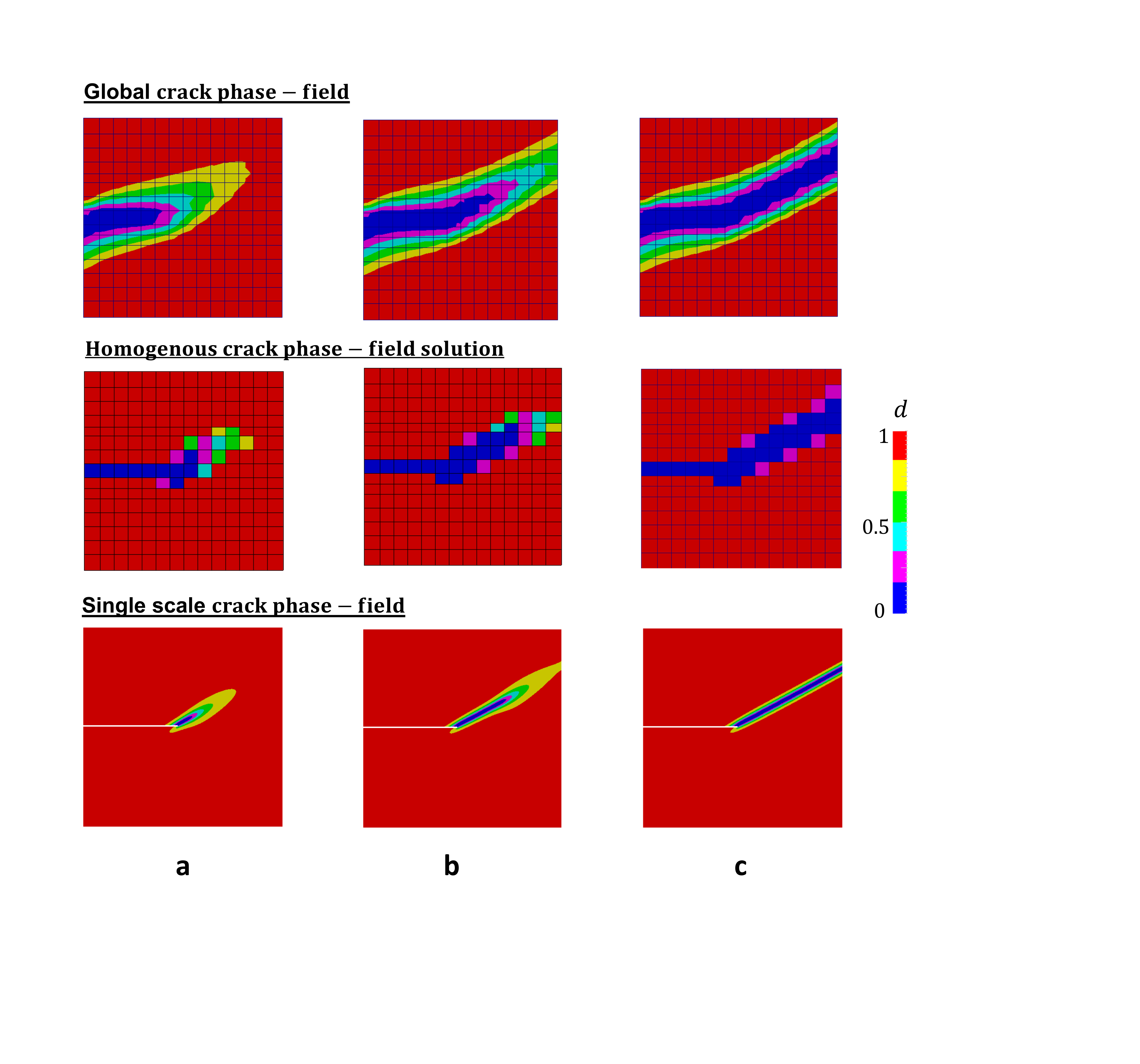}}  
	\caption{Example 2 (\ref{Section521} -- \textbf{Case c}). Homogenized global crack phase-field solution for $\phi=+30$. First row: the global crack phase-field; Second row: the homogeneous solution; Third row: single scale phase field solution per loading step for different deformations states as: (a) $\bar{\bm u}= 0.0088\;mm$ (b) $\bar{\bm u}= 0.0092\;mm$ and (c) $\bar{\bm u}= 0.010\;mm$.
	}
	\label{D2}
\end{figure}

\sectpc[Section522]{Fiber direction of $\phi=-30^\circ$}
This numerical example illustrates the transversely isotropic single-edge-notched tension test with $\phi=-30^\circ$. The crack phase-field resolution has a brutal fracture response in which a complete failure happens in one load increment. Thus the post-peak behavior is almost vertical, see Fig. \ref{E3}. The aim of this numerical example is to show the capability of the Global-Local approach to capture such a brutal fracture behavior. This is mainly possible due to the introduction of the corrector step in the adaptive scheme described in Section \ref{Section4}. 

Next, we investigate the effect of the $\texttt{TOL}_{d}$ in the case of the
brutal fracture behavior, by setting $\texttt{TOL}_{d}=(0.90, 0.80)$, as
illustrated in Fig. \ref{E1}. Similar as before, different $\texttt{TOL}_{d}$
lead to different fracture zones and hence different local domains. The crack
paths for both $\texttt{TOL}_{d}=0.90$ and $\texttt{TOL}_{d}=0.80$ are
identical, yet with $\texttt{TOL}_{d}=0.80$ a more narrow fracture
zone is observed (hence {a} reduction of computational time). 

\begin{figure}[!ht]
	\centering
	{\includegraphics[clip,trim=1cm 19cm 3cm 1cm, width=15cm]{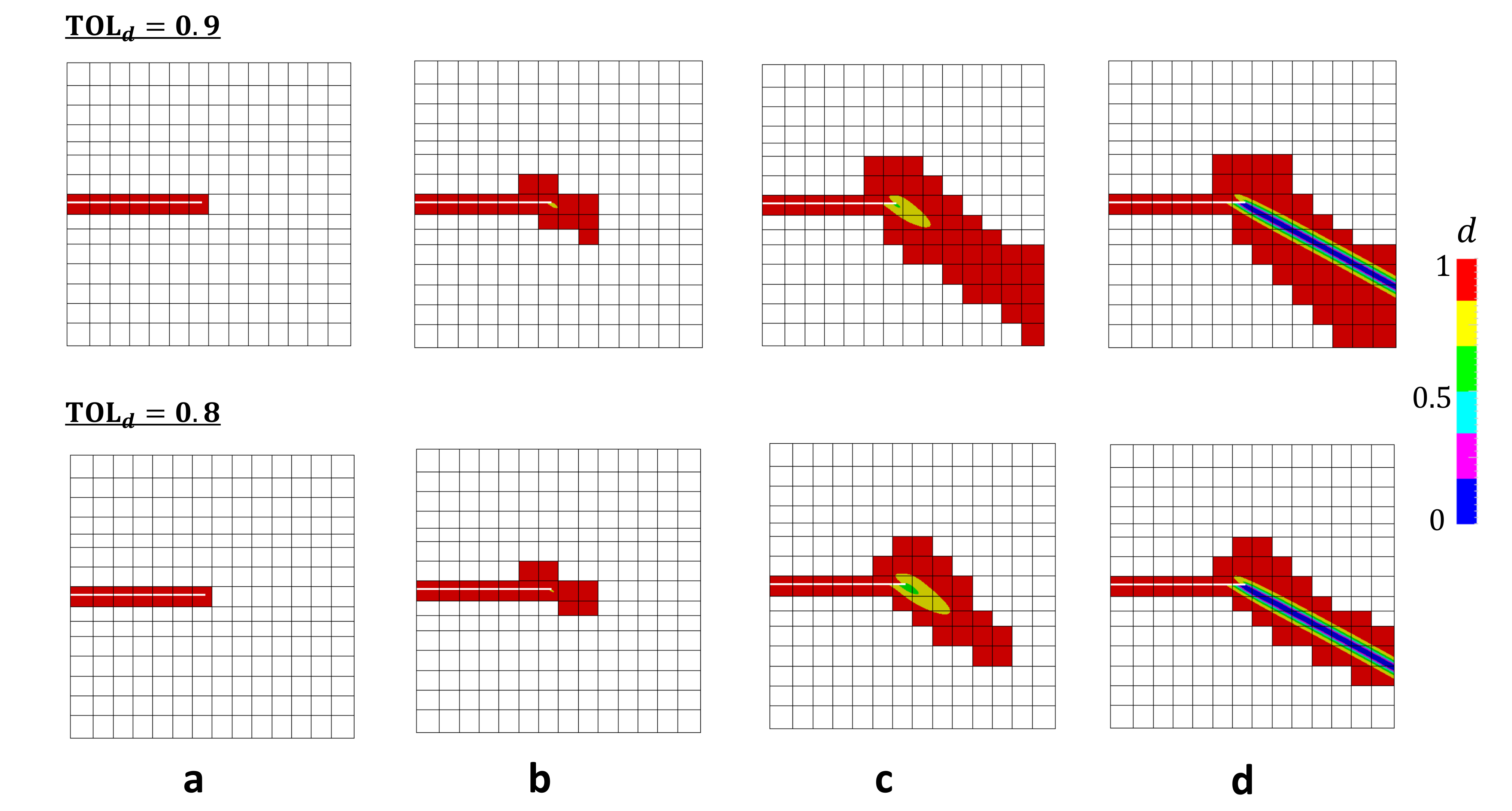}}  
	\caption{Example 2 (\ref{Section522}). Fracture of transversely isotropic single-edge-notched plate
		under tension for $\phi=-30^\circ$. Resulting, local crack phase-field by Global-Local adaptive scheme indicated in; First row with $\texttt{TOL}_{d}=0.90$; Second row with $\texttt{TOL}_{d}=0.80$ per loading steps at (a) $\bar{\bm u}= 0.0030\;mm$ (b) $\bar{\bm u}= 0.0062\;mm$ (c) $\bar{\bm u}= 0.0078\;mm$ (d) $\bar{\bm u}= 0.0085\;mm$.}
	\label{E1}
\end{figure}

{A} comparison of the load-displacement curves are shown in Fig. \ref{E3}a. 
The effect of $\texttt{TOL}_{d}$ on the load-displacement curve with zoom-in to the framed region of the left plot is shown in Fig. \ref{E3}b. Despite of its brutal fracture behavior, the load-displacement curve with the higher value of $\texttt{TOL}_{d}=0.90$ has a good agreement when compared with a single scale solution.
\begin{figure}[!ht]
	\centering
	{\includegraphics[clip,trim=1cm 12cm 0cm 15cm, width=16cm]{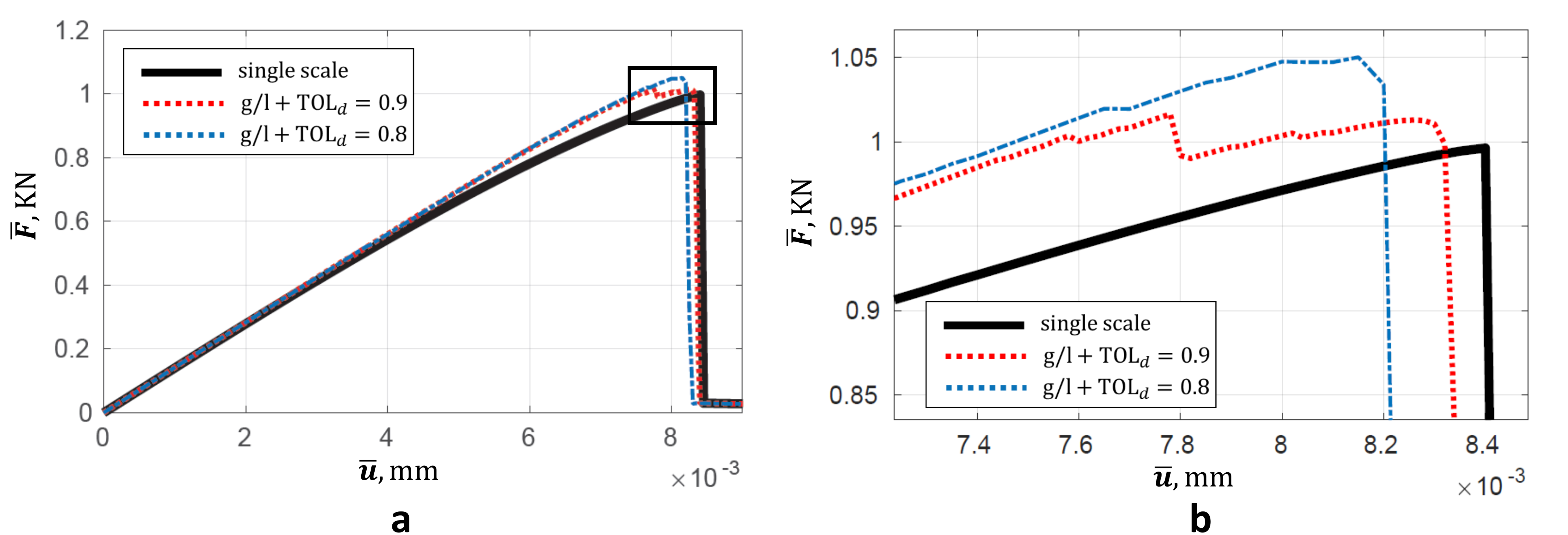}}  
	\caption{Example 2 (\ref{Section522}).  Effect of different $\texttt{TOL}_{d}$ on Global-Local approach by (a) Comparison of the load-displacement curve; (b) Zooming into the framed region of the left plot.}
	\label{E3}
\end{figure}
Following the approximation between the single scale and the Global-Local modeling, Fig. \ref{D5}b presents a very good agreement in the total elastic strain energy of both schemes during the load increments.


Figure \ref{E5} describes the efficiency of the proposed Global-Local approach. 
Here the accumulative computational time is plotted in Fig. \ref{E5}a and the number of unknowns
are plotted in Fig. \ref{E5}b versus the displacement and compared with the single scale domain. At each jump in Fig. \ref{E5}b, the predictor-corrector adaptive scheme is active and applied on the Global-Local scheme which increases the number of degrees of freedoms.

\begin{figure}[!ht]
	\centering
	{\includegraphics[clip,trim=1cm 11cm 0cm 14.5cm, width=16cm]{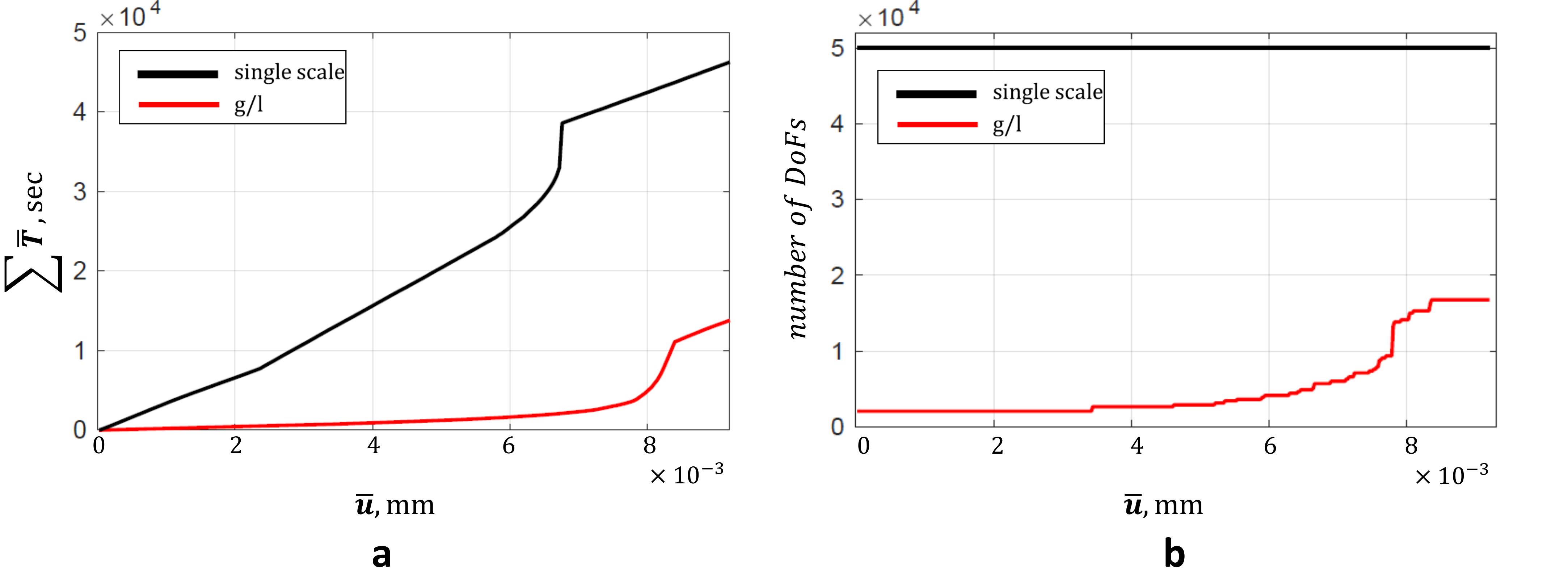}}  
	\caption{Example 2 (\ref{Section522}). Time efficiency for the transversely isotropic single-edge-notched tension test (a)  Time-displacement curves in terms of the accumulated time; (b) Number of degrees of freedom.
	}
	\label{E5}
\end{figure}

The homogenized global phase-field solution for this numerical setting is indicated in Fig. \ref{E2}a at the complete fracture state. The single scale resolution is indicated in third row of Fig. \ref{E2}b. It is evident that the homogenized phase-field solution is able to (\texttt{i}) capture the initial crack of the notch plate located at the local level and accordingly, (\texttt{ii}) the evolution of the fracture state which follows the preferred fiber direction.

\begin{figure}[!ht]
	\centering
	{\includegraphics[clip,trim=1cm 23.5cm 2cm 7cm, width=16cm]{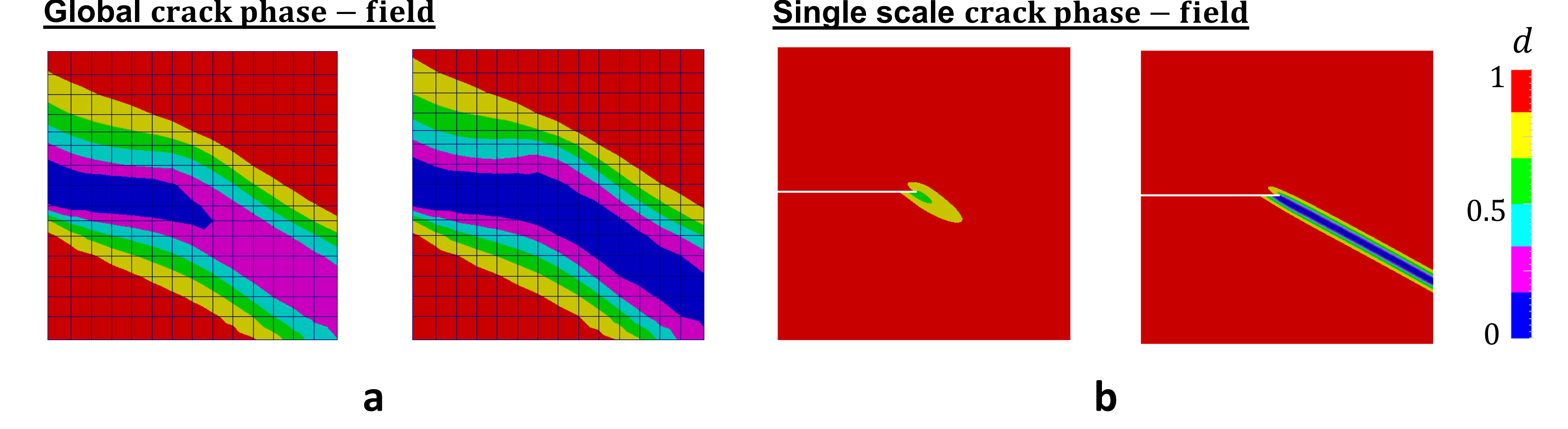}}  
	\caption{Example 2 (\ref{Section522}). Homogenized global crack phase-field solution for $\phi=-30^\circ$. (a) Global crack phase-field and (b) single scale solution for two deformation states: $\bar{\bm u}= 0.0078\;mm$ and $\bar{\bm u}= 0.0085\;mm$.
	}
	\label{E2}
\end{figure}

\sectpb[Section53]{{Example 3:} Investigation of transversely isotropic heterogeneous L-shaped panel test}

The third model problem is concerned with anisotropic brittle
fracture of 
a heterogeneous L-shaped panel test. The homogeneous isotropic counterpart 
setting for this benchmark problem has been reported by many authors, see
e.g. \cite{MESGARNEJAD2015420,UNGER20074087,ambati15,Wi17_SISC}. 
We demonstrate the performance of the Global-Local approach to predict crack propagation without any given initial crack region.
In this case the initial local domain needs to be determined based on the
critical stress state at the global level as outlined in Fig \ref{q1}c. To increase the order of complexity, a heterogeneous structure is considered by means of randomly distributed hard inclusions as plotted in Fig \ref{q1}b. Furthermore transversely isotropic material behavior is assumed. The structural director $\bm{a}$ is inclined under $\phi=-15^\circ$.

Geometry and loading conditions are depicted in Fig. \ref{q1}a. The size of the specimen is chosen to be: $A=B=250\;mm$ and $H=30\;mm$. The bottom edge of the specimen is fixed in both directions and a vertical displacement is applied until final failure, see Fig. \ref{q1}a. One third of the 
specimen is covered by hard inclusions, as shown in Fig. \ref{q1}b. Here, crack propagation is expected. The remaining parts of the domain are supposed to be homogeneous, but affected by the transversely isotropic behavior.

The material parameters used in the simulation are the same as
in \cite{UNGER20074087} and set as: $\lambda = 6.16$ kN/mm$^2$, $\mu = 10.95$
kN/mm$^2$, $G_c = 9 \times 10^{-5}$ kN/mm, $\alpha= \chi = 50$ and $\Xi =0$. 
The dimensionless mismatch ratio is denoted by 
$m=E_{inclusion}/E_{matrix}$ (here, $E$ refers to Young's modulus) and 
 set as $m=10$. 
Thus, we deal with 
$E_{matrix}=25.85$ kN/mm$^2$ and $E_{inclusion}=258.5$ kN/mm$^2$. 

\begin{figure}[!ht]
	\centering
	{\includegraphics[clip,trim=1cm 19cm 0cm 9.5cm, width=17cm]{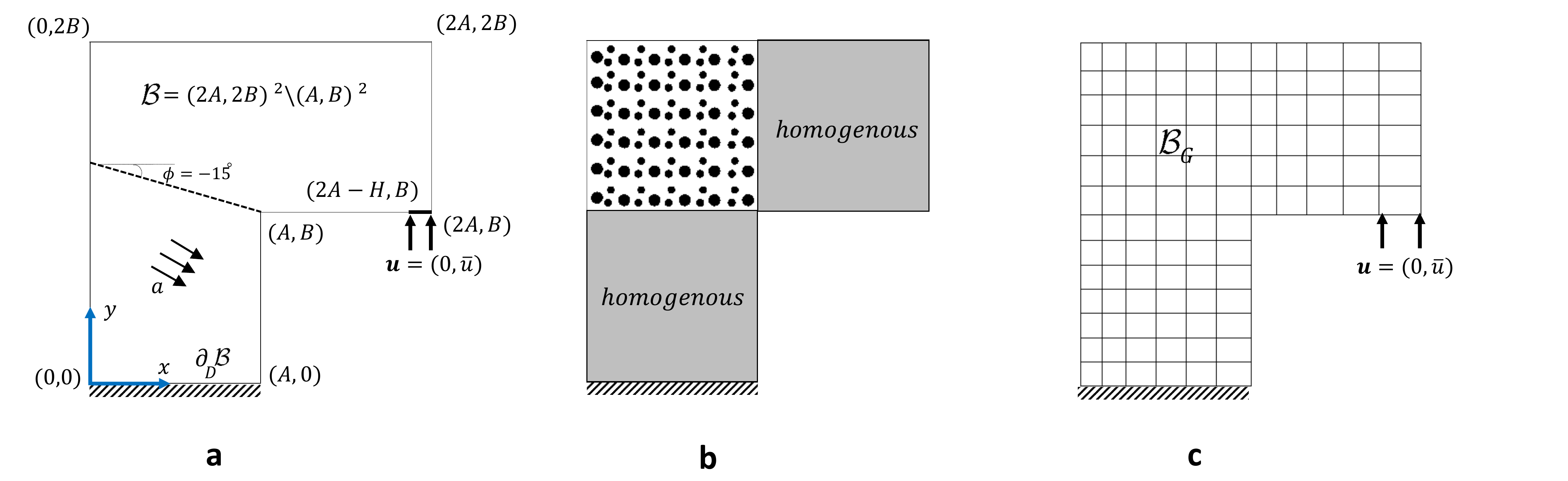}}  
	\caption{Example 3. Heterogeneous L-shaped panel test. (a) Geometry and loading setup with a structural director $\bm{a}$ inclined under an angle $\phi=-15$ (b) partitioning of domain into the heterogeneity and homogeneity counterparts and (c) global finite element mesh without potential fictitious zones.}
	\label{q1}
\end{figure}

In order to determine an initial fictitious domain, which 
has to be used for the Global-Local approach, 
an idea of the phase-field formulation with threshold state is considered. 
Here, a critical stress state on the global level is employed by extending the critical stress value of the isotropic phase-field formulation in \cite{borden+verhoosel+scott+hughes+landis12} to an anisotropic heterogeneous setting. Hence, the critical values for the stress and corresponding strain, are obtained as
\begin{equation}
{\bm{\varepsilon}}_c=\frac{\sqrt{3}}{3\,}\sqrt{\frac{G_{c}}{l \left(\mathrm{\chi}\, a^4 + \bar{E}\right)}}\quad\mbox{and}  \quad
\bm{\sigma}_c=\frac{3\, \sqrt{3}}{16}{\sqrt{\frac{G_c\, \left(\mathrm{\chi}\, a^4 + \bar{E}\right)}{l}}} .
\label{critical_stress}
\end{equation}
This result is based on (\ref{eq29_H3}) and (\ref{homo_phf_sol}) where $a = \sin(\phi)$, see \cite{Fadianiso17}. The effective Young's modulus $\bar{E}$ for the heterogeneous domain is defined as
\begin{equation*}
\bar{E}:=\frac{{E}_{Voigt}+{E}_{Reuss}}{2}
\end{equation*}
with
\begin{equation*}
{E}_{Voigt}:= \frac{1}{V} \int_{\calB} E \,\mathrm{d}{\Bx},\;\;\mbox{and}\;\;
{E}^{-1}_{Reuss}:= \frac{1}{V} \int_{\calB} E^{-1} \,\mathrm{d}{\Bx} .
\end{equation*}
For $\chi=0$ in (\ref{critical_stress}) the isotropic case of \cite{borden+verhoosel+scott+hughes+landis12} is recovered.
The critical stress state increases as $l$ decreases. Additionally, if the
length scale $l$ goes to 0 in the limit, the crack nucleation stress
tends to infinity. This is in agreement with Griffith's theory, which allows 
to have crack nucleation in  stress singularities. Then, the critical stresses based on the effective Young's modulus $\bar{E}=36.136$ kN/mm$^2$ yields $\bm{\sigma}_c=9.2603$ N/mm$^2$ with the Voigt average ${E}_{Voigt}=44.423$ kN/mm$^2$ and the Reuss average ${E}_{Reuss} = 27.850$ kN/mm$^2$.


Figure \ref{q3} shows the maximum stress state distribution of the
heterogeneous L-shaped panel test. Here, the maximum stress is observed in the
corner point of the specimen where the singularities are located. Hence, this is the potential candidate for the local domain. The Global-Local approach is then started after the stress state on the global domain reaches $75 \%$ of $\bm{\sigma}_c$. This percentage of the critical stress is chosen to be on the safe side when starting the Global-Local formulations.

\begin{figure}[!t]
	\centering
	{\includegraphics[clip,trim=7.5cm 16.5cm 2cm 14cm, width=17.5cm]{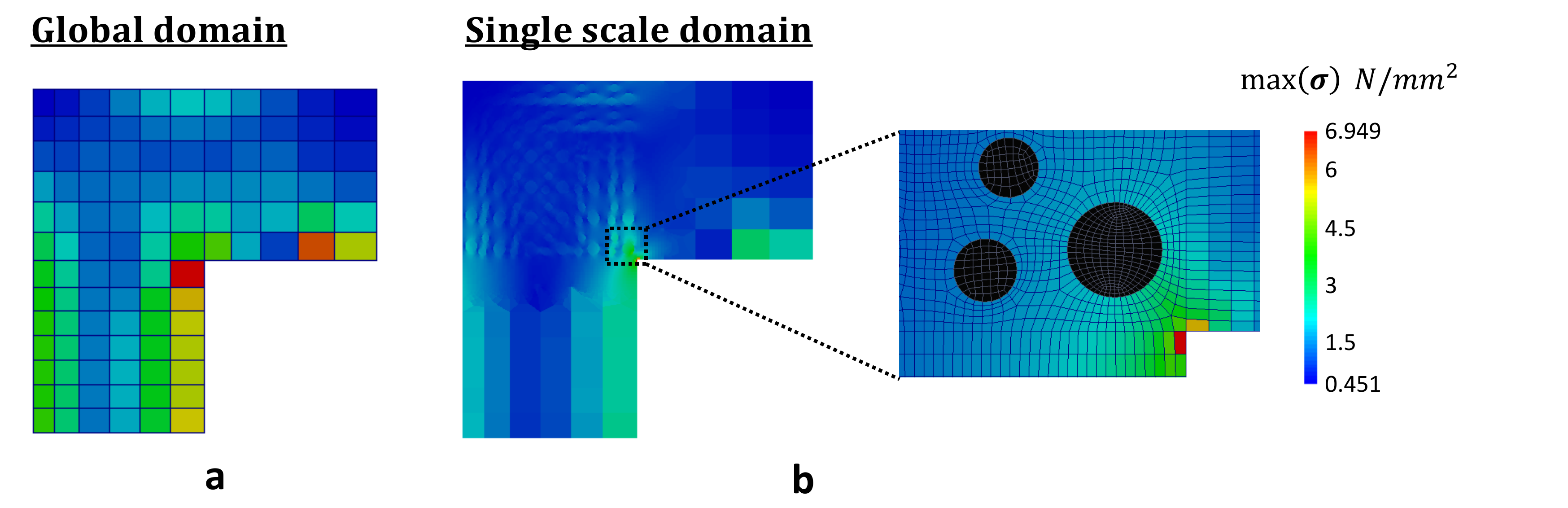}}  
	\caption{Example 3. Maximum stress state in the heterogeneous L-shaped panel test. (a) Global stress state and (b) single scale stress state.
	}
	\label{q3}
\end{figure}
\begin{figure}[!t]
	\centering
	{\includegraphics[clip,trim=0cm 6cm 0cm 2.2cm, width=17cm]{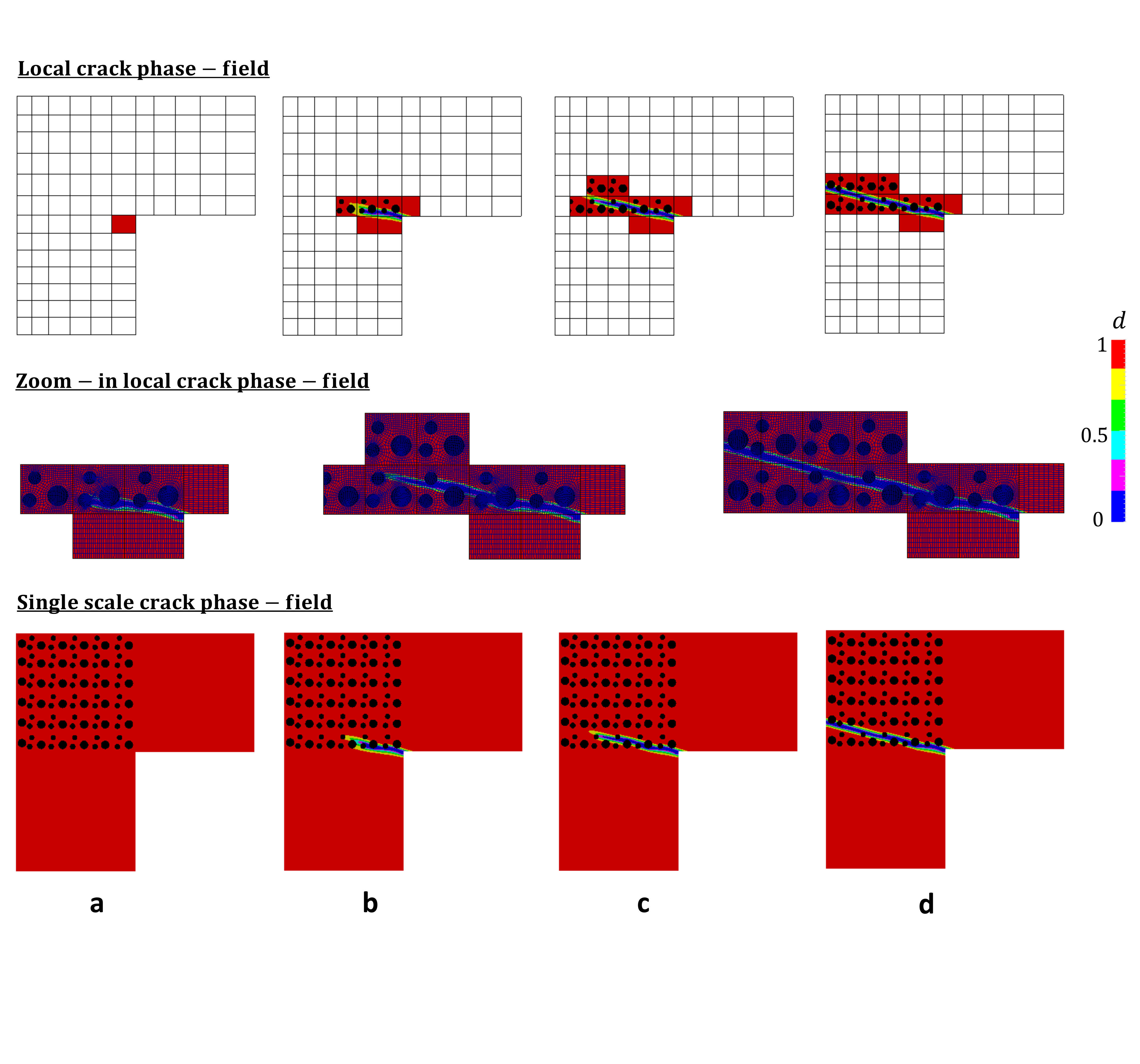}}  
	\caption{Example 3. Crack phase-field pattern for the transversely isotropic heterogeneous L-shaped panel test with fiber direction angle of $\phi=-15^\circ$.
		First row: local crack phase-field based on the adaptive scheme; Second row: mesh evolution for local domain by considering the influence of inclusions; Third row: resulting single scale  phase-field solution at (a) $\bar{\bm u}= 0.15\;mm$, (b) $\bar{\bm u}= 0.324\;mm$, (c) $\bar{\bm u}= 0.333\;mm$ and (d) $\bar{\bm u}= 0.58\;mm$.
	}
	\label{q2}
\end{figure}
The evolution of the local crack phase-field with the corresponding mesh is
depicted in Fig. \ref{q2} for different deformation stages. Specifically, the
second row in Fig. \ref{q2} corresponds to the deformations $\bar{\bm u}=
0.324\;mm$, $\bar{\bm u}= 0.333\;mm$ and $\bar{\bm u}= 0.58\;mm$,
respectively. Due to the existing hard/stiff inclusions, the crack phase-field
propagates around the inclusions. The
resulting crack pattern indicated in Fig. \ref{q2} demonstrates an excellent
agreement with the single scale simulation, with the advantage that
the Global-Local approach requires significantly less 
degrees of freedom.  

A comparison of the load-displacement curves is shown in Fig. \ref{q4}a. 
Therein, a good agreement of the Global-Local approach with the single scale solution was observed for the heterogeneous L-shaped panel test. Figure \ref{q4}b
illustrates the efficiency of the Global-Local approach. 
Here the accumulative computational time is reduced by a factor of eight.
\begin{figure}[!t]
	\centering
	{\includegraphics[clip,trim=1cm 15.5cm 0cm 11cm, width=17cm]{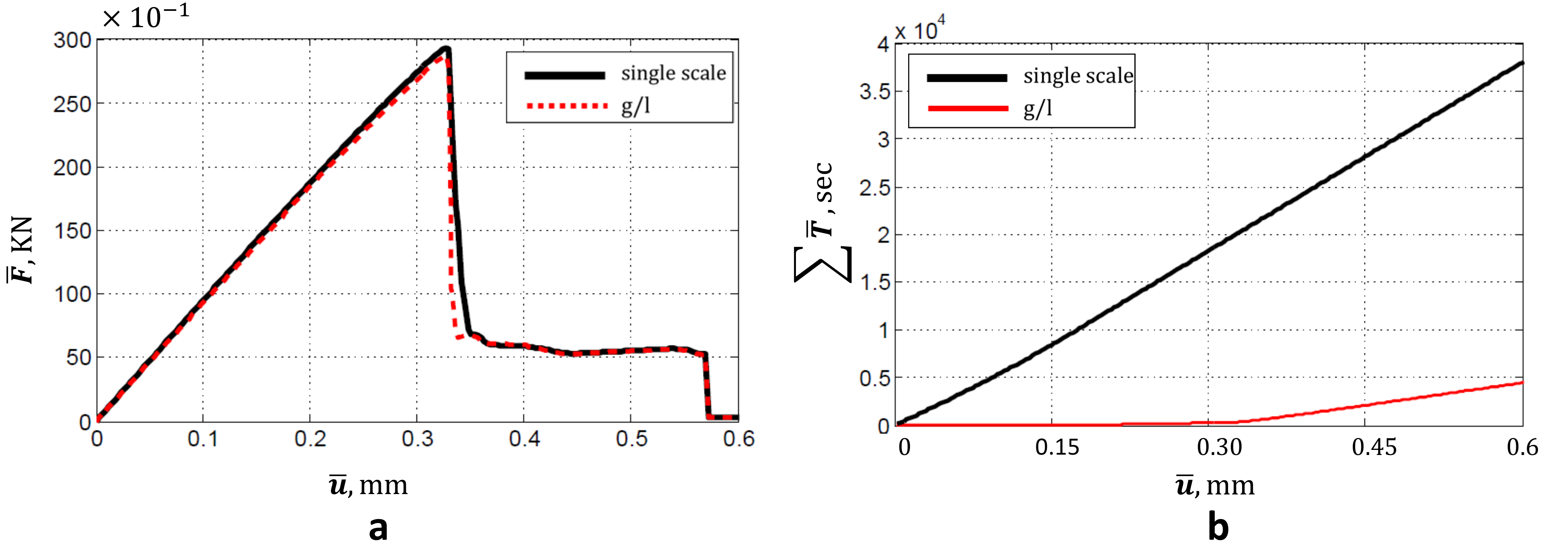}}  
	\caption{Example 3. Heterogeneous L-shaped panel test. (a) Comparison of the load-displacement curves and (b) accumulated time-displacement curves.
	}
	\label{q4}
\end{figure}
\sectpb[Section53]{{Example 4:} Investigation of transversely isotropic double-edge-notched tension}
The last example is concerned with the capability of the proposed Global-Local approach for handling \textit{coalescence} and \textit{merging} of \textit{crack paths} in the local domains. Specifically, the following numerical test aims to illustrate the effects of the double notch shaped specimen. Here crack-initiation and curved-crack-propagation, representing a mixed-mode fracture, are predicted with a Global-Local formulation.  Additionally a transversely isotropic material behavior given by the structural director $\bm{a}$ that is inclined under $\phi=-15^\circ$ is assumed.

The geometrical setup and the loading conditions of the notched specimen is depicted in Fig. \ref{F1}a. The bottom edge of the plate is fixed in the $x$ and $y$ directions. A vertical displacement is applied at the top edge until final failure. We set $A=20\;mm$ and $B=10\;mm$ hence $\calB=(20,10)^2$ $mm^2$. For the double-edge-notches, let $H_1=5.5\;mm$ and $H_2=3.5\;mm$ with the predefined crack length of $l_0=5\;mm$, see Fig. \ref{F1}a. The material parameters used in the simulation are the same as in \cite{aldakheel+blaz+wriggers18} and set as: $\lambda = 12$ kN/mm$^2$, $\mu = 8$ kN/mm$^2$, $G_c = 1 \times 10^{-3}$ kN/mm, $\alpha= \chi = 50$ and $\Xi =0$. 

The global finite element mesh includes two potential fictitious zones $\calB_{F,1}$ and $\calB_{F,2}$ with the interfaces $\Gamma_1$ and $\Gamma_2$ shown in Fig, \ref{F1}b.

\begin{figure}[!t]
	\centering
	{\includegraphics[clip,trim=1cm 19cm 0cm 10cm, width=16cm]{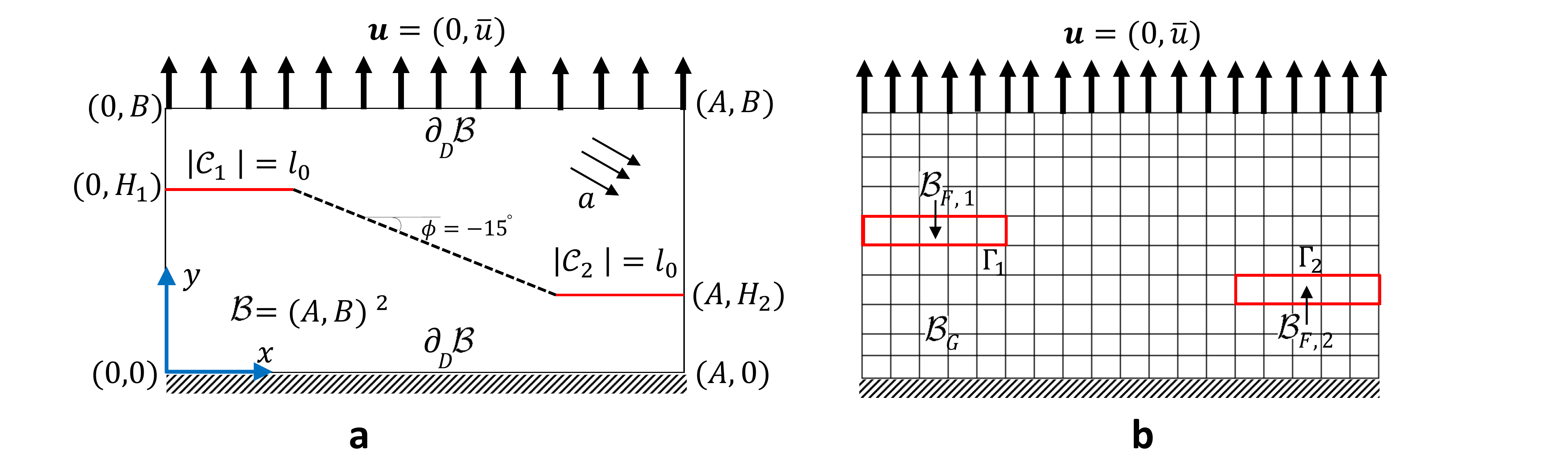}}  
	\caption{Example 4. Double-edge-notched tensile test. (a) Geometry and loading setup with a structural director $\bm{a}$ inclined under an angle $\phi=-15$ and (b) global finite element mesh with two potential fictitious zone $\calB_{F,1}$ and $\calB_{F,2}$.
	}
	\label{F1}
\end{figure}
\begin{figure}[!b]
	\centering
	{\includegraphics[clip,trim=0cm 16cm 0cm 6cm, width=17cm]{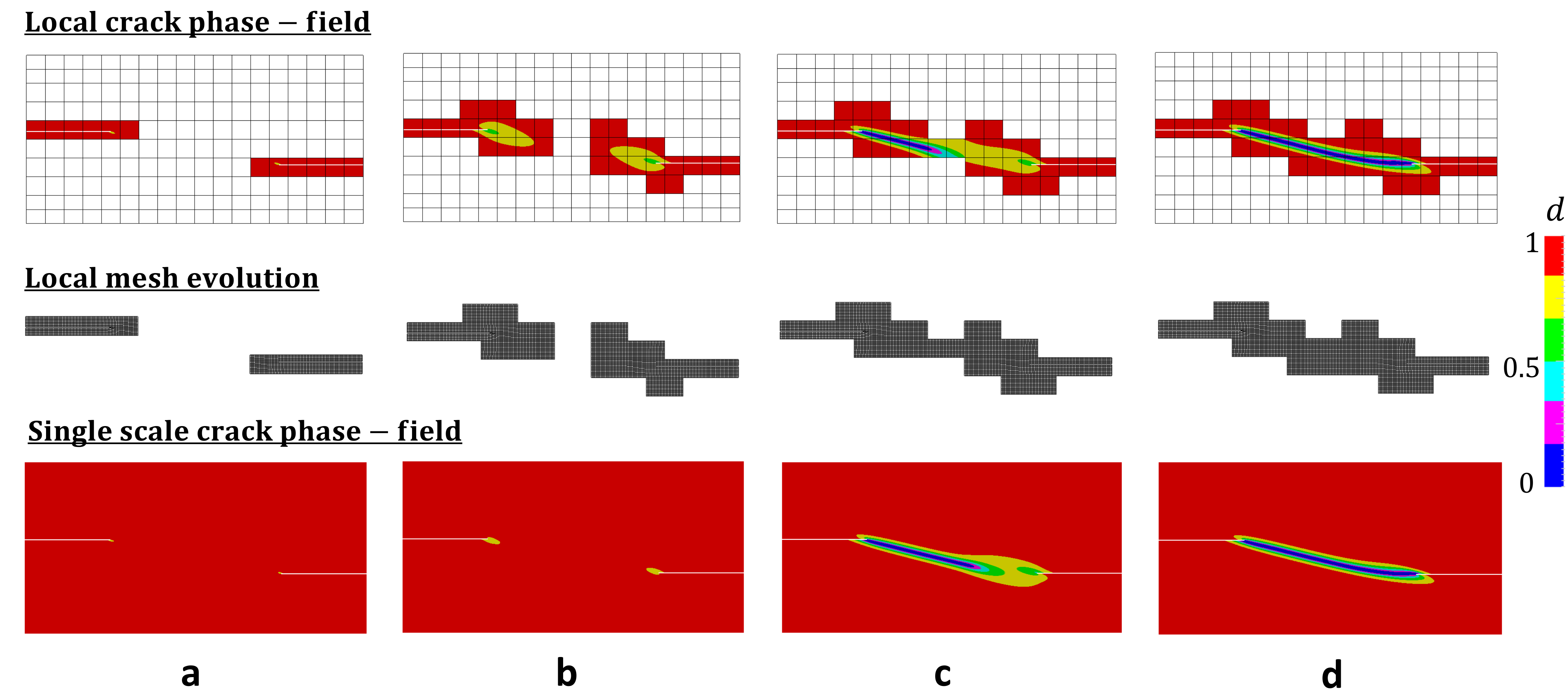}}  
	\caption{Example 4. Crack phase-field pattern for transversely isotropic double-edge-notched plate
		under tension with fiber direction angle of $\phi=-15^\circ$.
		First row: local crack phase-field based on the adaptive scheme; Second row: mesh evolution for local domain; Third row: resulting single scale  phase-field solution at (a) $\bar{\bm u}= 0.01\;mm$, (b) $\bar{\bm u}= 0.0122\;mm$, (c) $\bar{\bm u}= 0.0143\;mm$ and (d) $\bar{\bm u}= 0.0145\;mm$.
	}
	\label{F2}
\end{figure}

The evolution of the crack phase-field resulting from the local domain along with its mesh evolution are indicted in Fig. \ref{F2} for different deformation stages up to final failure. The single scale crack phase-field simulation is shown in Fig. \ref{F2} in the third row. At the loading step $\bar{\bm u}= 0.0143\;mm$ the coalescence and merging of the two local domains are observed. Furthermore, a considerable reduction of the unknowns versus the single scale problem, can be clearly noticed due to the very narrow $\calB_L$ in Fig. \ref{F2}. 

These results demonstrate the feasibility/efficiency of using the proposed adaptive Global-Local approach for different local domains and more complicated structures.

\sectpa[Section7]{Conclusion}
In this work, we developed a robust and efficient Global-Local approach for treating
phase-field fracture problems. Anisotropic heterogeneous
materials, due to the highly oscillating response on the local level, often
require {multi-scale} treatments such that the full resolution on the
local scale must be taken into account. {Multi-scale approaches are} advantageous when large-scale problems are considered in which the fracture state
only develops in smaller, localized, regions. Our first aim 
was to extend the phase-field model to anisotropic constitutive laws.
Next, {we developed} an adaptive scheme in which 
the local domains are advanced dynamically during a computation. This approach 
is realized in terms of a predictor-corrector scheme. First, the new local domains 
are predicted and then the solution is again computed considering the new local regions. An important challenge are interface conditions for the coupling of the two-nested models denoted by the local and global domains in a {variational} consistent way. For that purpose Robin-type boundary conditions were formulated. Moreover, non-matching grids are allowed such that the local and global domains can be updated independently of any additional conditions on the mesh.
Several numerical examples are presented to substantiate 
our algorithmic developments. Here, we considered isotropic and 
anisotropic constitutive materials equations. The focus was on 
crack patterns, load-displacement curves and investigations of the
thermodynamic consistency in terms of the energy functional, the efficiency of the adaptive schemes and the choice of different threshold values for detecting the local
domain. The overall response of the Global-Local approach in terms of
accuracy, robustness and efficiency for the two nested finite element models
was verified using {single scale problems}. {In all examples, 
an excellent performance of the proposed framework was observed.}

\subsection*{Acknowledgment}

NN {was partially supported} by the \textsc{Priority Program DFG - SPP 1748} 
under the project \textsc{WI 4367/2-1}. FA was founded by the \textsc{Priority Program DFG - SPP 2020} under the project \textsc{WR 19/58-1}. TW and PW were funded by the
Deutsche Forschungsgemeinschaft (DFG, German Research Foundation)
under Germany's Excellence Strategy within the \textsc{Cluster of
Excellence PhoenixD (EXC 2122)}, project ID 390833453.

\begin{Appendix}
	\setcounter{equation}{0}
	\renewcommand{\theequation}{A.\arabic{equation}}
\subsection*{Appendix A. Finite Element Discretization}
Let $\mathcal{P}$ be a finite element partition of $\calB$ into quadrilaterals, and $N_i$ denotes the nodal shape function associated with the node $i$. The scalar-valued quantity $\hat{\bullet}_i$ represents the nodal value. For the global-local formulation, we assume the existence of the partitions $\mathcal{P}_G$ and $\mathcal{P}_L$. The solution discretization are given by
\begin{equation}
\bm u_G=\bm N_u^G\hat{\bm u}_G, \quad
\bm u_L=\bm N_u^L\hat{\bm u}_L, \quad
d_L=\bm{N}_d^L\hat{\bm d}_L,
\label{Discr1}
\end{equation}
and its derivative
\begin{equation}
\bm\varepsilon(\bm u_G)=\bm B_u^G\hat{\bm u}_G, \quad
\bm\varepsilon(\bm u_L)=\bm B_u^L\hat{\bm u}_L, \quad
\nabla d_L=\bm B_d^L\hat{\bm d}_L.
\label{Discr2}
\end{equation}
To construct the discretization of the Lagrange multipliers $\bm\lambda_C$, $\bm\lambda_L$, $\bm u_\Gamma$ and the supplementary quantity $\bm\lambda_F$ on $\Gamma$, we write
\begin{equation}
\bm\lambda_C=\bm N_\lambda^G\hat{\bm\lambda}_C, \quad
\bm\lambda_L=\bm N_\lambda^L\hat{\bm\lambda}_L, \quad
\bm u_\Gamma=\bm N_u^\Gamma\hat{\bm u}_\Gamma, \quad
\bm\lambda_F=\bm N_\lambda^G\hat{\bm\lambda}_F.
\label{Discr3}
\end{equation}

{We assume that $\bm N_{\lambda}^G= \bm N_{u}^G=\bm N_{u}^\Gamma$ 
and $\bm N_\lambda^L:=\bm N_{u}^L$. 
{This discretization for the Lagrange multipliers satisfies the {\em inf-sup} condition}, see e.g.\ \cite{Wohlmuth,Wohlmuth2003}.}

Thus, coupling terms are discretized by
\begin{equation}\label{var13}
\begin{aligned}
&{{\bf{J}}_G= \int_{{\mathit{\calB}}_G}({{\bm N_\lambda^G})^{T} {\bm N_u^G} \ \mathrm{d}{\Bx}} ,\ \ \;\;\;\; {\bf{J}}_L= \int_{{\mathit{\calB}}_L}({{\bm N_\lambda^L})^{T} {\bm N_u^L} \ {\,\mathrm{d}{\Bx}}}}, \\ 
&{{\bf{L}}_G= \int_{{\Gamma}_{G}}({\bm N_\lambda^G})^{T}{\bm N_u^{\Gamma}}\ {\mathrm{d}s}} , \ \ \;\;\;\;{{\bf{L}}_L= \int_{{\Gamma}_{L}}({\bm N_\lambda^L})^{T}{\bm N_u^{\Gamma}}\ {\mathrm{d}s}} .
\end{aligned}
\end{equation}

Following our discretization, ${\bf{J}}_G:\calB_G\rightarrow\Gamma_G$ and
${\bf{J}}_L:\calB_L\rightarrow\Gamma_L$ 
become signed Boolean mappings \cite{Belytschko} 
which are used to project the entire domain to interface contributions, such that
\begin{equation}\label{bool_map}
{\hat{\bm u}_{G,b}}:={\bf{J}}_G{\hat{\bm u}_G}\AND{\hat{\bm u}_{L,b}}:={\bf{J}}_L{\hat{\bm u}_L}.
\end{equation}

Here, $b$ are denoted as interface nodes. 
{In order to handle a non-matching finite element discretization on the interface, more specifically to compute ${\bf L}_L$ and ${\bf L}_G$, a dual mortar
method \cite{Wohlmuth} is implemented. This provides sufficient regularity of the underlying FE meshes.} 

	\setcounter{equation}{0}
	\renewcommand{\theequation}{B.\arabic{equation}}
\subsection*{Appendix B. Derivation of Robin-Type Boundary conditions}

In this section, we {investigate} the relationship between ${\Delta {\hat{\bm u}}}$ and $\Delta{\hat{\bm\lambda}}$ (in the incremental sense) for the complementary, fictitious and local domain at the converged solution state. Doing that, Robin-type boundary conditions can be derived such that all coupling terms given in (\ref{Coupl1}), (\ref{Coupl2}) and (\ref{Coupl3}) satisfied, simultaneously, at a Global-Local iteration $k$.

Recall the complementary term used in Eq. \ref{E11pE21pE3} and let $\bm u_C$ and $\bm\lambda_C$ be the stationary of the following functional,		
\begin{equation}
{\mathcal L}={\mathcal L}({\bm u}_C,\bm\lambda_C;\BM):
= \int_{\calB_C} w(\Bve_C, 1, 1;\BM) \,\mathrm{d}{\bm x}
+ \int_\Gamma \bm\lambda_C\cdot({{\bm u}_{\Gamma}}-{\bm u}_C) \,\mathrm{d}s
- \int_{\Gamma_{N,C}} {\bm {\bar\tau}} \cdot \bm u_C\,\mathrm{d}s.
\label{eqdefA1}
\end{equation}	
Here, $\Gamma\in \R^{\delta-1}\subset \calB_C$ is denoted as an interface and $\;{{\bm u}_{\Gamma}}:=\mbox{tr}\;{\bm u}_C\in{\bf H}^{1/2}(\Gamma)$ can be given implicitly, i.e. (\ref{Coupl7nL})$+$(\ref{Coupl6nL}) or explicitly, i.e. (\ref{Coupl6nG}). Recall Eq. \ref{eqdefA1} lives in $\calB_C$ (the following description holds true for $\calB_F$ except ${\bm {\bar \tau}}={\bm 0}$). The stationary points of the energy functional for the ${\mathcal L}$ is characterized 	by the first-order necessary conditions through ${\bf {\mathcal L}_1}={{\mathcal L}_{\bm u}({\bm u}_C,\bm\lambda_C;{\bm w})=\bm 0}\;\mbox{and}\;{{\bf {\mathcal L}_2}={\mathcal L}_{\bm \lambda}({\bm u}_C,\bm\lambda_C;{\bm \kappa})=\bm 0}$. We split ${\bf {\mathcal L}_1}$ into inner nodes and interface nodes denoted as, $\{a,\;b\}$, respectively, by

\begin{equation}
\begin{aligned}
&\displaystyle {{{\bf {\mathcal L}}^a_1}({\bm u})={\bm{ \texttt{f}}^{\:a}}-{ \bar{\bf F}}\overset{!}{=}{\bm 0}} \quad\;
\quad\;\;{\bm x \in \calB\backslash\Gamma},\\ 
&{{{\bf {\mathcal L}}^{b}_1}({\bm u}){=}{\bm{ \texttt{f}}^{\:b}-{{\bf{L}}^{T}_G}{\hat{\bm\lambda}}_C}\overset{!}{=}{\bm 0}} \quad\;\;{\bm x \in \Gamma},\\
&\displaystyle {\bf {\mathcal L}_2}={{\bf{L}}_G}\hat{\bm u}_{\Gamma}-{{\bf{J}}_G}{\hat{\bm u}_C}\overset{!}={\bm 0} \;\;\;\;{\bm x \in \Gamma}.
\label{eqdefA2}
\end{aligned}
\end{equation}	
Here, ${\bm{ \texttt{f}}}=\displaystyle\int_{\calB} {(\bm B_u^G)^{T}}{\bm \sigma}({{\bm u}_C}) \,\mathrm{d}{\bm x}$ is an internal nodal force vector and $\bar{{\bf F}}=\displaystyle\int_{\Gamma_{N}} {(\bm N_u^G)^{T}}{\bm {\bar \tau}} \,\mathrm{d}{s}$ stands for the external force vector. It is trivial that the Lagrange multiplier acts as an external {force} on the interface. A Newton-type solution for the residual based system of equations for $({\bm u}_C,\bm\lambda_C)$ is provided by the linearization

\begin{equation}
\begin{aligned}
&\displaystyle {({\bm{ \texttt{f}}^{\:a}}-{\bar{ \bf F}})+{{\bf K}_{aa}}{\Delta {\hat{\bm u}}_{C,a}}+{{\bf K}_{ab}}{\Delta {\hat{\bm u}}_{C,b}}={\bm 0}},\\ 
&{({\bm{ \texttt{f}}^{\:b}-{{\bf{L}}^{T}_G}{\hat{\bm\lambda}_C}})+{{\bf K}_{b a}}{\Delta {\hat{\bm u}}_{C,a}}+{{\bf K}_{bb}}{\Delta {{\hat{\bm u}}}_{C,b}}-{{\bf{L}}^{T}_G}{\Delta\hat{\bm\lambda}_C}={\bm 0}},\\
&({{\bf{L}}_G}{{\hat{\bm u}}}_{\Gamma}-{{\bf{L}}_G}{\Delta{{\hat{\bm u}}}_b})+{{\bf{L}}_G}{\Delta{\hat{\bm u}}}_{\Gamma}-{{\bf{L}}_G}{\Delta{{\hat{\bm u}}}_{C,b}}=0.
\label{eqdefA3}
\end{aligned}
\end{equation}
where ${\bf K}:={\partial{\bm {\texttt{f}}}}/{\partial {\hat{\bm u}_C}}$ is
the standard tangent stiffness matrix. {Here, we have the} following iterative update
\begin{equation}
\begin{aligned}
{{\hat{\bm u}}}_{C,a}\leftarrow{{\hat{\bm u}}}_{C,a}+{\Delta {{\hat{\bm u}}}_{C,a}},\;\;{{\hat{\bm u}}}_{C,b}\leftarrow{{\hat{\bm u}}}_{C,b}+{\Delta {{\hat{\bm u}}}_{C,b}},\;\;\mbox{and}\;\;{\hat{\bm\lambda}_C}\leftarrow{\hat{\bm\lambda}_C}+{\Delta {\hat{\bm\lambda}}_C}.
\label{eqdefA4}
\end{aligned}
\end{equation}
Let us assume that the equilibrium state is achieved such that ${{\bf {\mathcal L}}^a_1}={\bm 0}, {{\bf {\mathcal L}}^b_1}={\bm 0}$ and ${\bf {\mathcal L}_2}={\bm 0}$. Thus, Eq. \ref{eqdefA3} takes the form

\begin{equation}
\begin{aligned}
{\bm{\mathcal{S}}_C}{\Delta{\hat{\bm u}}}_{C,b}
={\bm{\mathcal{S}}_C}{\Delta{\hat{\bm u}}}_{\Gamma}
=:{{\bf{L}}^{T}_G}{\Delta {\hat{\bm\lambda}}_C}\quad\mbox{with}\quad{\bm{\mathcal{S}}}_C:={\bm{\mathcal{S}}}({\bf K}_C)
={\bf K}_{bb}-{\bf K}_{bb a}{\bf K}^{-1}_{aa}{\bf K}_{ab},
\label{eqdefA5}
\end{aligned}
\end{equation}

where ${\bm{\mathcal{S}}}$ refers to {the} \textit{Steklov-Poincar\'e mapping}  \cite{Steklov}. By means of Eq. \ref{eqdefA5}, displacement ${\hat{\bm u}_C}$ is extracted from the interface $\Gamma$ and through the Poincar\'e-Steklov mapping ${\bm{\mathcal{S}}}$ returns the outward normal stress derivative with respect to the trace of the displacement. That is called Dirichlet-to-Neumann mapping \cite{DtN1,DtN2}.

In a similar way, we have the following identity

\begin{equation}
\begin{aligned}
{\bm{\mathcal{S}}_L}{\Delta{\hat{\bm u}}}_{L,b}={{\bf{T}}^{T}_L}{\Delta {\hat{\bm\lambda}_L}}\quad\mbox{and}\quad{\bm{\mathcal{S}}_F}{\Delta{\hat{\bm u}}}_{F,b}={{\bf{L}}^{T}_G}{\Delta {\hat{\bm\lambda}}_F}.
\label{eqdefA51}
\end{aligned}
\end{equation}

Here, ${{\bf{T}}_L}:={\bf{J}}_L|_{\Gamma_L}$ is the restriction of ${\bf{J}}_L$ from $\calB_L$ to $\Gamma_L$. Furthermore, we define ${\bm{\mathcal{S}}_L}:={\bm{\mathcal{S}}}({\bf K}_L)$ and ${\bm{\mathcal{S}}_F}:={\bm{\mathcal{S}}}({\bf K}_F)$ in Eq. \ref{eqdefA51}.

\begin{theorem}
	\label{theorem1}
	Let the global solutions be at the converged state and let the following identity holds true:
	\begin{equation}
	\bm u_\Gamma^{k,\frac{1}{2}}=\bm u_\Gamma^{k}\;\in\;\Gamma,
	\label{robin_cond}
	\end{equation}
	then the Global-Local formulation is converged. In addition, (\ref{robin_cond}) holds true if and only if
	\begin{equation}
	\Delta {\bm\Lambda}_L={\bm\Lambda}^{k}_L-{\bm\Lambda}^{k-1}_L=0.
	\label{robin_cond2}
	\end{equation}	
	
\end{theorem} 

\textbf{Proof.} The proof constitutes of two parts. Note, the Global-Local procedure is in the convergence state if, all the coupling terms (\ref{Coupl1}), (\ref{Coupl2}) and (\ref{Coupl3}) holds true at iteration $k$.

(\textit{a}) Let condition $\bm u_\Gamma^{k,\frac{1}{2}}=\bm u_\Gamma^{k}$ hold, then it is evident (\ref{Coupl2}) and (\ref{Coupl3}) are satisfied in iteration $k$. 
Accordingly, replacing Eq. \ref{rhs_robin_L} in (\ref{Coupl7nG}) yields

\begin{equation}
\int_\Gamma (\bm\lambda^{k}_C+\bm\lambda^{k}_L) \cdot \delta {\bm u}_\Gamma \,\mathrm{d}s+\KIA_G\int_\Gamma (\bm u_\Gamma^{k}-\bm u^{k}_L) \cdot {\delta {\bm \lambda}_L} \,\mathrm{d}s=0,
\label{proofa1}
\end{equation}

where the second term due to (\ref{Coupl6nL}) and the identity of $\bm u_\Gamma^{k,\frac{1}{2}}=\bm u_\Gamma^{k}$ becomes zero. Hence this results in the continuity of tractions at iteration $k$, i.e. (\ref{Coupl1}) is satisfied.

(\textit{b}) Let $\Delta {\bm\Lambda}_L={\bm\Lambda}^{k}_L-{\bm\Lambda}^{k-1}_L=0$ hold, then (\ref{rhs_robin_G}) can be restated as
\begin{equation}
\begin{aligned}
{\bm\Lambda}^{k-1}_L={\bm\Lambda}^{k}_L
&=\KIA_L\int_\Gamma \bm u_G^{k} \cdot {\delta {\bm \lambda}_C} \,\mathrm{d}s-\int_\Gamma \bm\lambda^{k}_C\cdot \bm v_\Gamma \,\mathrm{d}s,
\end{aligned}
\label{rhs_L_k}
\end{equation}
and therefore (\ref{Coupl7nL}) is restated as,
\begin{equation}
\int_\Gamma (\bm\lambda^{k}_C+\bm\lambda^{k}_L) \cdot {\delta {\bm u}_\Gamma} \,\mathrm{d}s+\KIA_L\int_\Gamma (\bm u_\Gamma^{k,\frac{1}{2}}-\bm u^{k}_G) \cdot {\delta {\bm \lambda}_C} \,\mathrm{d}s=\int_\Gamma (\bm\lambda^{k}_C+\bm\lambda^{k}_L) \cdot {\delta {\bm u}_\Gamma} \,\mathrm{d}s=0,
\label{Coupl5nL_recal}
\end{equation}
where (\ref{Coupl6nG}) is used. Subtracting Eq. \ref{Coupl5nL_recal} from (\ref{Coupl7nG}) yields
\begin{equation}
\KIA_G\int_\Gamma (\bm u_\Gamma^{k}-\bm u^{k}_L) \cdot {\delta {\bm \lambda}_L} \,\mathrm{d}s=\KIA_G\int_\Gamma (\bm u_\Gamma^{k}-\bm u_\Gamma^{k,\frac{1}{2}}) \cdot {\delta {\bm \lambda}_L} \,\mathrm{d}s=0,
\label{Coupl5nG_recal}
\end{equation}
which results in $\bm u_\Gamma^{k}=\bm u_\Gamma^{k,\frac{1}{2}}$ . Herein, Eq. \ref{Coupl6nL} is used. Reciprocally, if $\bm u_\Gamma^{k,\frac{1}{2}}=\bm u_\Gamma^{k}\;\in\;\Gamma$ satisfied then $\Delta {\bm\Lambda}_L=0$ holds. The proof is left for the readers. $\qedsymbol$

\textbf{Remark B.1.}
\label{indentical_gl_ref}
{\it Let the Global-Local approach {be in the} converged state such that $\bm u_\Gamma^{k,\frac{1}{2}}=\bm u_\Gamma^{k}\;\in\;\Gamma$. 
We underline if $\calB=\calB_C\cup\Gamma\cup\calB_L$ holds, 
then {the} Global-Local mesh compared with {a} single scale mesh is one to one, 
{and the} resulting Global-Local solutions are identical with the single scale solutions.
	{This} is because {the} two finite element meshes which correspond to 
the complementary and local domains are `exactly' identical 
{to} {the} single scale mesh; {we refer the reader to} Section \ref{Section521} for case a.
}

We now {determine} specific Robin-type boundary conditions such that $\Delta {\bm\Lambda}_L=0$ holds which results in $\bm u_\Gamma^{k}=\bm u_\Gamma^{k,\frac{1}{2}}$ by means of Proposition \ref{theorem1} and that yields the Global-Local iterative process to be in the converged state. 

Recall (\ref{rhs_robin_G}) and find $\KIA_L$ such that $\Delta {\bm\Lambda}_L={\bm\Lambda}^{k}_L-{\bm\Lambda}^{k-1}_L=0$. Hence we have,
\begin{equation}\label{proofa11}
\Delta {\bm\Lambda}_L=\KIA_L\int_\Gamma \Delta \bm u_G \cdot {\delta {\bm \lambda}_C} \,\mathrm{d}s-\int_\Gamma \Delta \bm\lambda_C\cdot {\delta {\bm u}_\Gamma} \,\mathrm{d}s=0,
\end{equation}
resulting to
\begin{equation}\label{proofa12}
\KIA_L\int_\Gamma \Delta \bm u_G \cdot {\delta {\bm \lambda}_C} \,\mathrm{d}s=\int_\Gamma \Delta \bm\lambda_C\cdot {\delta {\bm u}_\Gamma} \,\mathrm{d}s\;\;\rightarrow\;\;\KIA_L{\bf{J}}_G{\Delta\hat{\bm u}_G}={{\bf{L}}^{T}_G}{\Delta \hat{\bm\lambda}_C},
\end{equation}
By means of (\ref{eqdefA5}) and considering (\ref{bool_map})$_1$, this equality holds if $\KIA_L:={\bm {\mathcal{S}}}_C$ which is the Dirichlet-to-Neumann operator assigned to $\Gamma_G \in \calB_C$. 

In a similar manner, let us find $\KIA_G$ such that $\Delta {\bm\Lambda}_G={\bm\Lambda}^{k}_G-{\bm\Lambda}^{k-1}_G=0$. This yields
\begin{equation}
\KIA_G\int_\Gamma \Delta \bm u_L \cdot {\delta {\bm \lambda}_L} \,\mathrm{d}s=\int_\Gamma \Delta \bm\lambda_L\cdot {\delta {\bm u}_\Gamma} \,\mathrm{d}s\;\;\rightarrow\;\;\KIA_G{\bf{J}}_L{\Delta\hat{\bm u}_L}={{\bf{L}}^{T}_L}{\Delta \hat{\bm\lambda}_L}.
\label{proofb9}
\end{equation}
By means of (\ref{eqdefA51})$_1$ and (\ref{bool_map})$_2$, we have $\KIA_G={\bm L^T_L}{\bm T^{-T}_L}{\bm {\mathcal{S}}}_L$. Based on Proposition \ref{theorem1} the converged state of the Global-Local iteration is independent of the choice of $\KIA_G$ hence one can simply replace  ${\bm {\mathcal{S}}}_L$ by the identity tensor. Hence, there is no need to access ${\bm {\mathcal{S}}}_L$ at the global level.

\end{Appendix}
\bibliographystyle{bibls1}
\bibliography{./lit.bib}

\begin{thebibliography}{100}
\setlength{\itemsep}{0mm}
\providecommand{\natexlab}[1]{#1}
\providecommand{\url}[1]{\texttt{#1}}
\expandafter\ifx\csname urlstyle\endcsname\relax
  \providecommand{\doi}[1]{doi: #1}\else
  \providecommand{\doi}{doi: \begingroup \urlstyle{rm}\Url}\fi

\bibitem[\textsc{Aldakheel et~al.}(2014)Aldakheel, Mauthe, \&
  Miehe]{aldakheel+mauthe+miehe14}
\textsc{Aldakheel, F.; Mauthe, S.; Miehe, C.} [2014]: \emph{Towards phase field
  modeling of ductile fracture in gradient-extended elastic-plastic solids}.
\newblock Proceedings in Applied Mathematics and Mechanics, 14: 411--412.

\bibitem[\textsc{Aldakheel et~al.}(2018{\natexlab{a}})Aldakheel, Hudobivnik,
  Hussein, \& Wriggers]{aldakheel+blaz+wriggers18}
\textsc{Aldakheel, F.; Hudobivnik, B.; Hussein, A.; Wriggers, P.} [2018]:
  \emph{Phase-field modeling of brittle fracture using an efficient virtual
  element scheme}.
\newblock Computer Methods in Applied Mechanics and Engineering, 341: 443--466.

\bibitem[\textsc{Aldakheel et~al.}(2018{\natexlab{b}})Aldakheel, Hudobivnik, \&
  Wriggers]{aldakheeletal18}
\textsc{Aldakheel, F.; Hudobivnik, B.; Wriggers, P.} [2018]: \emph{Virtual
  element formulation for phase-field modeling of ductile fracture}.
\newblock International Journal for Multiscale Computational Engineering.
\newblock DOI: 10.1615/IntJMultCompEng.2018026804.

\bibitem[\textsc{Aldakheel}(2016)]{aldakheel16}
\textsc{Aldakheel, F.} [2016]: \emph{Mechanics of Nonlocal Dissipative Solids:
  Gradient Plasticity and Phase Field Modeling of Ductile Fracture}.
\newblock Ph.D. Thesis, Institute of Applied Mechanics (CE), Chair I,
  University of Stuttgart.
\newblock http://dx.doi.org/10.18419/opus-8803.

\bibitem[\textsc{Aldakheel et~al.}(2018{\natexlab{c}})Aldakheel, Wriggers, \&
  Miehe]{aldakheel+wriggers+miehe18}
\textsc{Aldakheel, F.; Wriggers, P.; Miehe, C.} [2018]: \emph{A modified
  gurson-type plasticity model at finite strains: Formulation, numerical
  analysis and phase-field coupling}.
\newblock Computational Mechanics, 62: 815--833.

\bibitem[\textsc{Ambati et~al.}(2015)Ambati, Gerasimov, \&
  De~Lorenzis]{ambati15}
\textsc{Ambati, M.; Gerasimov, T.; De~Lorenzis, L.} [2015]: \emph{A review on
  phase-field models of brittle fracture and a new fast hybrid formulation}.
\newblock Computational Mechanics, 55(2): 383--405.

\bibitem[\textsc{Ambrosio \& Tortorelli}(1990)]{AmTo90}
\textsc{Ambrosio, L.; Tortorelli, V.} [1990]: \emph{Approximation of
  functionals depending on jumps by elliptic functionals via
  $\gamma$-convergence}.
\newblock Communications on Pure and Applied Mathematics, 43: 999--1036.

\bibitem[\textsc{Ambrosio \& Tortorelli}(1992)]{AmTo92}
\textsc{Ambrosio, L.; Tortorelli, V.} [1992]: \emph{On the approximation of
  free discontinuity problems}.
\newblock Bollettino dell'Unione Matematica Italiana, 6: 105--123.

\bibitem[\textsc{Amor et~al.}(2009)Amor, Marigo, \&
  Maurini]{amor+marigo+maurini09}
\textsc{Amor, H.; Marigo, J.; Maurini, C.} [2009]: \emph{Regularized
  formulation of the variational brittle fracture with unilateral contact:
  Numerical experiments}.
\newblock Journal of the Mechanics and Physics of Solids, 57: 1209--1229.

\bibitem[\textsc{Balzani et~al.}(2006)Balzani, Neff, Schr\"oder, \&
  Holzapfel]{balzani+etal06}
\textsc{Balzani, D.; Neff, P.; Schr\"oder, J.; Holzapfel, G.} [2006]: \emph{A
  polyconvex framework for soft biological tissues. adjustment to experimental
  data}.
\newblock International Journal of Solids and Structures, 43: 6052--6070.

\bibitem[\textsc{Bangerth et~al.}(2007)Bangerth, Hartmann, \&
  Kanschat]{BangerthHartmannKanschat2007}
\textsc{Bangerth, W.; Hartmann, R.; Kanschat, G.} [2007]: \emph{{deal.II} -- a
  general purpose object oriented finite element library}.
\newblock ACM Transactions on Mathematical Software, 33(4): 24/1--24/27.

\bibitem[\textsc{Belytschko et~al.}(2014)Belytschko, Liu, Moran, \&
  Elkhodary]{Belytschko}
\textsc{Belytschko, T.; Liu, W.~K.; Moran, B.; Elkhodary, K.} [2014]:
  \emph{Nonlinear Finite Elements for Continua and Structure}.
\newblock John Wiley and Sons, Ltd., United Kingdom.

\bibitem[\textsc{Bleyer \& Alessi}(2018)]{bleyer+alessi18}
\textsc{Bleyer, J.; Alessi, R.} [2018]: \emph{Phase-field modeling of
  anisotropic brittle fracture including several damage mechanisms}.
\newblock Computer Methods in Applied Mechanics and Engineering.
\newblock https://doi.org/10.1016/j.cma.2018.03.012.

\bibitem[\textsc{Borden et~al.}(2012)Borden, Verhoosel, Scott, Hughes, \&
  Landis]{borden+verhoosel+scott+hughes+landis12}
\textsc{Borden, M.~J.; Verhoosel, C.~V.; Scott, M.~A.; Hughes, T. J.~R.;
  Landis, C.~M.} [2012]: \emph{A phase-field description of dynamic brittle
  fracture}.
\newblock Computer Methods in Applied Mechanics and Engineering, 217-220:
  77--95.

\bibitem[\textsc{Bourdin}(1999)]{Bou99}
\textsc{Bourdin, B.} [1999]: \emph{Image segmentation with a finite element
  method}.
\newblock Mathematical Modelling and Numerical Analysis, 33(2): 229--244.

\bibitem[\textsc{Bourdin}(2007)]{Bour07}
\textsc{Bourdin, B.} [2007]: \emph{Numerical implementation of the variational
  formulation for quasi-static brittle fracture}.
\newblock Interfaces and free boundaries, 9: 411--430.

\bibitem[\textsc{Bourdin et~al.}(2000)Bourdin, Francfort, \&
  Marigo]{BourFraMar00}
\textsc{Bourdin, B.; Francfort, G.; Marigo, J.-J.} [2000]: \emph{Numerical
  experiments in revisited brittle fracture}.
\newblock Journal of the Mechanics and Physics of Solids, 48(4): 797--826.

\bibitem[\textsc{Bourdin et~al.}(2008)Bourdin, Francfort, \&
  Marigo]{BourFraMar08}
\textsc{Bourdin, B.; Francfort, G.; Marigo, J.-J.} [2008]: \emph{The
  variational approach to fracture}.
\newblock Journal of Elasticity, 91: 5--148.

\bibitem[\textsc{Chapman}(2003)]{Chapman03}
\textsc{Chapman, S.~J.} [2003]: \emph{Fortran 90/95 for Scientists and
  Engineers}.
\newblock McGraw-Hill, Inc., New York, NY, USA, 2 Edition.

\bibitem[\textsc{Chevreuil et~al.}(2013)Chevreuil, Nouy, \& Safatly]{Chevreuil}
\textsc{Chevreuil, M.; Nouy, A.; Safatly, E.} [2013]: \emph{A multiscale method
  with patch for the solution of stochastic partial differential equations with
  localized uncertainties}.
\newblock Computer Methods in Applied Mechanics and Engineering, 255: 255--274.

\bibitem[\textsc{Ciarlet}(1987)]{Cia87}
\textsc{Ciarlet, P.~G.} [1987]: \emph{The finite element method for elliptic
  problems}.
\newblock North-Holland, Amsterdam [u.a.], 2. pr. Edition.

\bibitem[\textsc{Deparis et~al.}(2007)Deparis, Discacciati, Fourestey, \&
  Quarteroni]{DtN1}
\textsc{Deparis, S.; Discacciati, M.; Fourestey, G.; Quarteroni, A.} [2007]:
  \emph{Heterogeneous domain decomposition methods for fluid-structure
  interaction problems}.
\newblock Domain Decomposition Methods in Science and Engineering XVI. Lecture
  Notes in Computational Science and Engineering, 55: 41--52.

\bibitem[\textsc{Dittmann et~al.}(2018)Dittmann, Aldakheel, Schulte, Wriggers,
  \& Hesch]{dittmann+fadi+etal18}
\textsc{Dittmann, M.; Aldakheel, F.; Schulte, J.; Wriggers, P.; Hesch, C.}
  [2018]: \emph{Variational phase-field formulation of non-linear ductile
  fracture}.
\newblock Computer Methods in Applied Mechanics and Engineering, 342: 71--94.

\bibitem[\textsc{Donath}(1961)]{DONATH61}
\textsc{Donath, F.~A.} [1961]: \emph{{Experimental study of shear failure in
  anisotropic rocks}}.
\newblock GSA Bulletin, 72(6): 985--989.

\bibitem[\textsc{Farhat \& Roux}(1991)]{Farhat1991}
\textsc{Farhat, C.; Roux, F.} [1991]: \emph{A method of finite element tearing
  and interconnecting and its parallel solution algorithm}.
\newblock International Journal for Numerical Methods in Engineering, 32:
  1205--1227.

\bibitem[\textsc{Farhat et~al.}(2000)Farhat, Macedo, Lesoinne, Roux,
  Magoul\'es, \& Bourdonnaie]{Farhat2LM20002}
\textsc{Farhat, C.; Macedo, A.; Lesoinne, M.; Roux, F.-X.; Magoul\'es, F.;
  Bourdonnaie, A. D.~L.} [2000]: \emph{Two-level domain decomposition methods
  with lagrange multipliers for the fast iterative solution of acoustic
  scattering problems}.
\newblock Computer Methods in Applied Mechanics and Engineering, 184: 213--239.

\bibitem[\textsc{Farrell \& Maurini}(2017)]{FaMau17}
\textsc{Farrell, P.~E.; Maurini, C.} [2017]: \emph{Linear and nonlinear solvers
  for variational phase-field models of brittle fracture}.
\newblock International Journal for Numerical Methods in Engineering, 109:
  648--667.

\bibitem[\textsc{Fish \& Wagiman}(1993)]{Fish3}
\textsc{Fish, J.; Wagiman, A.} [1993]: \emph{Multiscale finite element method
  for a locally nonperiodic heterogeneous medium}.
\newblock Computational Mechanics, 12(3): 164--180.

\bibitem[\textsc{Fish}(2014)]{Fish2014}
\textsc{Fish, J.} [2014]: \emph{Practical Multiscaling}.
\newblock John Wiley and Sons, Ltd., United Kingdom.

\bibitem[\textsc{Francfort \& Marigo}(1998)]{FraMar98}
\textsc{Francfort, G.; Marigo, J.-J.} [1998]: \emph{Revisiting brittle fracture
  as an energy minimization problem}.
\newblock Journal of the Mechanics and Physics of Solids, 46(8): 1319--1342.

\bibitem[\textsc{Gander et~al.}(2007)Gander, Halpern, \& Magoules]{Gander2007}
\textsc{Gander, M.; Halpern, L.; Magoules, F.} [2007]: \emph{An optimized
  schwarz method with two-sided robin transmission conditions for the helmholtz
  equation}.
\newblock International Journal for Numerical Methods in Fluids, 55: 163--175.

\bibitem[\textsc{Gendre et~al.}(2009)Gendre, Allix, Gosselet, \&
  Comte]{Allix09}
\textsc{Gendre, L.; Allix, O.; Gosselet, P.; Comte, F.} [2009]:
  \emph{Non-intrusive and exact global/local techniques for structural problems
  with local plasticity}.
\newblock Computational Mechanics, 44: 233--245.

\bibitem[\textsc{Gendre et~al.}(2011)Gendre, Allix, \& Gosselet]{Gendre2}
\textsc{Gendre, L.; Allix, O.; Gosselet, P.} [2011]: \emph{A two-scale
  approximation of the schur complement and its use for non-intrusive
  coupling}.
\newblock International Journal for Numerical Methods in Engineering, 87:
  889--905.

\bibitem[\textsc{Gerasimov et~al.}(2018)Gerasimov, Noii, Allix, \&
  De~Lorenzis]{NoiiGL18}
\textsc{Gerasimov, T.; Noii, N.; Allix, O.; De~Lorenzis, L.} [2018]: \emph{A
  non-intrusive global/local approach applied to phase-field modeling of
  brittle fracture}.
\newblock Advanced Modeling and Simulation in Engineering Sciences.
\newblock https://doi.org/10.1186/s40323-018-0105-8.

\bibitem[\textsc{Gosselet \& Rey}(2006)]{DDRey06}
\textsc{Gosselet, P.; Rey, C.} [2006]: \emph{Non-overlapping domain
  decomposition methods in structural mechanics}.
\newblock Archives of Computational Methods in Engineering, 13: 515--572.

\bibitem[\textsc{Greer \& Loisel}(2015)]{DtN2}
\textsc{Greer, N.; Loisel, S.} [2015]: \emph{The optimised schwarz method and
  the two-lagrange multiplier method for heterogeneous problems in general
  domains with two general subdomains}.
\newblock Numerical Algorithms, 69: 737--762.

\bibitem[\textsc{G\"ultekin et~al.}(2016)G\"ultekin, Dal, \&
  Holzapfel]{GULTEKIN16}
\textsc{G\"ultekin, O.; Dal, H.; Holzapfel, G.~A.} [2016]: \emph{A phase-field
  approach to model fracture of arterial walls: Theory and finite element
  analysis}.
\newblock Computer Methods in Applied Mechanics and Engineering, 312: 542 --
  566.
\newblock Phase Field Approaches to Fracture.

\bibitem[\textsc{Hakim \& Karma}(2009)]{HAKIM2009}
\textsc{Hakim, V.; Karma, A.} [2009]: \emph{Laws of crack motion and
  phase-field models of fracture}.
\newblock Journal of the Mechanics and Physics of Solids, 57(2): 342 -- 368.

\bibitem[\textsc{Hautefeuille et~al.}(2012)Hautefeuille, Colliat,
  Ibrahimbegovic, Matthies, \& Villon]{Hautefeuille}
\textsc{Hautefeuille, M.; Colliat, J.-B.; Ibrahimbegovic, A.; Matthies, H.;
  Villon, P.} [2012]: \emph{A multi-scale approach to model localized failure
  with softening}.
\newblock Computers \& Structures, 94-95: 83--95.

\bibitem[\textsc{Hecht et~al.}(2009)Hecht, Lozinski, \& Pironneau]{Lozinski09}
\textsc{Hecht, F.; Lozinski, A.; Pironneau, O.} [2009]: \emph{Numerical zoom
  and the schwarz algorithm}.
\newblock Domain Decomposition Methods in Science and Engineering XVIII, 70:
  63--73.

\bibitem[\textsc{Heider \& Markert}(2017)]{heider17}
\textsc{Heider, Y.; Markert, B.} [2017]: \emph{A phase-field modeling approach
  of hydraulic fracture in saturated porous media}.
\newblock Mechanics Research Communications, 80: 38--46.

\bibitem[\textsc{Heister et~al.}(2015)Heister, Wheeler, \& Wick]{Wick15Adapt}
\textsc{Heister, T.; Wheeler, M.~F.; Wick, T.} [2015]: \emph{A primal-dual
  active set method and predictor-corrector mesh adaptivity for computing
  fracture propagation using a phase-field approach}.
\newblock Computer Methods in Applied Mechanics and Engineering, 290: 466 --
  495.

\bibitem[\textsc{Heister \& Wick}(2018)]{HeiWi18_pamm}
\textsc{Heister, T.; Wick, T.} [2018]: \emph{Parallel solution, adaptivity,
  computational convergence, and open-source code of 2d and 3d pressurized
  phase-field fracture problems}.
\newblock PAMM, 18(1): e201800353.

\bibitem[\textsc{Hesch \& Weinberg}(2014)]{hesch+weinberg14}
\textsc{Hesch, C.; Weinberg, K.} [2014]: \emph{Thermodynamically consistent
  algorithms for a finite-deformation phase-field approach to fracture}.
\newblock International Journal for Numerical Methods in Engineering, 99:
  906--924.

\bibitem[\textsc{Hill}(1965)]{Hill}
\textsc{Hill, R.} [1965]: \emph{{A self-consistent mechanics of composite
  materials}}.
\newblock Journal of the Mechanics and Physics of Solids, 13(4): 213--222.

\bibitem[\textsc{Hinojosa et~al.}(2014)Hinojosa, Allix, Guidault, \&
  Cresta]{allix1}
\textsc{Hinojosa, J.; Allix, O.; Guidault, P.-A.; Cresta, P.} [2014]:
  \emph{Domain decomposition methods with nonlinear localization for the
  buckling and post-buckling analyses of large structures}.
\newblock Advances in Engineering Software, 70: 13--24.

\bibitem[\textsc{Holzapfel et~al.}(2000)Holzapfel, Gasser, \&
  Ogden]{holzapfel+gasser+ogden00}
\textsc{Holzapfel, G.~A.; Gasser, T., C.; Ogden, R., W.} [2000]: \emph{A new
  constitutive framework for arterial wall mechanics and a comperative study of
  material models}.
\newblock Journal of Elasticity, 61: 1--48.

\bibitem[\textsc{Hughes et~al.}(1998)Hughes, Feij\'oo, Mazzei, \&
  Quincy]{Hughes98}
\textsc{Hughes, T.~J.; Feij\'oo, G.~R.; Mazzei, L.; Quincy, J.-B.} [1998]:
  \emph{The variational multiscale method-a paradigm for computational
  mechanics}.
\newblock Computer Methods in Applied Mechanics and Engineering, 166(1): 3--24.
\newblock Advances in Stabilized Methods in Computational Mechanics.

\bibitem[\textsc{Hughes et~al.}(2005)Hughes, Cottrell, \& Bazilevs]{hughes05}
\textsc{Hughes, T.; Cottrell, J.; Bazilevs, Y.} [2005]: \emph{Isogeometric
  analysis: Cad, finite elements, nurbs, exact geometry and mesh refinement}.
\newblock Computer Methods in Applied Mechanics and Engineering, 194(39): 4135
  -- 4195.

\bibitem[\textsc{Khoromskij \& Wittum}(1998)]{Steklov}
\textsc{Khoromskij, B.; Wittum, G.} [1998]: \emph{Robust Interface Reduction
  for Highly Anisotropic Elliptic Equations}.
\newblock Springer-Verlag Berlin Heidelberg.

\bibitem[\textsc{Kikuchi \& Oden}(1988)]{KiOd88}
\textsc{Kikuchi, N.; Oden, J.} [1988]: \emph{Contact problems in elasticity}.
\newblock Studies in Applied Mathematics. Society for Industrial and Applied
  Mathematics (SIAM), Philadelphia, PA.

\bibitem[\textsc{Kinderlehrer \& Stampacchia}(2000)]{KiStam00}
\textsc{Kinderlehrer, D.; Stampacchia, G.} [2000]: \emph{An Introduction to
  Variational Inequalities and Their Applications}.
\newblock Classics in Applied Mathematics. Society for Industrial and Applied
  Mathematics.

\bibitem[\textsc{Kuhn \& M\"uller}(2010)]{KUHN10}
\textsc{Kuhn, C.; M\"uller, R.} [2010]: \emph{A continuum phase field model for
  fracture}.
\newblock Engineering Fracture Mechanics, 77(18): 3625 -- 3634.

\bibitem[\textsc{Lee et~al.}(2016)Lee, Wheeler, \& Wick]{LeeWheWi16}
\textsc{Lee, S.; Wheeler, M.~F.; Wick, T.} [2016]: \emph{Pressure and
  fluid-driven fracture propagation in porous media using an adaptive finite
  element phase field model}.
\newblock Computer Methods in Applied Mechanics and Engineering, 305: 111 --
  132.

\bibitem[\textsc{Li \& Maurini}(2019)]{LI19}
\textsc{Li, B.; Maurini, C.} [2019]: \emph{Crack kinking in a variational
  phase-field model of brittle fracture with strongly anisotropic surface
  energy}.
\newblock Journal of the Mechanics and Physics of Solids, 125: 502 -- 522.

\bibitem[\textsc{Lloberas-Valls et~al.}(2012)Lloberas-Valls, Rixen, Simone, \&
  Sluys]{Lloberas}
\textsc{Lloberas-Valls, O.; Rixen, D.~J.; Simone, A.; Sluys, L.~J.} [2012]:
  \emph{Multiscale domain decomposition analysis of quasi-brittle heterogeneous
  materials}.
\newblock International Journal for Numerical Methods in Engineering, 83:
  1337--1366.

\bibitem[\textsc{Loehnert \& Belytschko}(2007)]{Loehnert07}
\textsc{Loehnert, S.; Belytschko, T.} [2007]: \emph{A multiscale projection
  method for macro/microcrack simulations}.
\newblock International Journal for Numerical Methods in Engineering, 71(12):
  1466--1482.

\bibitem[\textsc{Maday \& Magoul\'es}(2006)]{Magoules2007}
\textsc{Maday, Y.; Magoul\'es, F.} [2006]: \emph{Absorbing interface conditions
  for domain decomposition methods: A general presentation}.
\newblock International Journal for Numerical Methods in Fluids, 195:
  3880--3900.

\bibitem[\textsc{Magoules et~al.}(2006)Magoules, Roux, \& Series]{Magoules2006}
\textsc{Magoules, F.; Roux, F.~X.; Series, L.} [2006]: \emph{Algebraic
  approximation of dirichlet-to-neumann maps for the equations of linear
  elasticity}.
\newblock Computer Methods in Applied Mechanics and Engineering, 195:
  3742--3759.

\bibitem[\textsc{Mandel}(1993)]{Mandel1993}
\textsc{Mandel, J.} [1993]: \emph{Balancing domain decomposition}.
\newblock Communications in Applied Numerical Methods, 9(4): 233--241.

\bibitem[\textsc{Mang \& Wick}(2019)]{MaWi19}
\textsc{Mang, K.; Wick, T.} [2019]: \emph{Numerical methods for variational
  phase-field fracture problems}.
\newblock Lecture notes at Leibniz University Hannover.

\bibitem[\textsc{Markovic \& Ibrahimbegovic}(2004)]{Markovic}
\textsc{Markovic, D.; Ibrahimbegovic, A.} [2004]: \emph{On micro-macro
  interface conditions for micro scale based fem for inelastic behavior of
  heterogeneous materials}.
\newblock Computer Methods in Applied Mechanics and Engineering, 193(48):
  5503--5523.
\newblock Advances in Computational Plasticity.

\bibitem[\textsc{MATLAB}(2018)]{MATLAB18b}
\textsc{MATLAB} [2018]: \emph{version 9.5.0.944444 (R2018b)}.
\newblock The MathWorks Inc., Natick, Massachusetts.

\bibitem[\textsc{Mesgarnejad et~al.}(2015)Mesgarnejad, Bourdin, \&
  Khonsari]{MESGARNEJAD2015420}
\textsc{Mesgarnejad, A.; Bourdin, B.; Khonsari, M.} [2015]: \emph{Validation
  simulations for the variational approach to fracture}.
\newblock Computer Methods in Applied Mechanics and Engineering, 290: 420--437.

\bibitem[\textsc{Michel et~al.}(1999)Michel, Moulinec, \& Suquet]{Michel}
\textsc{Michel, J.; Moulinec, H.; Suquet, P.} [1999]: \emph{Effective
  properties of composite materials with periodic microstructure: a
  computational approach}.
\newblock Computer Methods in Applied Mechanics and Engineering, 172(1):
  109--143.

\bibitem[\textsc{Miehe \& Bayreuther}(2007)]{Miehe5}
\textsc{Miehe, C.; Bayreuther, C.} [2007]: \emph{On multiscale fe analyses of
  heterogeneous structures: from homogenization to multigrid solvers}.
\newblock International Journal for Numerical Methods in Engineering, 71:
  1135--1180.

\bibitem[\textsc{Miehe et~al.}(2010{\natexlab{a}})Miehe, Hofacker, \&
  Welschinger]{MieWelHof10b}
\textsc{Miehe, C.; Hofacker, M.; Welschinger, F.} [2010]: \emph{A phase field
  model for rate-independent crack propagation: {R}obust algorithmic
  implementation based on operator splits}.
\newblock Comput. Meth. Appl. Mech. Engrg., 199: 2765--2778.

\bibitem[\textsc{Miehe et~al.}(2010{\natexlab{b}})Miehe, Welschinger, \&
  Hofacker]{miehe+welschinger+hofacker10a}
\textsc{Miehe, C.; Welschinger, F.; Hofacker, M.} [2010]:
  \emph{Thermodynamically consistent phase-field models of fracture:
  Variational principles and multi-field fe implementations}.
\newblock International Journal for Numerical Methods in Engineering, 83:
  1273--1311.

\bibitem[\textsc{Miehe et~al.}(2015{\natexlab{a}})Miehe, Hofacker, Sch\"anzel,
  \& Aldakheel]{miehe+hofacker+schaenzel+aldakheel15}
\textsc{Miehe, C.; Hofacker, M.; Sch\"anzel, L.-M.; Aldakheel, F.} [2015]:
  \emph{Phase field modeling of fracture in multi-physics problems. {P}art
  {II}. brittle-to-ductile failure mode transition and crack propagation in
  thermo-elastic-plastic solids}.
\newblock Computer Methods in Applied Mechanics and Engineering, 294: 486--522.

\bibitem[\textsc{Miehe et~al.}(2015{\natexlab{b}})Miehe, Sch\"anzel, \&
  Ulmer]{miehe+schaenzel+ulmer15}
\textsc{Miehe, C.; Sch\"anzel, L.; Ulmer, H.} [2015]: \emph{Phase field
  modeling of fracture in multi-physics problems. {P}art {I}. balance of crack
  surface and failure criteria for brittle crack propagation in thermo-elastic
  solids}.
\newblock Computer Methods in Applied Mechanics and Engineering, 294: 449--485.

\bibitem[\textsc{Miehe et~al.}(2015{\natexlab{c}})Miehe, Mauthe, \&
  Teichtmeister]{Miehe2015186}
\textsc{Miehe, C.; Mauthe, S.; Teichtmeister, S.} [2015]: \emph{Minimization
  principles for the coupled problem of darcy-biot-type fluid transport in
  porous media linked to phase field modeling of fracture}.
\newblock Journal of the Mechanics and Physics of Solids, 82: 186 -- 217.

\bibitem[\textsc{Mikeli\'c et~al.}(2015{\natexlab{a}})Mikeli\'c, Wheeler, \&
  Wick]{MiWheWi14}
\textsc{Mikeli\'c, A.; Wheeler, M.~F.; Wick, T.} [2015]: \emph{A phase-field
  method for propagating fluid-filled fractures coupled to a surrounding porous
  medium}.
\newblock SIAM Multiscale Model. Simul., 13(1): 367--398.

\bibitem[\textsc{Mikeli\'c et~al.}(2015{\natexlab{b}})Mikeli\'c, Wheeler, \&
  Wick]{MiWheWi15b}
\textsc{Mikeli\'c, A.; Wheeler, M.~F.; Wick, T.} [2015]: \emph{A quasi-static
  phase-field approach to pressurized fractures}.
\newblock Nonlinearity, 28(5): 1371--1399.

\bibitem[\textsc{Mikeli{\'{c}} et~al.}(2019)Mikeli{\'{c}}, Wheeler, \&
  Wick]{MiWheWi18}
\textsc{Mikeli{\'{c}}, A.; Wheeler, M.~F.; Wick, T.} [2019]: \emph{Phase-field
  modeling through iterative splitting of hydraulic fractures in a poroelastic
  medium}.
\newblock GEM - International Journal on Geomathematics, 10(1).
\newblock https://doi.org/10.1007/s13137-019-0113-y.

\bibitem[\textsc{Mota et~al.}(2017)Mota, Tezaur, \& Alleman]{alternatMota17}
\textsc{Mota, A.; Tezaur, I.; Alleman, C.} [2017]: \emph{The schwarz
  alternating method in solid mechanics}.
\newblock Computer Methods in Applied Mechanics and Engineering, 319: 19--51.

\bibitem[\textsc{Na \& Sun}(2018)]{na+sun18}
\textsc{Na, S.; Sun, W.} [2018]: \emph{Computational thermomechanics of
  crystalline rock, part i: A combined multi-phase-field/crystal plasticity
  approach for single crystal simulations}.
\newblock Computer Methods in Applied Mechanics and Engineering, 338: 657 --
  691.

\bibitem[\textsc{Nasseri \& Mohanty}(2008)]{NASSERI08}
\textsc{Nasseri, M.; Mohanty, B.} [2008]: \emph{Fracture toughness anisotropy
  in granitic rocks}.
\newblock International Journal of Rock Mechanics and Mining Sciences, 45(2):
  167 -- 193.

\bibitem[\textsc{Nguyen et~al.}(2017)Nguyen, R\'ethor\'e, \& Baietto]{NGUYEN17}
\textsc{Nguyen, T.~T.; R\'ethor\'e, J.; Baietto, M.-C.} [2017]: \emph{Phase
  field modelling of anisotropic crack propagation}.
\newblock European Journal of Mechanics - A/Solids, 65: 279 -- 288.

\bibitem[\textsc{Noii \& Wick}(2019)]{NoiiWick2019}
\textsc{Noii, N.; Wick, T.} [2019]: \emph{A phase-field description for
  pressurized and non-isothermal propagating fractures}.
\newblock Computer Methods in Applied Mechanics and Engineering, 351: 860 --
  890.

\bibitem[\textsc{Paggi \& Reinoso}(2017)]{paggi17}
\textsc{Paggi, M.; Reinoso, J.} [2017]: \emph{Revisiting the problem of a crack
  impinging on an interface:a modeling framework for the interaction between
  the phase field approach for brittle fracture and the interface cohesive zone
  model}.
\newblock Computer Methods in Applied Mechanics and Engineering, 321: 145 --
  172.

\bibitem[\textsc{Park \& Felippa}(2000)]{Park2000}
\textsc{Park, K.; Felippa, C.} [2000]: \emph{A variational principle for the
  formulation of partitioned structural systems}.
\newblock International Journal for Numerical Methods in Engineering, 47:
  395--418.

\bibitem[\textsc{Park \& Felippa}(2002)]{Park2002}
\textsc{Park, K.; Felippa, C.} [2002]: \emph{A simple algorithm for localized
  construction of non-matching structural interfaces}.
\newblock International Journal for Numerical Methods in Engineering, 53(9):
  2117--2142.

\bibitem[\textsc{Reis \& Pires}(2014)]{REIS2014168}
\textsc{Reis, F.; Pires, F.~A.} [2014]: \emph{A mortar based approach for the
  enforcement of periodic boundary conditions on arbitrarily generated meshes}.
\newblock Computer Methods in Applied Mechanics and Engineering, 274: 168--191.

\bibitem[\textsc{Rice}(1968)]{ricecond}
\textsc{Rice, J.} [1968]: \emph{Mathematical analysis in the mechanics of
  fracture}.
\newblock Academic Press, New York, chapter 3 of fracture: An advanced treatise
  edition, 3: 191--311.

\bibitem[\textsc{Sargado et~al.}(2017)Sargado, Keilegavlen, Berre, \&
  Nordbotten]{sargado17}
\textsc{Sargado, J.; Keilegavlen, E.; Berre, I.; Nordbotten, J.} [2017]:
  \emph{High-accuracy phase-field models for brittle fracture based on a new
  family of degradation functions}.
\newblock Journal of the Mechanics and Physics of Solids.

\bibitem[\textsc{Song et~al.}(2015)Song, Youn, \& Park]{Park2015}
\textsc{Song, Y.; Youn, S.; Park, K.} [2015]: \emph{A gap element for treating
  non-matching discrete interfaces}.
\newblock International Journal for Numerical Methods in Engineering, 56(3):
  551--563.

\bibitem[\textsc{Takei et~al.}(2013)Takei, Roman, Bico, Hamm, \& Melo]{Takei13}
\textsc{Takei, A.; Roman, B.; Bico, J.; Hamm, E.; Melo, F.} [2013]:
  \emph{Forbidden directions for the fracture of thin anisotropic sheets: An
  analogy with the wulff plot}.
\newblock Physical Review Letters, 110: 144301.

\bibitem[\textsc{Teichtmeister et~al.}(2017)Teichtmeister, Kienle, Aldakheel,
  \& Keip]{Fadianiso17}
\textsc{Teichtmeister, S.; Kienle, D.; Aldakheel, F.; Keip, M.-A.} [2017]:
  \emph{Phase field modeling of fracture in anisotropic brittle solids}.
\newblock International Journal of Non-Linear Mechanics, 97: 1--21.

\bibitem[\textsc{Unger et~al.}(2007)Unger, Eckardt, \& K\"onke]{UNGER20074087}
\textsc{Unger, J.~F.; Eckardt, S.; K\"onke, C.} [2007]: \emph{Modelling of
  cohesive crack growth in concrete structures with the extended finite element
  method}.
\newblock Computer Methods in Applied Mechanics and Engineering, 196(41):
  4087--4100.

\bibitem[\textsc{Verhoosel \& De~Borst}(2013)]{verhoosel+deborst13}
\textsc{Verhoosel, C.~V.; De~Borst, R.} [2013]: \emph{A phase-field model for
  cohesive fracture}.
\newblock International Journal for Numerical Methods in Engineering, 96:
  43--62.

\bibitem[\textsc{Wheeler \& Wheeler}(2019)]{IPARS}
\textsc{Wheeler, J.; Wheeler, M.} [2019]: \emph{IPARS, A New Generation
  Framework for Petroleum Reservoir Simulation, Technical Reference}.
\newblock {http://csm.ices.utexas.edu/ipars/}.

\bibitem[\textsc{Wheeler et~al.}(2014)Wheeler, Wick, \&
  Wollner]{WickLagrange2014}
\textsc{Wheeler, M.; Wick, T.; Wollner, W.} [2014]: \emph{An
  augmented-lagrangian method for the phase-field approach for pressurized
  fractures}.
\newblock Computer Methods in Applied Mechanics and Engineering, 271: 69--85.

\bibitem[\textsc{Wick}(2016)]{Wi16_dwr_pff}
\textsc{Wick, T.} [2016]: \emph{Goal functional evaluations for phase-field
  fracture using {PU}-based {DWR} mesh adaptivity}.
\newblock Computational Mechanics, 57(6): 1017--1035.

\bibitem[\textsc{Wick}(2017)]{Wi17_SISC}
\textsc{Wick, T.} [2017]: \emph{An error-oriented {N}ewton/inexact augmented
  {L}agrangian approach for fully monolithic phase-field fracture propagation}.
\newblock SIAM Journal on Scientific Computing, 39(4): B589--B617.

\bibitem[\textsc{Wick et~al.}(2016)Wick, Singh, \& Wheeler]{WiSiWhe15}
\textsc{Wick, T.; Singh, G.; Wheeler, M.} [2016]: \emph{Fluid-filled fracture
  propagation using a phase-field approach and coupling to a reservoir
  simulator}.
\newblock SPE Journal, 21(03): 981--999.

\bibitem[\textsc{Wohlmuth}(2000)]{Wohlmuth}
\textsc{Wohlmuth, B.} [2000]: \emph{A mortar finite element method using dual
  spaces for the lagrange multiplier}.
\newblock SIAM Journal on Numerical Analysis, 38(3): 989--1012.

\bibitem[\textsc{Wohlmuth}(2002)]{Wohlmuth2003}
\textsc{Wohlmuth, B.} [2002]: \emph{A comparison of dual lagrange multiplier
  spaces for mortar finite element discretizations}.
\newblock ESAIM: Mathematical Modelling and Numerical Analysis - Mod\'elisation
  Math\'ematique et Analyse Num\'erique, 36(6): 995--1012.

\bibitem[\textsc{Wriggers}(2008)]{wriggers06}
\textsc{Wriggers, P.} [2008]: \emph{{N}onlinear {F}inite {E}lements}.
\newblock Springer, Berlin, Heidelberg, New York.

\bibitem[\textsc{Wriggers et~al.}(2019)Wriggers, Aldakheel, \&
  Hudobivnik]{Wriggers+aldakheel+blaz19}
\textsc{Wriggers, P.; Aldakheel, F.; Hudobivnik, B.} [2019]: \emph{Application
  of the virtual element method in mechanics}.
\newblock GAMM-Rundbriefe, pp. 4--10.
\newblock ISSN: 2196-3789.

\bibitem[\textsc{Zhang \& Oskay}(2015)]{Zhang}
\textsc{Zhang, S.; Oskay, C.} [2015]: \emph{Variational multiscale enrichment
  method with mixed boundary conditions for elasto-viscoplastic problems}.
\newblock Computational Mechanics, 55(4): 771--787.

\end{thebibliography}

\end{document}